\documentclass{article}
\usepackage{amsmath}
\usepackage{amssymb}
\begin{document}
\def\e#1\e{\begin{equation}#1\end{equation}}
\def\ea#1\ea{\begin{align}#1\end{align}}
\def\eq#1{{\rm(\ref{#1})}}
\newtheorem{thm}{Theorem}[section]
\newtheorem{lem}[thm]{Lemma}
\newtheorem{prop}[thm]{Proposition}
\newtheorem{cor}[thm]{Corollary}
\newenvironment{dfn}{\medskip\refstepcounter{thm}
\noindent{\bf Definition \thesection.\arabic{thm}\ }}{\medskip}
\newenvironment{proof}[1][,]{\medskip\ifcat,#1
\noindent{\it Proof.\ }\else\noindent{\it Proof of #1.\ }\fi}
{\hfill$\square$\medskip}
\def\dim{\mathop{\rm dim}}
\def\Re{\mathop{\rm Re}}
\def\Im{\mathop{\rm Im}}
\def\Image{\mathop{\rm Image}}
\def\ind{\mathop{\rm ind}}
\def\Ker{\mathop{\rm Ker}}
\def\Coker{\mathop{\rm Coker}}
\def\vol{\mathop{\rm vol}}
\def\supp{\mathop{\rm supp}}
\def\id{\mathop{\rm id}}
\def\U{\mathbin{\rm U}}
\def\SU{\mathop{\rm SU}}
\def\g{\mathfrak{g}} 
\def\su{\mathfrak{su}} 
\def\so{\mathfrak{so}} 
\def\ge{\geqslant} 
\def\le{\leqslant} 
\def\R{\mathbin{\mathbb R}}
\def\N{\mathbin{\mathbb N}}
\def\Z{\mathbin{\mathbb Z}}
\def\C{\mathbin{\mathbb C}}
\def\al{\alpha}
\def\be{\beta}
\def\ga{\gamma}
\def\de{\delta}
\def\ep{\epsilon}
\def\la{\lambda}
\def\ka{\kappa}
\def\th{\theta}
\def\ze{\zeta}
\def\up{\upsilon}
\def\vp{\varphi}
\def\si{\sigma}
\def\om{\omega}
\def\De{\Delta}
\def\La{\Lambda}
\def\Om{\Omega}
\def\Ga{\Gamma}
\def\Si{\Sigma}
\def\Th{\Theta}
\def\Up{\Upsilon}
\def\d{{\rm d}}
\def\pd{\partial}
\def\ts{\textstyle}
\def\sst{\scriptscriptstyle}
\def\w{\wedge}
\def\sm{\setminus}
\def\op{\oplus}
\def\ot{\otimes}
\def\iy{\infty}
\def\ra{\rightarrow}
\def\hookra{\hookrightarrow}
\def\longra{\longrightarrow}
\def\t{\times}
\def\na{\nabla}
\def\ha{{\textstyle\frac{1}{2}}}
\def\ti{\tilde}
\def\bs{\boldsymbol}
\def\ov{\overline}
\def\ovB{\,\overline{\!B}}
\def\D{{\cal D}}
\def\sSi{{\smash{\sst\Si}}}
\def\sSii{{\smash{\sst\Si_i}}}
\def\F{{\cal F}}
\def\sF{{\smash{\sst\cal F}}}
\def\sX{{\smash{\sst X}}}
\def\sXp{{\smash{\sst X'}}}
\def\ms#1{\vert#1\vert^2}
\def\bms#1{\bigl\vert#1\bigr\vert^2}
\def\md#1{\vert #1 \vert}
\def\bmd#1{\big\vert #1 \big\vert}
\def\nm#1{\Vert #1 \Vert}
\def\cnm#1#2{\Vert #1 \Vert_{C^{#2}}} 
\def\lnm#1#2{\Vert #1 \Vert_{L^{#2}}} 
\def\snm#1#2#3{\Vert #1 \Vert_{L^{#2}_{#3}}} 
\def\an#1{\langle#1\rangle}
\def\ban#1{\bigl\langle#1\bigr\rangle}
\title{Special Lagrangian submanifolds with isolated \\
conical singularities. I. Regularity}
\author{Dominic Joyce \\ Lincoln College, Oxford}
\date{}
\maketitle

\section{Introduction}
\label{cr1}

{\it Special Lagrangian $m$-folds (SL\/ $m$-folds)} are a
distinguished class of real $m$-dimensional minimal submanifolds
which may be defined in $\C^m$, or in {\it Calabi--Yau $m$-folds},
or more generally in {\it almost Calabi--Yau $m$-folds} (compact
K\"ahler $m$-folds with trivial canonical bundle).

This is the first in a series of five papers
\cite{Joyc3,Joyc4,Joyc5,Joyc6} studying SL $m$-folds with
{\it isolated conical singularities}. That is, we consider an
SL $m$-fold $X$ in $M$ with singularities at $x_1,\ldots,x_n$
in $M$, such that for some SL cones $C_i$ in $T_{\smash{x_i}}M
\cong\C^m$ with $C_i\sm\{0\}$ nonsingular, $X$ approaches $C_i$
near $x_i$ in an asymptotic $C^1$ sense. Readers are advised to
begin with the final paper \cite{Joyc6}, which surveys the
series, and applies the results to prove some conjectures.

Having a good understanding of the singularities of special
Lagrangian submanifolds will be essential in clarifying the
Strominger--Yau--Zaslow conjecture on the Mirror Symmetry
of Calabi--Yau 3-folds \cite{SYZ}, and also in resolving
conjectures made by the author \cite{Joyc1} on defining
new invariants of Calabi--Yau 3-folds by counting special
Lagrangian homology 3-spheres with weights. The series of
papers aims to develop such an understanding for simple
kinds of singularities of SL $m$-folds.

This first paper lays the foundations for
\cite{Joyc3,Joyc4,Joyc5,Joyc6}, setting up definitions and
notation, and proving some auxiliary results in symplectic
geometry and asymptotic analysis that will be needed in
\cite{Joyc3,Joyc4,Joyc5}. However, we also prove results of
independent interest on the {\it regularity} of SL $m$-folds
with conical singularities, and also of {\it Asymptotically
Conical\/} SL $m$-folds in~$\C^m$.

We initially define SL $m$-folds $X$ with conical singularities
$x$ in Definition \ref{cr3def5} below, such that $X$ approaches
the cone $C$ near $x$ like $O(r^{\mu-1})$ in a $C^1$ sense for
some {\it rate} $\mu\in(2,3)$, where $r$ is the distance to $x$
in $M$. In \S\ref{cr5} we use elliptic regularity to prove an
$O(r^{\mu-1-k})$ asymptotic estimate on the $k^{\rm th}$
derivative of the difference between $X$ and $C$ near $x$, for
all~$k\ge 0$.

We also show that the rate $\mu\in(2,3)$ can be improved, up
to a limit depending on the eigenvalues of the Laplacian on
$\Si=C\cap{\cal S}^{2m-1}$. These results in effect {\it
strengthen} the definition of conical singularities of
SL $m$-folds, showing that it is equivalent to a rather
stronger definition.

Section \ref{cr6} relates special Lagrangian geometry to
{\it Geometric Measure Theory}. Our main result here is that
a special Lagrangian integral current whose tangent cones are
`Jacobi integrable' and of multiplicity one is actually an SL
$m$-fold with conical singularities. Thus we {\it weaken}
the definition of conical singularities of SL $m$-folds.

In \cite{Joyc3} we will study the {\it deformation theory}
of compact SL $m$-folds $X$ with conical singularities in
an almost Calabi--Yau $m$-fold $M$. We will show that the
moduli space ${\cal M}_\sX$ of deformations of $X$
as an SL $m$-fold with conical singularities in $M$ is
locally homeomorphic to the zeroes of a smooth map
$\Phi:{\cal I}_\sXp\ra{\cal O}_\sXp$ between
finite-dimensional vector spaces, and if the {\it obstruction
space} ${\cal O}_\sXp$ is zero then ${\cal M}_\sX$
is a smooth manifold.

Then \cite{Joyc4,Joyc5} will consider {\it desingularizations}
of a compact SL $m$-fold $X$ with conical singularities
$x_1,\ldots,x_n$ with cones $C_1,\ldots,C_n$ in an
almost Calabi--Yau $m$-fold $M$. We will take nonsingular
Asymptotically Conical SL $m$-folds $L_1,\ldots,L_n$ in
$\C^m$ asymptotic to $C_1,\ldots,C_n$ at infinity, and
glue them in to $X$ at $x_1,\ldots,x_n$ to get a smooth
family of compact, {\it nonsingular} SL $m$-folds $\ti N$
in $M$ which converge to~$X$.

We begin in \S\ref{cr2} by defining {\it Riemannian manifolds with
conical singularities}, and developing the theory of {\it weighted
Sobolev spaces} upon them, and the {\it Fredholm properties} of the
Laplacian on these spaces, adapting results of Lockhart and McOwen
\cite{Lock,LoMc}. We give a detailed treatment, in the hope that
\S\ref{cr2} will be a useful reference for further work on
manifolds with conical singularities.

Almost Calabi--Yau manifolds and special Lagrangian geometry
are introduced in \S\ref{cr3}, and SL $m$-folds with conical
singularities defined in \S\ref{cr33}. Then \S\ref{cr4} proves
{\it Lagrangian Neighbourhood Theorems} for SL $m$-folds $X$
with conical singularities in almost Calabi--Yau $m$-folds $M$.
Essentially these are special coordinate systems on $M$ near
$X$, in which the symplectic form $\om$ on $M$ has a canonical
form, and which satisfy asymptotic conditions near the singular
points $x_1,\ldots,x_n$ of $X$. These theorems will be important
tools in~\cite{Joyc3,Joyc4,Joyc5}.

In \S\ref{cr5} we prove regularity results for the convergence
of $X$ to its cone $C_i$ near a singular point $x_i$, with all
derivatives. Section \ref{cr6} introduces Geometric Measure
Theory, recalls results on {\it tangent cones} due to Adams
and Simon, and shows that under some conditions on its tangent
cones, a special Lagrangian integral current is an SL $m$-fold
with conical singularities, in the sense of~\S\ref{cr33}.

We finish in \S\ref{cr7} by extending many of the results
of \S\ref{cr4}--\S\ref{cr5} to {\it Asymptotically Conical\/}
SL $m$-folds in $\C^m$, which are asymptotic to an SL cone
$C$ in $\C^m$ at some rate $\la$. These results will be
needed in~\cite{Joyc4,Joyc5}.

Throughout we shall for simplicity take all submanifolds to be
{\it embedded}. Nearly all of our results generalize immediately
to {\it immersed\/} submanifolds, with only cosmetic changes.
However, this does {\it not\/} apply to the Geometric Measure
Theory material in \S\ref{cr6}, where the tangent cones really
do have to be embedded rather than immersed.
\medskip

\noindent{\it Acknowledgements.} I would like to thank Stephen
Marshall for many discussions on the material of \S\ref{cr2}
and his thesis \cite{Mars}, Mark Haskins for essential help
with \S\ref{cr6}, and Tadashi Tokieda and Ivan Smith for useful
conversations. I was supported by an EPSRC Advanced Research
Fellowship whilst writing this paper.

\section{Manifolds with conical singularities}
\label{cr2}

We shall study a class of singular Riemannian manifolds
with isolated singularities modelled on cones.

\begin{dfn} Let $(X,d)$ be a metric space and $x_1,\ldots,x_n$
be distinct points in $X$, and define $X'=X\sm\{x_1,\ldots,x_n\}$.
We call $X$ a {\it Riemannian $m$-manifold with conical
singularities\/} $x_1,\ldots,x_n$ if the following conditions hold:
\begin{itemize}
\item[(a)] $X'$ has the structure of a smooth, connected $m$-manifold
with a Riemannian metric $g$ inducing the metric $d$ on~$X'$.
\item[(b)] We are given $\ep\in(0,1)$ small such that
$d(x_i,x_j)>2\ep$ for $1\le i<j\le n$ and a compact, nonsingular
Riemannian $(m\!-\!1)$-manifold $(\Si_i,g_\sSii)$ for
$i=1,\ldots,n$. Write points in $\Si_i\t(0,\ep)$ as $(\si,r)$.
Define the {\it cone metric} $h_i$ on $\Si_i\t(0,\ep)$ to
be~$h_i=r^2g_\sSii+\d r^2$.
\item[(c)] For $i=1,\ldots,n$ there exist $\nu_i>0$ and a
diffeomorphism $\phi_i:\Si_i\t(0,\ep)\ra S_i=\{y\in X:
0<d(x_i,y)<\ep\}\subset X'$ such that
\e
\bmd{\na^k\bigl(\phi_i^*(g)-h_i\bigr)}=O(r^{\nu_i-k})
\quad\text{as $r\ra 0$, for all $k\ge 0$.}
\label{cr2eq1}
\e
Here the Levi-Civita connection $\na$ and $\md{\,.\,}$ are
computed using~$h_i$.
\end{itemize}
Let $C_\sSii$ be the {\it Riemannian cone} on
$(\Si_i,g_\sSii)$, to be defined in Definition
\ref{cr2def2}. We call $C_\sSii$ the {\it cone}
and $\nu_i$ the {\it rate} of the singular point~$x_i$.
\label{cr2def1}
\end{dfn}

Usually we will also assume that $X$ is {\it compact}. Equation
\eq{cr2eq1} implies that near $x_i$ the metric $g$ and its
derivatives are asymptotic to the cone metric $h_i$ on
$\Si_i\t(0,\ep)$. For applications it generally suffices for
\eq{cr2eq1} to hold when $k\le l$ for some $l$. However, we
will show in \S\ref{cr5} that for the singular SL $m$-folds
we are interested in \eq{cr2eq1} holds for all $k\ge 0$
automatically, so we may as well assume it.

Various authors have studied analysis of elliptic operators on
classes of spaces including manifolds with conical singularities.
We shall quote parts of their work, adapting it for our purposes
where necessary. We treat the subject at some length in the hope
that this will be a useful reference for future work on manifolds
with conical singularities.

We start in \S\ref{cr21} by discussing Riemannian cones and harmonic
functions on them. Section \ref{cr22} defines Banach spaces of
functions on $X'$ using weights, and \S\ref{cr23} gives elliptic
regularity results for the Laplacian on these spaces. Finally,
\S\ref{cr24} and \S\ref{cr25} discuss homology, cohomology and
Hodge theory on $X'$ and~$X$.

\subsection{Riemannian cones and harmonic functions}
\label{cr21}

{\it Riemannian cones} are a class of singular Riemannian
manifolds.

\begin{dfn} Let $(\Si,g_\sSi)$ be a compact Riemannian
$(m\!-\!1)$-manifold, not necessarily connected. Define the
{\it cone $C_\sSi$ on} $\Si$ to be $\{0\}\cup C_\sSi'$
where $C'_\sSi=\Si\t(0,\iy)$. Write points in $C'_\sSi$
as $(\si,r)$. Define a Riemannian metric, the {\it cone metric}
$g$ on $C'_\sSi$ by $g=\d r^2+r^2g_\sSi$.

Define a metric $d$ on $C_\sSi$ to be that induced by $g$
on the connected components of $C'_\sSi$, together with
$d\bigl(0,(\si,r)\bigr)=r$ for $(\si,r)\in C'_\sSi$ and
$d\bigl((\si,r),(\si',r')\bigr)=r+r'$ for $\si,\si'$ in
different connected components of $\Si$ and $r,r'>0$. Then
$(C_\sSi,d)$ is a metric space, called the {\it Riemannian
cone on} $\Si$. It is a singular Riemannian manifold, with an
isolated singularity at the {\it vertex} 0. Often we will take
$d$ as given and refer to $C_\sSi$ as a Riemannian cone.

For $t>0$, define the {\it dilation} $t:C_\sSi\ra C_\sSi$
by $t0=0$ and $t(\si,r)=(\si,tr)$. Then $t^*(d)=td$ and $t^*(g)=t^2g$.
For $\al\in\R$, we say that a function $u:C'_\sSi\ra\R$ is
{\it homogeneous of order} $\al$ if $u\circ t\equiv t^\al u$ for all
$t>0$. Equivalently, $u$ is homogeneous of order $\al$ if $u(\si,r)
\equiv r^\al v(\si)$ for some function~$v:\Si\ra\R$.
\label{cr2def2}
\end{dfn}

Clearly, a Riemannian cone $(C_\sSi,d)$ is an example
of a manifold with conical singularities. Here is an elementary
lemma on {\it harmonic functions} on cones.

\begin{lem} In the situation of Definition \ref{cr2def2}, let\/
$u(\si,r)\equiv r^\al v(\si)$ be a homogeneous function of order
$\al$ on $C'_\sSi=\Si\t(0,\iy)$, for $v\in C^2(\Si)$. Then
\e
\De u(\si,r)=r^{\al-2}\bigl(\De_\sSi v-\al(\al+m-2)v\bigr),
\label{cr2eq2}
\e
where $\De$, $\De_\sSi$ are the Laplacians on $(C'_\sSi,g)$
and\/ $(\Si,g_\sSi)$. Hence, $u$ is harmonic on $C'_\sSi$
if and only if\/ $v$ is an eigenfunction of\/ $\De_\sSi$ with
eigenvalue~$\al(\al+m-2)$.
\label{cr2lem1}
\end{lem}

Now for our later work we should also consider harmonic functions
$u$ on cones with more general scaling behaviour under dilations
than homogeneous of order $\al$. For example, $\R^2$ with its
Euclidean metric is the Riemannian cone on ${\cal S}^1$, and
$\log r$ is harmonic on $\R^2\sm\{0\}$. We shall show that
in dimension $m>2$, harmonic functions cannot scale like
$r^\al(\log r)^k$ for~$k>0$.

\begin{prop} In the situation of Definition \ref{cr2def2},
suppose $m>2$. Then there do not exist any harmonic functions
$u$ on $C'_\sSi=\Si\t(0,\iy)$ of the form
\e
u(\si,r)=r^\al(\log r)^kv_k(\si)+r^\al(\log r)^{k-1}v_{k-1}(\si)
+\cdots+r^\al v_0(\si),
\label{cr2eq3}
\e
where $\al\in\R$, $k>0$ and\/ $v_k,v_{k-1},\ldots,v_0\in C^2(\Si)$
with\/~$v_k\ne 0$.
\label{cr2prop1}
\end{prop}

\begin{proof} Suppose that $u$ in \eq{cr2eq3} is harmonic. By
applying infinitesimal dilations we see that $r\frac{\pd u}{\pd r}$
is also harmonic, and so
\begin{equation*}
r\frac{\pd u}{\pd r}-\al u=
r^\al(\log r)^{k-1}kv_k(\si)+r^\al(\log r)^{k-2}(k-1)v_{k-1}(\si)
+\cdots+r^\al v_1(\si)
\end{equation*}
is harmonic. So if there exist harmonic $u$ of the form \eq{cr2eq3}
for $k$, there also exist such $u$ for $k-1$. Thus by induction,
it is sufficient to prove the case~$k=1$.

Suppose for a contradiction that $u$ is harmonic of the form
\eq{cr2eq3} with $k=1$ and $v_1\ne 0$ in $C^2(\Si)$. A calculation
similar to Lemma \ref{cr2lem1} shows that
\begin{align*}
\De u(\si,r)&=r^{\al-2}\log r\bigl(\De_\sSi v_1-\al(\al+m-2)v_1
\bigr)\\
&\qquad +r^{\al-2}\bigl(\De_\sSi v_0-\al(\al+m-2)v_0
-(2\al+m-2)v_1\bigr).
\end{align*}
Thus, as $u$ is harmonic we have
\begin{equation*}
\De_\sSi v_1=\al(\al+m-2)v_1
\quad\text{and}\quad
\De_\sSi v_0=\al(\al+m-2)v_0+(2\al+m-2)v_1.
\end{equation*}
Integrating $v_0\bigl(\De_\sSi-\al(\al+m-2)\bigr)v_1$ over
$\Si$ by parts we get
\begin{align*}
0&=\int_\Si v_0\bigl(\De_\sSi-\al(\al+m-2)\bigr)v_1\,\d V_g
=\int_\Si v_1\bigl(\De_\sSi-\al(\al+m-2)\bigr)v_0\,\d V_g\\
&=\int_\Si v_1(2\al+m-2)v_1=(2\al+m-2)\lnm{v_1}{2}^2.
\end{align*}
Thus $2\al+m-2=0$, as $v_1\ne 0$, so $\al=\ha(2-m)$.
But then $\De_\sSi v_1=-\frac{1}{4}(m-2)^2v_1$
and $v_1\ne 0$, so that $-\frac{1}{4}(m-2)^2$ is an
eigenvalue of $\De_\sSi$. As $m>2$ this contradicts the
fact that eigenvalues of $\De_\sSi$ are nonnegative.
\end{proof}

Here is some more notation.

\begin{dfn} In the situation of Definition \ref{cr2def2},
suppose $m>2$ and define
\e
\D_\sSi=\bigl\{\al\in\R:\text{$\al(\al+m-2)$ is
an eigenvalue of $\De_\sSi$}\bigr\}.
\label{cr2eq4}
\e
By Lemma \ref{cr2lem1}, an equivalent definition is that
$\D_\sSi$ is the set of $\al\in\R$ for which there
exists a nonzero homogeneous harmonic function $u$ of order
$\al$ on $C'_\sSi$. By properties of the spectrum of
$\De_\sSi$, it follows that $\D_\sSi$ is a
countable, discrete subset of~$\R$.

Define $m_\sSi:\D_\sSi\ra\N$ by taking
$m_\sSi(\al)$ to be the multiplicity of the eigenvalue
$\al(\al+m-2)$ of $\De_\sSi$, or equivalently the
dimension of the vector space of homogeneous harmonic
functions $u$ of order $\al$ on $C'_\sSi$. Define
$N_\sSi:\R\ra\Z$ by
\e
N_\sSi(\de)=
-\sum_{\!\!\!\!\al\in\D_\sSi\cap(\de,0)\!\!\!\!}m_\sSi(\al)
\;\>\text{if $\de<0$, and}\;\>
N_\sSi(\de)=
\sum_{\!\!\!\!\al\in\D_\sSi\cap[0,\de]\!\!\!\!}m_\sSi(\al)
\;\>\text{if $\de\ge 0$.}
\label{cr2eq5}
\e
Then $N_\sSi$ is monotone increasing and upper semicontinuous,
and is discontinuous exactly on $\D_\sSi$, increasing by
$m_\sSi(\al)$ at each $\al\in\D_\sSi$. As the
eigenvalues of $\De_\sSi$ are nonnegative, we see that
$\D_\sSi\cap(2-m,0)=\emptyset$ and $N_\sSi\equiv 0$
on~$(2-m,0)$.
\label{cr2def3}
\end{dfn}

\subsection{Weighted Banach spaces}
\label{cr22}

We will need the following tool, a smoothed out version of the
distance from the singular set $\{x_1,\ldots,x_n\}$ in~$X$.

\begin{dfn} Let $(X,d)$ be a compact Riemannian manifold with
conical singularities $\{x_1,\ldots,x_n\}$, and use the notation
of Definition \ref{cr2def1}. Define a {\it radius function}
$\rho$ on $X'$ to be a smooth function $\rho:X'\ra(0,1]$ such
that $\rho(y)=d(x_i,y)$ whenever $0<d(x_i,y)\le\ha\ep$ and
$i=1,\ldots,n$, and $\rho(y)=1$ when $d(x_i,y)\ge\ep$ for all
$i=1,\ldots,n$. Radius functions always exist.

For $\bs{\be}=(\be_1,\ldots,\be_n)\in\R^n$, define a function
$\rho^{\bs\be}$ on $X'$ by $\rho^{\bs\be}(y)=\rho(y)^{\be_i}$
whenever $0<d(x_i,y)<\ep$ for some $i=1,\ldots,n$ and
$\rho^{\bs\be}(y)=1$ when $d(x_i,y)\ge\ep$ for all $i=1,\ldots,n$.
Then $\rho^{\bs\be}$ is well-defined and smooth on $X'$, and
equals $\rho^{\be_i}$ near $x_i$ in $X'$. If ${\bs\be},{\bs\ga}
\in\R^n$, write ${\bs\be}\ge{\bs\ga}$ if $\be_i\ge\ga_i$ and
${\bs\be}>{\bs\ga}$ if $\be_i>\ga_i$ for $i=1,\ldots,n$. If
${\bs\be}\in\R^n$ and $a\in\R$, write ${\bs\be}+a=(\be_1+a,
\ldots,\be_n+a)$ in~$\R^n$.
\label{cr2def4}
\end{dfn}

Now we define some Banach spaces of functions on~$X'$.

\begin{dfn} Let $(X,d)$ be a compact Riemannian $m$-manifold
with conical singularities $x_1,\ldots,x_n$, and use the notation
of Definition \ref{cr2def1}. Let $\rho$ be a radius function
on $X'$. For ${\bs\be}\in\R^n$ and $k\ge 0$ define
$C^k_{\smash{\bs\be}}(X')$ to be the space of continuous functions
$f$ on $X'$ with $k$ continuous derivatives, such that
$\bmd{\rho^{-{\bs\be}+j}\na^jf}$ is bounded on $X'$ for
$j=0,\dots,k$. Define the norm $\nm{\,.\,}_{\smash{C^k_{\bs\be}}}$
on $C^k_{\smash{\bs\be}}(X')$ by
\e
\nm{f}_{\smash{C^k_{\bs\be}}}=\sum_{j=0}^k\sup_{X'}
\bmd{\rho^{-{\bs\be}+j}\na^jf}.
\label{cr2eq6}
\e
Then $C^k_{\smash{\bs\be}}(X')$ is a Banach space.
Define~$C^\iy_{\smash{\bs\be}}(X')=\bigcap_{k\ge 0}
C^k_{\smash{\bs\be}}(X')$.

For $p\ge 1$, ${\bs\be}\in\R^n$ and $k\ge 0$ define the
{\it weighted Sobolev space} $L^p_{k,{\bs\be}}(X')$ to be
the set of functions $f$ on $X'$ that are locally integrable
and $k$ times weakly differentiable, and for which the norm
\e
\snm{f}p{k,{\bs\be}}=
\left(\sum_{j=0}^k\int_{X'}\bmd{\rho^{-{\bs\be}+j}\na^jf}^p
\rho^{-m}\d V_g\right)^{1/p}
\label{cr2eq7}
\e
is finite. Then $L^p_{k,{\bs\be}}(X')$ is a Banach space, and
$L^2_{k,{\bs\be}}(X')$ a Hilbert space. 
\label{cr2def5}
\end{dfn}

We call these {\it weighted Banach spaces} since the norms are
locally weighted by a power of $\rho$. Roughly speaking, if $f$
lies in $L^p_{k,{\bs\be}}(X')$ or $C^k_{\smash{\bs\be}}(X')$ then
$f$ grows at most like $\rho^{\be_i}$ near $x_i$ as $\rho\ra 0$,
and so the multi-index ${\bs\be}=(\be_1,\ldots,\be_n)$ should be
interpreted as an {\it order of growth}. Similarly, $\na^jf$ grows
at most like $\rho^{\be_i-j}$ near $x_i$ for $j=1,\dots,k$. The
vector spaces $L^p_{k,{\bs\be}}(X')$ and $C^k_{\smash{\bs\be}}(X')$
are independent of the choice of radius function $\rho$. Different
choices of $\rho$ give equivalent norms.

Our spaces $L^p_{k,{\bs\be}}(X')$ are part of the scheme of Lockhart
and McOwen \cite{Lock}, \cite{LoMc}. They consider a larger class
of metrics, called {\it admissible metrics} on manifolds with ends
\cite[\S 2]{Lock}, and they use two weight functions $z,\rho$
rather than one. In the notation of \cite[\S 4]{Lock}, our space
$L^p_{k,{\bs\be}}(X')$ coincides with Lockhart's space $W^p_{k,{\bs\de},a}
(X')$ if ${\bs\be}=a-{\bs\de}$. Definition \ref{cr2def5} is actually
based on Bartnik \cite[\S 1]{Bart} for asymptotically Euclidean manifolds.

The Banach space dual of $L^p_{0,\bs\be}(X')$ is another space of
the same form.

\begin{lem} In the situation above, let\/ $p,q>1$ with\/
$\frac{1}{p}+\frac{1}{q}=1$ and\/ ${\bs\be}\in\R^n$. Then the
map $\an{\,,\,}:L^p_{0,\bs\be}(X')\t L^q_{0,-{\bs\be}-m}(X')\ra\R$
given by $\an{u,v}=\int_{X'}uv\,\d V_g$ is well-defined
and continuous and defines a dual pairing, so that\/
$L^p_{0,\bs\be}(X')$, $L^q_{0,-{\bs\be}-m}(X')$ are the
Banach space duals of each other.
\label{cr2lem2}
\end{lem}

\begin{proof} Let $u\in L^p_{0,\bs\be}(X')$ and $v\in
L^q_{0,-{\bs\be}-m}(X')$. Then \eq{cr2eq7} and
$\frac{1}{p}+\frac{1}{q}=1$ imply that
$u\rho^{-\bs{\be}-m/p}\in L^p(X')$ and
$v\rho^{\bs{\be}+m/p}\in L^q(X')$ with
\begin{equation*}
\lnm{u\rho^{-\bs{\be}-m/p}}{p}=\snm{u}{p}{0,\bs\be}
\quad\text{and}\quad
\lnm{v\rho^{\bs{\be}+m/p}}{q}=\snm{u}{q}{0,-{\bs\be}-m}.
\end{equation*}
Here $L^p(X'),L^q(X')$ are the usual, unweighted Lebesgue
spaces on $X'$. Since $\frac{1}{p}+\frac{1}{q}=1$, H\"older's
inequality gives $uv\in L^1(X')$ with
\begin{equation*}
\lnm{uv}{1}\le\lnm{u\rho^{-\bs{\be}-m/p}}{p}\cdot
\lnm{v\rho^{\bs{\be}+m/p}}{q}=
\snm{u}{p}{0,\bs\be}\cdot\snm{v}{q}{0,-{\bs\be}-m}.
\end{equation*}
So $\an{u,v}=\int_{X'}uv\,\d V_g$ exists, and $\bmd{\an{u,v}}\le
\lnm{uv}{1}\le\snm{u}{p}{0,\bs\be}\cdot\snm{v}{q}{0,-{\bs\be}-m}$.
Thus $\an{\,,\,}$ is well-defined and continuous. The last part
follows from the well-known fact that $L^p(X')\t L^q(X')\ra\R$
is a dual pairing.
\end{proof}

Here is a weighted version of the {\it Sobolev Embedding Theorem}
and the {\it Kondrakov Theorem}, giving (compact) inclusions
between these spaces.

\begin{thm} In the situation above, suppose $k\ge l\ge 0$ are
integers, $p,q>1$ and\/ ${\bs\be},{\bs\ga}\in\R^n$. Then
\begin{itemize}
\item[{\rm(a)}] If\/ $\frac{1}{p}\le\frac{1}{q}+\frac{k-l}{m}$
and\/ ${\bs\be}\ge{\bs\ga}$, then $L^p_{k,{\bs\be}}(X')\hookra
L^q_{l,{\bs\ga}}(X')$ is a continuous inclusion. If\/ $\frac{1}{p}<
\frac{1}{q}+\frac{k-l}{m}$ and\/ ${\bs\be}>{\bs\ga}$, this
inclusion is compact.
\item[{\rm(b)}] If\/ $\frac{1}{p}<\frac{k-l}{m}$ and\/ ${\bs\be}\ge
{\bs\ga}$, then $L^p_{k,{\bs\be}}(X')\hookra C^l_{\bs\ga}(X')$ is a
continuous inclusion. If\/ $\frac{1}{p}<\frac{k-l}{m}$ and\/
${\bs\be}>{\bs\ga}$, this inclusion is compact.
\end{itemize}
\label{cr2thm1}
\end{thm}

\begin{proof} Part (a) follows from \cite[Th.~4.8 \& Th.~4.9]{Lock}
once the notation is disentangled. Inclusion in (b) when ${\bs\be}=
{\bs\ga}$ is proved by Bartnik \cite[Th.~1.2, eq.~(1.9)]{Bart} in
$\R^n$ using a scaling argument on annuli, and then generalized to
asymptotically Euclidean manifolds. The same method works in our
case, where instead of the annulus $A_R=B_{2R}\sm\ovB_R$ in $\R^n$
we substitute $\Si_i\t(R,2R)$ in $C_\sSii$. The rest of (b) follows
from~(a).
\end{proof}

\subsection{Elliptic regularity on weighted spaces}
\label{cr23}

Let $(X,d)$ be a compact Riemannian manifold with conical
singularities, and use the notation above. Let $\De=\d^*\d$
be the Laplacian on functions. We will study the map
\e
\De^p_{k,{\bs\be}}=\De:L^p_{k,{\bs\be}}(X')\ra
L^p_{k-2,{\bs\be}-2}(X')
\label{cr2eq8}
\e
for $p>1$, $k\ge 2$ and ${\bs\be}\in\R^n$. As a shorthand we
will refer to this map as $\De^p_{k,{\bs\be}}$. We will show
that under certain conditions on $\bs\be$ it is {\it Fredholm},
and describe its kernel and cokernel.

Here is an {\it elliptic regularity result\/} for~$\De^p_{k,{\bs\be}}$.

\begin{thm} Let\/ $(X,d)$ be a compact Riemannian manifold with
conical singularities. Then for all\/ $p>1$, $k\ge 2$ and\/
${\bs\be}\in\R^n$ there exists $C>0$ such that if\/
$u\in L^p_{0,{\bs\be}}(X')$ lies in $L^p_2$ locally and\/
$v\in L^p_{k-2,{\bs\be}-2}(X')$ with\/ $\De u=v$ then $u\in
L^p_{k,{\bs\be}}(X')$ and\/~$\snm{u}{p}{k,{\bs\be}}\le
C\bigl(\snm{u}{p}{0,{\bs\be}}+\snm{v}{p}{k-2,{\bs\be}-2}\bigr)$.
\label{cr2thm2}
\end{thm}

\begin{proof} Gilbarg and Trudinger \cite[Th.~9.19]{GiTr} show
that $u$ lies in $L^p_k$ locally in $X'$, and the result then
follows from Lockhart~\cite[Th.~3.7]{Lock}.
\end{proof}

Recall that a continuous linear map between Banach spaces is
{\it Fredholm} if it has finite-dimensional kernel and cokernel.

\begin{thm} Let\/ $(X,d)$ be a compact Riemannian $m$-manifold
with conical singularities $x_1,\ldots,x_n$. Then for all\/
$p>1$, $k\ge 2$ and\/ ${\bs\be}\in\R^n$, the map
$\De^p_{k,{\bs\be}}$ is Fredholm if and only if\/
$\be_i\notin\D_\sSii$ for all\/ $i=1,\ldots,n$, where
$\D_\sSii$ is defined in \eq{cr2eq4}, that is, if\/
$\bs\be$ lies in the subset
\e
\bigl(\R\sm\D_{\sst\Si_1}\bigr)\t\cdots\t
\bigl(\R\sm\D_{\sst\Si_n}\bigr)
\label{cr2eq9}
\e
in $\R^n$. Equivalently, $\De^p_{k,{\bs\be}}$ is Fredholm if and
only if for all\/ $i=1,\ldots,n$, there exists no nonzero homogeneous
harmonic function $u$ on $C_\sSii'$ with rate~$\be_i$.
\label{cr2thm3}
\end{thm}

\begin{proof} Translating our problem into his notation, Lockhart
\cite[Th.~5.2]{Lock} shows that $\De^p_{k,{\bs\be}}$ is Fredholm if
and only if $\be_i\notin\D_\sSii$ for $i=1,\ldots,n$,
where $\D_\sSii$ is a countable, discrete subset of
$\R$. Following the definition of $\D_\sSii$ back through
\cite[Th.~3.7]{Lock}, \cite[Th.~6.2]{LoMc} and \cite[p.~416-7]{LoMc}
we eventually find that it is given by \eq{cr2eq4}. In fact,
\cite[p.~416-7]{LoMc} defines $\D_\sSii$ as the
imaginary part of the spectrum of a complex eigenvalue problem,
but as the spectrum of $\De_\sSii$ is real and nonnegative,
it reduces to \eq{cr2eq4}. The final part follows from
Definition~\ref{cr2def3}.
\end{proof}

We study the dependence of the kernel of $\De^p_{k,{\bs\be}}$
on $p,k$ and~$\bs\be$.

\begin{thm} Let\/ $(X,d)$ be a compact Riemannian $m$-manifold
with conical singularities $x_1,\ldots,x_n$, and let\/ $p>1$,
$k\ge 2$ and\/ ${\bs\be}\in\R^n$. Then $\Ker\bigl(\De^p_{k,
{\bs\be}}\bigr)$ is independent of\/ $k\ge 2$, and is a
finite-dimensional subspace of\/ $C^\iy_{\bs\be}(X')$. If
also $\bs\be$ lies in \eq{cr2eq9} then $\Ker\bigl(\De^p_{k,
{\bs\be}}\bigr)$ is independent of\/ $p>1$, and depends only
on $(X,d)$ and the connected component of\/ \eq{cr2eq9}
containing~$\bs\be$.
\label{cr2thm4}
\end{thm}

\begin{proof} If $u\in L^p_{2,{\bs\be}}(X')$ and $\De u=0$
then Theorem \ref{cr2thm2} with $v=0$ shows that $u\in L^p_{k,
{\bs\be}}(X')$ for any $k\ge 2$. Part (b) of Theorem
\ref{cr2thm1} then implies that $u\in C^l_{{\bs\be}}(X')$
for all $l\ge 0$, and so $u\in C^\iy_{\bs\be}(X')$. Thus the
kernel of $\De:L^p_{k,{\bs\be}}(X')\ra L^p_{k-2,{\bs\be}-2}(X')$
is independent of $k$ and lies in $C^\iy_{\bs\be}(X')$.
Finite-dimensionality follows from~\cite[Cor.~5.6]{Lock}.

When $n=1$ and $\bs\be$ lies in \eq{cr2eq9}, Lockhart and
McOwen show \cite[Lem.~7.3]{LoMc} that the kernel of
$\De^p_{k,{\bs\be}}$ is independent of $p>1$, and depends only
on the connected component of $\bs\be$ in $\R\sm\D_{\sst\Si_1}$,
\cite[Lem.~7.1]{LoMc}. These are easily generalized to the case
$n>1$ as in~\cite[\S 8]{LoMc}.
\end{proof}

Here is an integration by parts formula in weighted Sobolev spaces.

\begin{lem} In the situation above, let\/ $p,q>1$ with\/
$\frac{1}{p}+\frac{1}{q}=1$ and\/ ${\bs\be}\in\R^n$. Then for all\/
$u\in L^p_{2,\bs\be}(X')$ and\/ $v\in L^q_{2,-{\bs\be}+2-m}(X')$
we have
\e
\an{\De u,v}=\int_{X'}(\De u)v\,\d V_g
=\int_{X'}(\d u\cdot\d v)\,\d V_g
=\int_{X'}u(\De v)\d V_g=\an{u,\De v}.
\label{cr2eq10}
\e
\label{cr2lem3}
\end{lem}

\begin{proof} First suppose $u,v$ are smooth with compact support
in $X'$. Then \eq{cr2eq10} is immediate by integration by parts.
But the smooth functions with compact support are dense in
$L^p_{2,{\bs\be}}(X')$ and $L^q_{2,-{\bs\be}+2-m}(X')$ by
\cite[Cor.~4.5]{Lock}, and Lemma \ref{cr2lem2} shows that
$(u,v)\mapsto\an{\De u,v}$, $(u,v)\mapsto\int_{X'}(\d u\cdot\d v)
\,\d V_g$ and $(u,v)\mapsto\an{u,\De v}$ are continuous maps
$L^p_{2,{\bs\be}}(X')\t L^q_{2,-{\bs\be}+2-m}(X')\ra\R$, so
the result follows.
\end{proof}

Now we can describe the cokernel of $\De^p_{k,{\bs\be}}$
when it is Fredholm.

\begin{thm} In the situation of Theorem \ref{cr2thm3}, suppose
$\De^p_{k,{\bs\be}}$ is Fredholm. Define $q>1$ by $\frac{1}{p}+
\frac{1}{q}=1$. Then $u\in L^p_{k-2,{\bs\be}-2}(X')$ lies in the
image of\/ $\De^p_{k,{\bs\be}}$ if and only if\/ $\an{u,v}=0$ for
all\/ $v$ in the kernel of\/ $\De^q_{2,-{\bs\be}+2-m}$. Hence,
the cokernel of\/ $\De^p_{k,{\bs\be}}$ is isomorphic to the dual
of the kernel of\/~$\De^q_{2,-{\bs\be}+2-m}$.
\label{cr2thm5}
\end{thm}

\begin{proof} If $u=\De w$ for $w\in L^p_{k,{\bs\be}}(X')$
and $v\in\Ker(\De^q_{2,-{\bs\be}+2-m})$ then
$\an{u,v}=\an{\De w,v}=\an{w,\De v}=0$ by Lemma
\ref{cr2lem3}. This proves the `only if' part. To prove the `if'
part, first suppose $k=2$, and let $f:L^p_{0,{\bs\be}-2}(X')\ra\R$ be
a linear map vanishing on the image of~$\De^p_{2,{\bs\be}}$.

As $\De^p_{2,{\bs\be}}$ is Fredholm this image is closed
and has finite codimension, so that $f$ is {\it continuous}.
Thus $f$ defines an element of the Banach space dual of
$L^p_{0,{\bs\be}-2}(X')$. Lemma \ref{cr2lem2} then gives
a unique $v\in L^q_{0,-{\bs\be}+2-m}(X')$ such that
$\int_{X'}uv\,\d V_g=\an{u,v}=f(u)$ for all $u\in
L^p_{0,{\bs\be}-2}(X')$. As $f=0$ on the image of
$\De^p_{2,{\bs\be}}$, this shows that
\e
\int_{X'}v(\De w)\,\d V_g=0
\quad\text{for all $w\in L^p_{2,{\bs\be}}(X')$.}
\label{cr2eq11}
\e

By an elliptic regularity result of Morrey \cite[Th.~6.4.3,
p.~246]{Morr}, if $v$ is $L^q$ locally for $q>1$ and \eq{cr2eq11}
holds for all compactly-supported smooth $w$ (which automatically
lie in $L^p_{2,{\bs\be}}(X')$), then $v$ is $L^q_l$ locally for
all $l\ge 0$. Thus $v$ is smooth, by the Sobolev Embedding Theorem,
and integration by parts shows that $v$ is harmonic. Theorem
\ref{cr2thm2} then proves that $v\in L^q_{2,-{\bs\be}+2-m}(X')$.
So $v$ lies in the kernel of~$\De^q_{2,-{\bs\be}+2-m}$.

We have shown that any linear map $f:L^p_{0,{\bs\be}-2}(X')\ra
\R$ vanishing on the image of $\De^p_{2,{\bs\be}}$ is of
the form $f(u)=\an{u,v}$ for some $v$ in the kernel of \eq{cr2eq11}.
As the image of $\De^p_{2,{\bs\be}}$ has finite codimension,
this proves the `if' part when~$k=2$.

For the $k>2$ case, suppose that $u\in L^p_{k-2,{\bs\be}-2}(X')$
and $\an{u,v}=0$ for all\/ $v$ in the kernel of $\De^q_{2,
-{\bs\be}+2-m}$. Then $u\in L^p_{0,{\bs\be}-2}(X')$, so by the $k=2$
case we have $u=\De w$ for some $w\in L^p_{2,{\bs\be}}(X')$.
Theorem \ref{cr2thm2} then implies that $w\in L^p_{k,{\bs\be}}(X')$,
so $u$ lies in the image of $\De^p_{k,{\bs\be}}$, and the
proof is complete.
\end{proof}

The {\it index} of a Fredholm operator $P$ is $\ind(P)=\dim\Ker(P)
-\dim\Coker(P)$.

\begin{thm} Let\/ $(X,d)$ be a compact Riemannian $m$-manifold
for $m>2$ with conical singularities $x_1,\ldots,x_n$, and let\/
$p>1$, $k\ge 2$ and\/ ${\bs\be}$ lie in \eq{cr2eq9}, so that\/
$\De^p_{k,{\bs\be}}$ is Fredholm. Then in the notation of
Definition \ref{cr2def3}, we have
\e
\ind\bigl(\De^p_{k,{\bs\be}}\bigr)=
-\sum_{i=1}^nN_\sSii(\be_i).
\label{cr2eq12}
\e
\label{cr2thm6}
\end{thm}

\begin{proof} First we prove the special case with $\be_i=\ha(2-m)$
for $i=1,\ldots,n$. From Definition \ref{cr2def3} we have
$\D_\sSii\cap(2-m,0)=\emptyset$ and $N_\sSii\equiv 0$
on $(2-m,0)$, so ${\bs\be}$ lies in \eq{cr2eq9} and $N_\sSii
(\be_i)=0$ for $i=1,\ldots,n$ as $m>2$. Theorems \ref{cr2thm4} and
\ref{cr2thm5} then show that $\Coker\bigl(\De^p_{k,{\bs\be}}\bigr)
\cong\Ker\bigl(\De^p_{k,{\bs\be}}\bigr)^*$. Hence $\ind\bigl(
\De^p_{k,{\bs\be}}\bigr)=0$. As $N_\sSii(\be_i)=0$ for
$i=1,\ldots,n$ this proves \eq{cr2eq12} when~$\be_i=\ha(2-m)$.

Now Lockhart and McOwen prove formulae for how $\ind\bigl(
\De^p_{k,{\bs\be}}\bigr)$ changes as $\bs\be$ varies in
\eq{cr2eq9}. They do this for the $n=1$ case in \cite[Th.~6.2]{LoMc},
and for all $n$ but with exactly rather than asymptotically conical
metrics near $x_i$  in \cite[Th.~8.1]{LoMc}; their proof is easily
generalized to the asymptotically conical case.

In the $n=1$ case, \cite[Th.~6.2]{LoMc} shows that if
$\be<\ga$ lie in $\R\sm\D_{\sst\Si_1}$ then
\begin{equation*}
\ind\bigl(\De^p_{k,\be}\bigr)-
\ind\bigl(\De^p_{k,\ga}\bigr)=N(\be,\ga),
\quad\text{where}\quad
N(\be,\ga)=\sum_{\al\in\D_{\sst\Si_1}\cap(\be,\ga)}d(\al).
\end{equation*}
Here $d(\al)$ is defined on \cite[p.~416]{LoMc}, and is
effectively the dimension of the space of all harmonic functions
on $C_{\sst\Si_1}'$ of the form \eq{cr2eq3}. But Proposition
\ref{cr2prop1} shows that such functions are homogeneous of order
$\al$ as $m>2$, so $d(\al)=m_\sSi(\al)$ in our notation.
Thus~$N(\be,\ga)=N_{\sst\Si_1}(\ga)-N_{\sst\Si_1}(\be)$.

Similarly, for $n\ge 1$ we find using \cite[Th.~8.1]{LoMc} that
if ${\bs\be},{\bs\ga}$ lie in \eq{cr2eq9} then
\begin{equation*}
\ind\bigl(\De^p_{k,{\bs\be}}\bigr)-\ind\bigl(
\De^p_{k,{\bs\ga}}\bigr)
=\sum_{i=1}^nN_\sSii(\ga_i)-\sum_{i=1}^nN_\sSii(\be_i).
\end{equation*}
Combining this with the case $\ga_i=\ha(2-m)$ for $i=1,\ldots,n$
proved above yields~\eq{cr2eq12}.
\end{proof}

We identify $\Ker(\De^p_{k,{\bs\be}})$ in simple cases.

\begin{lem} Let\/ $(X,d)$ be a compact Riemannian $m$-manifold
with conical singularities $x_1,\ldots,x_n$, and let\/ $p>1$,
$k\ge 2$ and\/ ${\bs\be}\in\R^n$. Then
\begin{itemize}
\item[{\rm(a)}] If\/ $\be_i>0$ for $i=1,\ldots,n$
then~$\Ker\bigl(\De^p_{k,{\bs\be}}\bigr)=\{0\}$.
\item[{\rm(b)}] If\/ $\be_i\in(2-m,0)$ for $i=1,\ldots,n$
then~$\Ker\bigl(\De^p_{k,{\bs\be}}\bigr)=\an{1}$.
\end{itemize}
\label{cr2lem4}
\end{lem}

\begin{proof} Let $u\in\Ker\bigl(\De^p_{k,{\bs\be}}\bigr)$,
so that $u\in C^\iy_{\bs\be}(X')$ by Theorem \ref{cr2thm4}.
If $\be_i>0$ this implies that $u(y)\ra 0$ as $y\ra x_i$
for $i=1,\ldots,n$. Applying the {\it maximum principle}
\cite[\S 3]{GiTr} then shows that $u=0$, proving (a). For
(b), first let $\be_i=\ha(2-m)$ for $i=1,\ldots,n$, and
suppose $u\in\Ker(\De^p_{k,{\bs\be}})$. Then Lemma
\ref{cr2lem3} gives
\begin{equation*}
0=\int_{X'}u(\De u)\,\d V_g=\int_{X'}\ms{\d u}\d V_g.
\end{equation*} 
Thus $\d u=0$, so $u$ is constant, and $u\in\an{1}$ as $X'$
is connected.

Conversely, $1\in\Ker(\smash{\De^p_{k,{\bs\be}}})$ as $\be_i<0$,
so $\Ker(\smash{\De^p_{k,{\bs\be}}})=\an{1}$. Now $(2-m,0)^n$ is
a connected subset of \eq{cr2eq9} containing $\bigl(\ha(2-m),
\ldots,\ha(2-m)\bigr)$, by Definition \ref{cr2def3}. Thus
Theorem \ref{cr2thm4} shows that $\Ker(\De^p_{k,{\bs\be}})$
is independent of $\bs\be$ for $\be_i\in(2-m,0)$, and part
(b) follows.
\end{proof}

The next inequality is in effect a lower bound for the positive
eigenvalues of the Laplacian $\De$ on $X'$. Here $\lnm{u}{2}$ is
the {\it unweighted\/} $L^2$-norm~$\bigl(\int_{X'}u^2\d V
\bigr)^{\smash{1/2}}$.

\begin{thm} Let\/ $(X,d)$ be a compact Riemannian $m$-manifold
for $m>2$ with conical singularities $x_1,\ldots,x_n$, and
suppose $X'=X\sm\{x_1,\ldots,x_n\}$ is connected. Then there
exists $C>0$ such that whenever $u\in C^2(X')$ is
compactly-supported with\/ $\int_{X'}u\,\d V_g=0$ we
have~$\lnm{u}{2}\le C\lnm{\d u}{2}\le C^2\lnm{\De u}{2}$.
\label{cr2thm7}
\end{thm}

\begin{proof} Let $\be_i=\ha(2-m)$ for $i=1,\ldots,n$. Then
$\bs\be$ lies in \eq{cr2eq9} by Definition \ref{cr2def3}, so 
$\De^2_{2,{\bs\be}}$ is Fredholm by Theorem \ref{cr2thm3},
and Theorem \ref{cr2thm5} and part (b) of Lemma \ref{cr2thm4}
show that $\De^2_{2,{\bs\be}}$ has kernel and cokernel $\an{1}$.
Therefore
\e
\De^2_{2,{\bs\be}}:\bigl\{u\in
L^2_{2,{\bs\be}}(X'):
\ts\int_{X'}u\,\d V_g=0\bigr\}\ra V
\label{cr2eq13}
\e
is a continuous vector space isomorphism between Banach spaces, where
\begin{equation*}
V=\bigl\{v\in L^2_{0,{\bs\be}-2}(X'):\ts\int_{X'}v\,\d V_g=0\bigr\}
\end{equation*}
is a Hilbert space, with the $L^2_{0,{\bs\be}-2}$ norm.

By the Open Mapping Theorem it follows that \eq{cr2eq13} has a
continuous inverse, $P$ say. Let $\iota:L^2_{2,{\bs\be}}(X')\ra
L^2_{0,{\bs\be}-2}(X')$ be the inclusion. Then $\iota$ is
continuous and compact, by part (a) of Theorem \ref{cr2thm1}.
Hence $\iota\circ P:V\ra V$ is a {\it continuous, injective,
compact, linear automorphism} of the Hilbert space~$V$.

We are interested in the ordinary, unweighted $L^2$-norm
$\lnm{\,.\,}{2}$ on functions. Define $\ga_i=-\ha m$ for
$i=1,\ldots,n$. Then the power of $\rho$ in \eq{cr2eq7} used to
define $L^2_{0,{\bs\ga}}(X')$ is trivial, so that $L^2(X')=L^2_{0,
{\bs\ga}}(X')$ with $\lnm{\,.\,}{2}=\snm{\,.\,}{2}{0,{\bs\ga}}$.
As $L^2_{2,{\bs\be}}(X')\hookra L^2_{0,{\bs\ga}}(X')$ is a
continuous inclusion by Theorem \ref{cr2thm1}, we see that
$\iota\circ P$ maps~$V\ra V\cap L^2(X')$.

Now $\iota\circ P$ is an inverse for the Laplacian $\De$, and
$\De$ is {\it self-adjoint\/} w.r.t.\ the $L^2$ inner product
(though not w.r.t.\ the $L^2_{0,{\bs\be}-2}$ inner product).
It easily follows that the restriction of $\iota\circ P$ to
$V\cap L^2(X')$ is self-adjoint w.r.t.\ the $L^2$ inner product
on~$V\cap L^2(X')$.

We can now apply the theory of compact self-adjoint operators
on Hilbert spaces. As $\iota\circ P$ is compact it has a
countable set of eigenvalues converging to zero, with finite
multiplicity. As $\iota\circ P$ is injective, the eigenvalues
are nonzero. Thus all eigenspaces lie in $\iota\circ P(V)\subset
V\cap L^2(X')$. As $\iota\circ P$ is self-adjoint in the $L^2$
inner product the eigenvalues are all real, there are no
nilpotency phenomena, and eigenvectors for distinct eigenvalues
are $L^2$-orthogonal.

Therefore there exists a sequence $(e_i)_{i=1}^\iy$ of eigenvectors
of $\iota\circ P$ in $V\cap L^2(X')$, such that $\iota\circ P(e_i)=
\la_i$ for $\la_i\in\R\sm\{0\}$ with $\la_i\ra 0$ as $i\ra\iy$, and
$(e_i)_{i=1}^\iy$ is orthonormal in the $L^2$ inner product, and
$(e_i)_{i=1}^\iy$ is a basis for the Hilbert space $V$. As $\iota
\circ P$ is an inverse for $\De$ we see that $\De e_i=\la_i^{-1}e_i$.
Thus
\begin{equation*}
0<\int_{X'}\ms{\d e_i}\,\d V_g=\int_{X'}e_i\De e_i\,\d V_g=
\la_i^{-1}\lnm{e_i}{2}^2=\la_i^{-1},
\end{equation*}
integrating by parts using Lemma \ref{cr2lem3}. So $\la_i>0$ for
all~$i$.

Now let $u\in C^2(X')$ be compactly-supported with $\int_{X'}u\,
\d V_g=0$. Then $u\in V\cap L^2(X')$ and $\De u\in V\cap L^2(X')$.
Thus $\De u=\sum_{i=1}^\iy x_ie_i$, where $x_i=\an{\De u,e_i}_{L^2}
\in\R$. As $\iota\circ P$ is an inverse for $\De$ we have
$u=\iota\circ P(\De u)$, so $u=\sum_{i=1}^\iy\la_ix_ie_i$. Hence
\begin{equation*}
\lnm{u}{2}^2=\sum_{i=1}^\iy\la_i^2x_i^2,\;\>
\lnm{\d u}{2}^2=\an{u,\De u}_{L^2}=\sum_{i=1}^\iy\la_ix_i^2
\;\>\text{and}\;\> \lnm{\De u}{2}^2=\sum_{i=1}^\iy x_i^2.
\end{equation*}
Since $\la_i>0$ with $\la_i\ra 0$ as $i\ra\iy$ this implies
that $\lnm{u}{2}\le C\lnm{\d u}{2}\le C^2\lnm{\De u}{2}$ with
$C=\min_{i=1}^\iy(\la_i^{-1})$. This completes the proof.
\end{proof}

\subsection{Homology and cohomology}
\label{cr24}

Next we discuss homology and cohomology of manifolds with conical
singularities. For a general reference on (co)homology of manifolds,
see for instance Bredon \cite{Bred}. If $Y$ is a manifold, write
$H^k(Y,\R)$ for the $k^{\rm th}$ {\it de Rham cohomology group} and
$H^k_{\rm cs}(Y,\R)$ for the $k^{\rm th}$ {\it compactly-supported
de Rham cohomology group} of $Y$. If $Y$ is compact then~$H^k(Y,\R)
=H^k_{\rm cs}(Y,\R)$.

Let $Y$ be a topological space, and $Z\subset Y$ a subspace. Write
$H_k(Y,\R)$ for the $k^{\rm th}$ {\it real singular homology group}
of $Y$, and $H_k(Y;Z,\R)$ for the $k^{\rm th}$ {\it real singular
relative homology group} of $(Y,Z)$. When $Y$ is a manifold and $Z$
a submanifold, we may define $H_k(Y,\R)$ and $H_k(Y;Z,\R)$ using
{\it smooth\/} simplices, as in \cite[\S V.5]{Bred}. Then the pairing
between (singular) homology and (de Rham) cohomology is defined at
the chain level by integrating $k$-forms over $k$-simplices.

Suppose $Y$ is a compact $m$-manifold with boundary, so that
$\pd Y$ is a compact $(m-1)$-manifold and $Y^\circ=Y\sm\pd Y$ is
an $m$-manifold without boundary, which is noncompact if $\pd Y
\ne\emptyset$. Then there is a natural long exact sequence
\e
\cdots\ra
H^k_{\rm cs}(Y^\circ,\R)\ra H^k(Y^\circ,\R)\ra
H^k(\pd Y,\R)\ra H^{k+1}_{\rm cs}(Y^\circ,\R)\ra\cdots.
\label{cr2eq14}
\e
Note that $H^k(Y^\circ,\R)=H^k(Y,\R)$. Suppose $Y$ is oriented.
Then by Poincar\'e--Lefschetz duality there are isomorphisms
\e
H_k(Y;\pd Y,\R)^*\cong H^k_{\rm cs}(Y^\circ,\R)\cong
H_{m-k}(Y^\circ,\R)\cong H^{m-k}(Y^\circ,\R)^*.
\label{cr2eq15}
\e

If $X$ is a compact Riemannian manifold with conical
singularities $x_1,\ldots,x_n$ then $X'=X\sm\{x_1,\ldots,x_n\}$
is the interior of a compact manifold $\bar X'$ with boundary
$\pd\bar X'$ the disjoint union $\coprod_{i=1}^n\Si_i$. Thus
\eq{cr2eq14} gives an exact sequence
\e
\cdots\ra
H^k_{\rm cs}(X',\R)\ra H^k(X',\R)\ra\bigoplus_{i=1}^n
H^k(\Si_i,\R)\ra H^{k+1}_{\rm cs}(X',\R)\ra\cdots.
\label{cr2eq16}
\e
If $X'$ is oriented then \eq{cr2eq15} gives isomorphisms
\e
H_k\bigl(X;\{x_1,\ldots,x_n\},\R\bigr)^*\cong
H^k_{\rm cs}(X',\R)\cong H_{m-k}(X',\R)\cong H^{m-k}(X',\R)^*,
\label{cr2eq17}
\e
as $H_{m-k}(X;\{x_1,\ldots,x_n\},\R)\cong H_{m-k}(\bar X';
\pd\bar X',\R)$ by excision.

Since $H_k\bigl(\{x_1,\ldots,x_n\},\R\bigr)=0$ for $k\ne 0$, the
long exact sequence
\begin{equation*}
\cdots\ra H_k\bigl(X,\R\bigr)\ra H_k\bigl(X;\{x_1,\ldots,x_n\},\R\bigr)
\ra H_{k-1}\bigl(\{x_1,\ldots,x_n\},\R\bigr)\ra\cdots
\end{equation*}
implies that $H_k\bigl(X,\R\bigr)\cong H_k\bigl(X;\{x_1,\ldots,
x_n\},\R\bigr)$ for $k>1$. Therefore \eq{cr2eq17} gives
\e
H^k_{\rm cs}(X',\R)\cong H_k(X,\R)^*
\quad\text{for all $k>1$.}
\label{cr2eq18}
\e

We can now study $\De^p_{k,{\bs\la}}$ when $\la_i$ is small and positive.

\begin{prop} Suppose $(X,d)$ is a compact Riemannian $m$-manifold for
$m>2$ with conical singularities $x_1,\ldots,x_n$, and use the notation
above. Let\/ ${\cal K}_\sXp\subset C^\iy(X')$ be a vector space of
smooth functions constant on $S_i$ for $i=1,\ldots,n$, such that\/
$v\mapsto[\d v]$ is an isomorphism from ${\cal K}_\sXp$ to the kernel
of the map $H^1_{\rm cs}(X',\R)\ra H^1(X',\R)$ in \eq{cr2eq16}. Let\/
$p>1$, $k\ge 2$ and\/ $\la_i>0$ with\/ $(0,\la_i]\cap\D_\sSii=
\emptyset$ for $i=1,\ldots,n$. Then
\e
\De:L^p_{k,{\bs\la}}(X')\op{\cal K}_\sXp\ra\bigl\{w\in
L^p_{k-2,{\bs\la}-2}(X'):\ts\int_{X'}w\,\d V_g=0\bigr\},
\label{cr2eq19}
\e
given by $(u,v)\mapsto\De^p_{k,{\bs\la}}u+\De v$, is an
isomorphism of topological vector spaces.
\label{cr2prop2}
\end{prop}

\begin{proof} Since ${\cal K}_\sXp$ is isomorphic to the kernel of
$H^1_{\rm cs}(X',\R)\ra H^1(X',\R)$, equation \eq{cr2eq16} gives
an exact sequence
\begin{equation*}
0\ra H^0(X',\R)\ra\bigoplus_{i=1}^n
H^0(\Si_i,\R)\ra{\cal K}_\sXp\ra 0,
\end{equation*}
and thus $\dim{\cal K}_\sXp=\sum_{i=1}^nb^0(\Si_i)-1$ as $X'$
is connected. As $\la_i\notin\D_\sSii$, Theorem
\ref{cr2thm6} shows that $\De^p_{k,{\bs\la}}$ is Fredholm with
\e
\ind(\De^p_{k,{\bs\la}})\!=\!\sum_{i=1}^nN_\sSii(\la_i)
\!=\!\sum_{i=1}^nN_\sSii(0)\!=\!\sum_{i=1}^nb^0(\Si_i)\!=\!
1\!+\!\dim{\cal K}_\sXp.
\label{cr2eq20}
\e
Here $N_\sSii(\la_i)=N_\sSii(0)$ as $(0,\la_i]
\cap\D_\sSii=\emptyset$ and $N_\sSii$ is upper
semicontinuous and locally constant on $\R\sm\D_\sSii$.
Also $N_\sSii(0)$ is the multiplicity of the eigenvalue 0 of
$\De_\sSii$ by Definition \ref{cr2def3}, which is~$b^0(\Si_i)$.

By integrating by parts as in \eq{cr2eq10} we see that
$\int_{X'}\De(u+v)\d V_g=0$ for $u\in L^p_{k,{\bs\la}}(X')$ and
$v\in{\cal K}_\sXp$, so $\De$ does map into the given r.h.s.\ in
\eq{cr2eq19}. Now \eq{cr2eq19} modifies the Fredholm map
$\De^p_{k,{\bs\la}}$, increasing the dimension of its domain by
$\dim{\cal K}_\sXp$, and decreasing the dimension of its range by
1. Therefore from \eq{cr2eq20} we see that \eq{cr2eq19} is
Fredholm with index 0. Thus \eq{cr2eq19} is an isomorphism of
topological vector spaces if and only if it is injective.

Suppose $(u,v)$ lies in the kernel of \eq{cr2eq19}. Then
$\De(u+v)\equiv 0$, so multiplying by $u+v$ and integrating by
parts as in \eq{cr2eq10} shows that $\int_{X'}\ms{\d u+\d v}\d V_g
=0$, so $\d(u+v)=0$ and $u+v\equiv c$ for some $c\in\R$. Now $u\in
C^0_{\bs\la}(X')$ by Theorem \ref{cr2thm1}, so that $u(x)\ra 0$ as
$x\ra x_i$, and $v$ is constant on~$S_i$.

Taking $x\ra x_i$ shows that $v\equiv c$ on $S_i$ for all $i$,
so $v-c$ is compactly supported. But then $[\d v]=[\d(v-c)]=0$ in
$H^1_{\rm cs}(X',\R)$, so $v=0$ as $v\mapsto[\d v]$ is an
isomorphism with a subspace of $H^1_{\rm cs}(X',\R)$. Hence
$c=0$, so $u=0$, and \eq{cr2eq19} is injective. This completes
the proof.
\end{proof}

\subsection{Hodge theory}
\label{cr25}

{\it Hodge theory} for a compact Riemannian manifold $(Y,g)$
shows that each class in $H^k(Y,\R)$ is represented by a unique
$k$-form $\al$ with $\d\al=\d^*\al=0$. Here is an analogue of
this on $X'$ for $k=1$, with decay conditions.

\begin{thm} Let\/ $X$ be a compact Riemannian $m$-manifold
for $m>2$ with conical singularities at\/ $x_1,\ldots,x_n$,
and let\/ $X',\ep,\Si_i,\phi_i,S_i$ and\/ $\nu_i$ be as in
Definition \ref{cr2def1}, $\D_\sSii$ as in Definition
\ref{cr2def3}, and\/ $\rho$ as in Definition \ref{cr2def4}. Define
\e
Y_\sXp\!=\!\bigl\{\al\in C^\iy(T^*X'):\text{$\d\al\!=\!0$,
$\d^*\al\!=\!0$, $\md{\na^k\al}\!=\!O(\rho^{-1-k})$ for $k\ge 0$}\bigr\}.
\label{cr2eq21}
\e
Then $\pi:Y_\sXp\!\ra\!H^1(X',\R)$ given by $\pi:\al\!\mapsto\![\al]$
is an isomorphism. Furthermore:
\begin{itemize}
\item[{\rm(a)}] Fix $\al\in Y_\sXp$. By Hodge theory there exists
a unique $\ga_i\in C^\iy(T^*\Si_i)$ with\/ $\d\ga_i=\d^*\ga_i=0$
for $i=1,\ldots,n$, such that the image of\/ $\pi(\al)$ under the
map $H^1(X',\R)\ra\bigoplus_{i=1}^nH^1(\Si_i,\R)$ of\/ \eq{cr2eq16}
is $\bigl([\ga_1],\ldots,[\ga_n]\bigr)$. There exist unique $T_i\in
C^\iy\bigl(\Si_i\t(0,\ep)\bigr)$ for $i=1,\ldots,n$ such that
\end{itemize}
\ea
\phi_i^*(\al)&=\pi_i^*(\ga_i)+\d T_i
\quad \text{on $\Si_i\t(0,\ep)$ for $i=1,\ldots,n$, and}
\label{cr2eq22}\\
\na^kT_i(\si,r)&=O(r^{\la_i-k})\qquad
\begin{aligned}
&\text{as\/ $r\ra 0$, for all\/ $k\ge 0$ and}\\
&\text{$\la_i\in(0,\nu_i)$ with\/ $(0,\la_i]\cap\D_\sSii=\emptyset$.}
\end{aligned}
\label{cr2eq23}
\ea
\begin{itemize}
\item[{\rm(b)}] Suppose $\ga_i\in C^\iy(T^*\Si_i)$ with\/
$\d\ga_i=\d^*\ga_i=0$ for $i=1,\ldots,n$, and the image
of\/ $\bigl([\ga_1],\ldots,[\ga_n]\bigr)$ under
$\bigoplus_{i=1}^nH^1(\Si_i,\R)\ra H^2_{\rm cs}(X',\R)$ in
\eq{cr2eq16} is $[\be]$ for some exact\/ $2$-form $\be$ on $X'$
supported on $X'\sm(S_1\cup\cdots\cup S_n)$. Then there exists
$\al\in C^\iy(T^*X')$ with\/ $\d\al=\be$, $\d^*\al=0$ and\/
$\md{\na^k\al}=O(\rho^{-1-k})$ for $k\ge 0$, such that\/
\eq{cr2eq22} and\/ \eq{cr2eq23} hold for~$T_i\in C^\iy
\bigl(\Si_i\t(0,\ep)\bigr)$.
\item[{\rm(c)}] Let\/ $f\in C^\iy(X')$ with\/ $\md{\na^kf}=
O(\rho^{{\bs\nu}-2-k})$ for $k\ge 0$ and\/ $\int_{X'}f\,\d V=0$.
Then there exists a unique exact\/ $1$-form $\al$ on $X'$ with\/
$\d^*\al=f$ and\/ $\md{\na^k\al}=O(\rho^{-1-k})$ for $k\ge 0$,
such that\/ \eq{cr2eq22} and\/ \eq{cr2eq23} hold for
$\ga_i=0$ and\/~$T_i\in C^\iy\bigl(\Si_i\t(0,\ep)\bigr)$.
\end{itemize}
\label{cr2thm8}
\end{thm}

\begin{proof} Clearly $Y_\sXp$ is a vector space and $\pi$ is
linear. We must show that $\pi$ is injective and surjective. Suppose
$\al\in Y_\sXp$ and $[\al]=0\in H^1(X',\R)$. Then $\al=\d\th$
for $\th\in C^\iy(X')$. Using $\na^k\th=\na^{k-1}\al$ for $k>0$, and
integrating $\md{\al}=O(\rho^{-1})$ to estimate $\md{\th}$, gives
\e
\th=O\bigl(1+\md{\log\rho}\bigr)\quad\text{and}\quad
\na^k\th=O(\rho^{-k})\quad\text{for all $k\ge 0$.}
\label{cr2eq24}
\e
Now $\d^*\d\th=0$. Multiplying this by $\th$ and integrating over $X'$
by parts, using \eq{cr2eq24} and arguing as in Lemma \ref{cr2lem3},
we can show that $\int_{X'}\ms{\d\th}\d V_g=0$. Thus $\al=\d\th=0$, so
if $[\al]=0$ then $\al=0$, and $\pi$ is injective.

Next we show $\pi$ is surjective, and at the same time prove part (a).
Let $\eta\in H^1(X',\R)$. By Hodge theory there exists a unique $\ga_i
\in C^\iy(T^*\Si_i)$ with $\d\ga_i=\d^*\ga_i=0$ for $i=1,\ldots,n$,
such that the image of $\eta$ under $H^1(X',\R)\ra\bigoplus_{i=1}^n
H^1(\Si_i,\R)$ in \eq{cr2eq16} is $\bigl([\ga_1],\ldots,[\ga_n]\bigr)$.

Choose a smooth, closed 1-form $\ga$ on $X'$ with $[\ga]=\eta$ and
$\phi_i^*(\ga)=\pi_i^*(\ga_i)$, where $\pi_i:\Si_i\t(0,\ep)\ra\Si_i$
is the obvious projection. Note that the condition $d(x_i,x_j)>2\ep$
for $i\ne j$ in part (b) of Definition \ref{cr2def1} implies that
the closures $\bar S_1,\ldots,\bar S_n$ are disjoint in $X$, and
using this we can show that $\ga$ exists.

As $\phi_i^*(\ga)=\pi_i^*(\ga_i)$ we can regard $\ga$ as independent
of $r$ on $S_i\cong\Si_i\t(0,\ep)\ni(\si,r)$. Since the metric $g$
on $\Si_i\t(0,\ep)$ is approximately the cone metric by \eq{cr2eq1},
we find that $\md{\ga}=O(\rho^{-1})$, and more generally
\e
\md{\na^k\ga}=O(\rho^{-1-k})\quad\text{for $k\ge 0$.}
\label{cr2eq25}
\e
This suggests that $\d^*\ga=O(\rho^{-2})$. However, because
$\d^*\ga_i=0$ on $\Si_i$ we have $\d^*(\pi_i^*(\ga_i))=0$ on
$\Si_i\t(0,\ep)$, computing $\d^*$ w.r.t.\ the cone metric on
$\Si_i\t(0,\ep)$. Since $g$ approximates the cone metric on
$S_i$, calculation using \eq{cr2eq1} shows that
\e
\bmd{\na^k(\d^*\ga)}=O(\rho^{{\bs\nu}-2-k})\quad\text{for $k\ge 0$,}
\label{cr2eq26}
\e
where ${\bs\nu}=(\nu_1,\ldots,\nu_n)$. Thus~$\d^*\ga\in
C^\iy_{{\bs\nu}-2}(X')$.

Choose $\la_i\in(0,\nu_i)$ with $(0,\la_i]\cap\D_\sSii
=\emptyset$ for $i=1,\ldots,n$. Let $p>1$ and $k\ge 2$. Then
\eq{cr2eq26} implies that $\d^*\ga\in L^p_{k-2,{\bs\la}-2}(X')$.
Integrating by parts as in \eq{cr2eq10} shows that $\int_{X'}
(\d^*\ga)\d V_g=\int_{X'}1(\d^*\ga)\d V_g=\int_{X'}(\d 1)\cdot
\ga\,\d V_g=0$, using $m>2$ and \eq{cr2eq25} for~$k=0,1$.

Thus $\d^*\ga$ lies in the r.h.s.\ of \eq{cr2eq19}, and by
Proposition \ref{cr2prop2} there exist unique $u\in
L^p_{k,{\bs\la}}(X')$ and $v\in{\cal K}_\sXp$ with
$\d^*\ga=\De u+\De v$. As $\De u=\d^*\ga-\De v$ and $\d^*\ga,
\De v\in L^p_{k-2,{\bs\la}-2}(X')$ for all $k\ge 2$, Theorem
\ref{cr2thm2} shows that $u\in L^p_{k,{\bs\la}}(X')$ for all
$k\ge 2$, so that $u\in C^\iy_{\bs\la}(X')$ by Theorem~\ref{cr2thm1}.

Define $\al=\ga-\d u-\d v$. Then $\d\al=0$ and $\d^*\al=\d^*\ga-\De
u-\De v=0$. As $\ga$ satisfies \eq{cr2eq25}, $u\in C^\iy_{\bs\la}(X')$
with $\la_i>0$ and $\d v$ is compactly-supported, we see that
$\md{\na^k\al}=O(\rho^{-1-k})$ for all $k\ge 0$. Hence $\al\in
Y_\sXp$ and $[\al]=[\ga]=\eta$, so $\pi:Y_\sXp\ra H^1(X',\R)$
is surjective. This proves the first part of the theorem. Define
$T_i=-\phi_i^*(u)$. Then \eq{cr2eq22} holds as $\al=\ga-\d u-\d v$,
$\d v=0$ on $S_i$ and $\phi_i^*(\ga)=\pi_i^*(\ga_i)$, and \eq{cr2eq23}
holds as $u\in C^\iy_{\bs\la}(X')$ whenever $\la_i\in(0,\nu_i)$
with $(0,\la_i]\cap\D_\sSii=\emptyset$. This proves part~(a).

For part (b), let $\ga_i$ and $\be$ be as in the theorem.
Choose $\ga\in C^\iy(T^*X')$ with $\phi_i^*(\ga)=\pi_i^*(\ga_i)$
for $i=1,\ldots,n$. Then $\d\ga$ is supported in $X'\sm(S_1\cup
\cdots\cup S_n)$ as $\d\ga_i=0$. By construction $[\d\ga]\in
H^2_{\rm cs}(X',\R)$ is the image of $\bigl([\ga_1],\ldots[\ga_n]
\bigr)$ under $\bigoplus_{i=1}^nH^1(\Si_i,\R)\ra H^2_{\rm cs}
(X',\R)$, so~$[\d\ga]=[\be]$.

Thus $\be=\d\ga+\d\de$, for some compactly-supported 1-form $\de$
on $X'$. Since $\be,\ga$ are supported on $X'\sm(S_1\cup\cdots\cup
S_n)$ we can choose $\de$ supported there too. As in part (a) there
exist unique $u\in C^\iy_{\bs\la}(X')$ and $v\in{\cal K}_\sXp$ with
$\De u+\De v=\d^*(\ga+\de)$. Then $\al=\ga+\de-\d u-\d v$ and $T_i=
-\phi_i^*(u)$ satisfy the conditions in~(b).

For part (c), choose $\la_i\in(0,\nu_i)$ with $(0,\la_i]\cap
\D_\sSii=\emptyset$ for $i=1,\ldots,n$. Let $p>1$ and $k\ge 2$.
Then $\md{\na^jf}=O(\rho^{{\bs\nu}-2-j})$ and $\int_{X'}f\,\d V=0$
imply that $f$ lies in the r.h.s.\ of \eq{cr2eq19}. So by Proposition
\ref{cr2prop2} there exist unique $u\in L^p_{k,{\bs\la}}(X')$ and
$v\in{\cal K}_\sXp$ with $\De u+\De v=\d^*f$. As $u,v$ are independent
of $k\ge 2$, Theorem \ref{cr2thm1} shows that $u\in C^\iy_{\bs\la}(X')$.
Thus $\al=\d u+\d v$ and $T_i=\phi_i^*(u)$ satisfy the conditions
in (c). This completes the proof.
\end{proof}

\section{Special Lagrangian geometry}
\label{cr3}

We now introduce special Lagrangian submanifolds (SL $m$-folds)
in two different geometric contexts. First, in \S\ref{cr31}, we
define SL $m$-folds in $\C^m$. Then \S\ref{cr32} discusses SL
$m$-folds in {\it almost Calabi--Yau $m$-folds}, compact
K\"ahler manifolds with a holomorphic volume form which
generalize Calabi--Yau manifolds.

Then \S\ref{cr33} defines {\it special Lagrangian $m$-folds
with conical singularities} in almost Calabi--Yau $m$-folds,
which are the subject of the paper. Some references for
\S\ref{cr31}--\S\ref{cr32} are Harvey and Lawson \cite{HaLa}
and the author~\cite{Joyc2}.

\subsection{Special Lagrangian submanifolds in $\C^m$}
\label{cr31}

We begin by defining {\it calibrations} and {\it calibrated 
submanifolds}, following Harvey and Lawson~\cite{HaLa}.

\begin{dfn} Let $(M,g)$ be a Riemannian manifold. An {\it oriented
tangent $k$-plane} $V$ on $M$ is a vector subspace $V$ of
some tangent space $T_xM$ to $M$ with $\dim V=k$, equipped
with an orientation. If $V$ is an oriented tangent $k$-plane
on $M$ then $g\vert_V$ is a Euclidean metric on $V$, so 
combining $g\vert_V$ with the orientation on $V$ gives a 
natural {\it volume form} $\vol_V$ on $V$, which is a 
$k$-form on~$V$.

Now let $\vp$ be a closed $k$-form on $M$. We say that
$\vp$ is a {\it calibration} on $M$ if for every oriented
$k$-plane $V$ on $M$ we have $\vp\vert_V\le \vol_V$. Here
$\vp\vert_V=\al\cdot\vol_V$ for some $\al\in\R$, and 
$\vp\vert_V\le\vol_V$ if $\al\le 1$. Let $N$ be an 
oriented submanifold of $M$ with dimension $k$. Then 
each tangent space $T_xN$ for $x\in N$ is an oriented
tangent $k$-plane. We say that $N$ is a {\it calibrated 
submanifold\/} if $\vp\vert_{T_xN}=\vol_{T_xN}$ for all~$x\in N$.
\label{cr3def1}
\end{dfn}

It is easy to show that calibrated submanifolds are automatically
{\it minimal submanifolds} \cite[Th.~II.4.2]{HaLa}. Here is the 
definition of special Lagrangian submanifolds in $\C^m$, taken
from~\cite[\S III]{HaLa}.

\begin{dfn} Let $\C^m$ have complex coordinates $(z_1,\dots,z_m)$, 
and define a metric $g'$, a real 2-form $\om'$ and a complex $m$-form 
$\Om'$ on $\C^m$ by
\e
\begin{split}
g'=\ms{\d z_1}+\cdots+\ms{\d z_m},\quad
\om'&=\ts\frac{i}{2}(\d z_1\w\d\bar z_1+\cdots+\d z_m\w\d\bar z_m),\\
\text{and}\quad\Om'&=\d z_1\w\cdots\w\d z_m.
\end{split}
\label{cr3eq1}
\e
Then $\Re\Om'$ and $\Im\Om'$ are real $m$-forms on $\C^m$. Let $L$
be an oriented real submanifold of $\C^m$ of real dimension $m$. We
say that $L$ is a {\it special Lagrangian submanifold\/} of $\C^m$,
or {\it SL\/ $m$-fold}\/ for short, if $L$ is calibrated with respect
to $\Re\Om'$, in the sense of Definition~\ref{cr3def1}.
\label{cr3def2}
\end{dfn}

Harvey and Lawson \cite[Cor.~III.1.11]{HaLa} give the following
alternative characterization of special Lagrangian submanifolds:

\begin{prop} Let\/ $L$ be a real $m$-dimensional submanifold 
of\/ $\C^m$. Then $L$ admits an orientation making it into an
SL submanifold of\/ $\C^m$ if and only if\/ $\om'\vert_L\equiv 0$ 
and\/~$\Im\Om'\vert_L\equiv 0$.
\label{cr3prop1}
\end{prop}

Thus SL $m$-folds are {\it Lagrangian submanifolds} in
$\R^{2m}\cong\C^m$ satisfying the extra condition that
$\Im\Om'\vert_L\equiv 0$, which is how they get their name.

\subsection{Almost Calabi--Yau $m$-folds and SL $m$-folds} 
\label{cr32}

We shall define special Lagrangian submanifolds not just in
Calabi--Yau manifolds, as usual, but in the much larger
class of {\it almost Calabi--Yau manifolds}.

\begin{dfn} Let $m\ge 2$. An {\it almost Calabi--Yau $m$-fold\/}
is a quadruple $(M,J,\om,\Om)$ such that $(M,J)$ is a compact
$m$-dimensional complex manifold, $\om$ is the K\"ahler form
of a K\"ahler metric $g$ on $M$, and $\Om$ is a non-vanishing
holomorphic $(m,0)$-form on~$M$.

We call $(M,J,\om,\Om)$ a {\it Calabi--Yau $m$-fold\/} if in
addition $\om$ and $\Om$ satisfy
\e
\om^m/m!=(-1)^{m(m-1)/2}(i/2)^m\Om\w\bar\Om.
\label{cr3eq2}
\e
Then for each $x\in M$ there exists an isomorphism $T_xM\cong\C^m$
that identifies $g_x,\om_x$ and $\Om_x$ with the flat versions
$g',\om',\Om'$ on $\C^m$ in \eq{cr3eq1}. Furthermore, $g$ is
Ricci-flat and its holonomy group is a subgroup of~$\SU(m)$.
\label{cr3def3}
\end{dfn}

This is not the usual definition of a Calabi--Yau manifold, but
is essentially equivalent to it.

\begin{dfn} Let $(M,J,\om,\Om)$ be an almost Calabi--Yau $m$-fold,
and $N$ a real $m$-dimensional submanifold of $M$. We call $N$ a
{\it special Lagrangian submanifold}, or {\it SL $m$-fold\/} for
short, if $\om\vert_N\equiv\Im\Om\vert_N\equiv 0$. It easily
follows that $\Re\Om\vert_N$ is a nonvanishing $m$-form on $N$.
Thus $N$ is orientable, with a unique orientation in which
$\Re\Om\vert_N$ is positive.
\label{cr3def4}
\end{dfn}

Again, this is not the usual definition of SL $m$-fold, but is
essentially equivalent to it. Suppose $(M,J,\om,\Om)$ is an
almost Calabi--Yau $m$-fold, with metric $g$. Let
$\psi:M\ra(0,\iy)$ be the unique smooth function such that
\e
\psi^{2m}\om^m/m!=(-1)^{m(m-1)/2}(i/2)^m\Om\w\bar\Om,
\label{cr3eq3}
\e
and define $\ti g$ to be the conformally equivalent metric $\psi^2g$
on $M$. Then $\Re\Om$ is a {\it calibration} on the Riemannian manifold
$(M,\ti g)$, and SL $m$-folds $N$ in $(M,J,\om,\Om)$ are calibrated
with respect to it, so that they are minimal with respect to~$\ti g$.

If $M$ is a Calabi--Yau $m$-fold then $\psi\equiv 1$ by \eq{cr3eq2},
so $\ti g=g$, and an $m$-submanifold $N$ in $M$ is special Lagrangian
if and only if it is calibrated w.r.t.\ $\Re\Om$ on $(M,g)$, as in
Definition \ref{cr3def2}. This recovers the usual definition of
special Lagrangian $m$-folds in Calabi--Yau $m$-folds.

\subsection{Special Lagrangian $m$-folds with conical singularities}
\label{cr33}

Now we can define {\it conical singularities} of SL $m$-folds.

\begin{dfn} Let $(M,J,\om,\Om)$ be an almost Calabi--Yau $m$-fold
for $m>2$, and define $\psi:M\ra(0,\iy)$ as in \eq{cr3eq3}. Suppose
$X$ is a compact singular SL $m$-fold in $M$ with singularities at
distinct points $x_1,\ldots,x_n\in X$, and no other singularities.

Fix isomorphisms $\up_i:\C^m\ra T_{x_i}M$ for $i=1,\ldots,n$
such that $\up_i^*(\om)=\om'$ and $\up_i^*(\Om)=\psi(x_i)^m\Om'$,
where $\om',\Om'$ are as in \eq{cr3eq1}. Let $C_1,\ldots,C_n$
be SL cones in $\C^m$ with isolated singularities at 0. For
$i=1,\ldots,n$ let $\Si_i=C_i\cap{\cal S}^{2m-1}$, and let
$\mu_i\in(2,3)$ with
\e
(2,\mu_i]\cap\D_\sSii=\emptyset,
\quad\text{where $\D_\sSii$ is defined in \eq{cr2eq4}.}
\label{cr3eq4}
\e
Then we say that $X$ has a {\it conical singularity} at $x_i$,
with {\it rate} $\mu_i$ and {\it cone} $C_i$ for $i=1,\ldots,n$,
if the following holds.

By Darboux's Theorem \cite[Th.~3.15]{McSa} there exist
embeddings $\Up_i:B_R\ra M$ for $i=1,\ldots,n$ satisfying
$\Up_i(0)=x_i$, $\d\Up_i\vert_0=\up_i$ and $\Up_i^*(\om)=\om'$,
where $B_R$ is the open ball of radius $R$ about 0 in $\C^m$ for
some small $R>0$. Define $\iota_i:\Si_i\t(0,R)\ra B_R$ by
$\iota_i(\si,r)=r\si$ for~$i=1,\ldots,n$.

Define $X'=X\sm\{x_1,\ldots,x_n\}$. Then there should exist a
compact subset $K\subset X'$ such that $X'\sm K$ is a union of
open sets $S_1,\ldots,S_n$ with $S_i\subset\Up_i(B_R)$, whose
closures $\bar S_1,\ldots,\bar S_n$ are disjoint in $X$. For
$i=1,\ldots,n$ and some $R'\in(0,R]$ there should exist a smooth
$\phi_i:\Si_i\t(0,R')\ra B_R$ such that $\Up_i\circ\phi_i:\Si_i
\t(0,R')\ra M$ is a diffeomorphism $\Si_i\t(0,R')\ra S_i$, and
\e
\bmd{\na^k(\phi_i-\iota_i)}=O(r^{\mu_i-1-k})
\quad\text{as $r\ra 0$ for $k=0,1$.}
\label{cr3eq5}
\e
Here $\na,\md{\,.\,}$ are computed using the cone metric
$\iota_i^*(g')$ on~$\Si_i\t(0,R')$.
\label{cr3def5}
\end{dfn}

We will show in Theorem \ref{cr5thm1} that if \eq{cr3eq5} holds
for $k=0,1$ then we can choose a natural $\phi_i$ for which
\eq{cr3eq5} holds for all $k\ge 0$. Thus the number of derivatives
required in \eq{cr3eq5} makes little difference, and we choose
$k=0,1$ to make the definition as weak as possible. We will also
show in Theorem \ref{cr5thm5} that if \eq{cr3eq5} holds for
{\it some} choice of rates $\mu_i$ satisfying the conditions of
the definition, then it holds for {\it all\/} choices of rates
$\mu_i$ satisfying the conditions, for the $\phi_i$ in Theorem
\ref{cr5thm1}. Thus the choice of rates $\mu_i$ again makes
little difference.

We restrict to $m>2$ in Definition \ref{cr3def5} for two reasons.
Firstly, the only SL cones $C$ in $\C^2$ are finite unions of SL
planes $\R^2$ in $\C^2$ intersecting only at 0. Therefore, SL
2-folds with conical singularities are actually {\it nonsingular}
as immersed 2-folds, so there is really no point in studying them.
Secondly, parts of the analysis in \S\ref{cr2} do not hold when $m=2$,
in particular Proposition \ref{cr2prop1} and Theorem \ref{cr2thm6}.
Therefore, in the rest of the paper we shall suppose~$m>2$.

Here are the reasons for the conditions on $\mu_i$ in
Definition~\ref{cr3def5}:
\begin{itemize}
\item We need $\mu_i>2$, or else \eq{cr3eq5} does not force
$X$ to approach $C_i$ near~$x_i$.
\item The definition involves a choice of $\Up_i:B_R\ra M$.
If we replace $\Up_i$ by a different choice $\ti\Up_i$ then
we should replace $\phi_i$ by $\ti\phi_i=(\ti\Up_i^{-1}\circ
\Up_i)\circ\phi_i$ near 0 in $B_R$. Calculation shows that
as $\Up_i,\ti\Up_i$ agree up to second order, we
have~$\bmd{\na^k(\ti\phi_i-\phi_i)}=O(r^{2-k})$.

Therefore if $\mu_i\le 3$ then \eq{cr3eq5} for $\phi_i$ is
equivalent to \eq{cr3eq5} for $\ti\phi_i$, and the definition
is independent of the choice of $\Up_i$. However, if $\mu_i>3$
then the definition would depend on the choice of $\Up_i$,
which we do not want. We also exclude $\mu_i=3$ for technical
reasons, to prevent $O(r^{2-k})$ terms from $\Up_i$ dominating
$\na^k(\phi_i-\iota_i)$, so we require $\mu_i<3$.
\item If we omit condition \eq{cr3eq4} then the proof of
Theorem \ref{cr5thm5} below would fail. Also, extra obstructions
would appear in the deformation theory of compact SL $m$-folds
with conical singularities studied in~\cite{Joyc3}.
\end{itemize}

To avoid proliferation of indices we have chosen $R,R'$ above
to be independent of $i=1,\ldots,n$. This is valid as we may
take $R=\min(R_1,\ldots,R_n)$, and so on. We will do this
without remark for other variables in later proofs.

\section{Lagrangian Neighbourhood Theorems}
\label{cr4}

Let $N$ be a real $m$-manifold. Then its tangent bundle $T^*N$ has
a canonical symplectic form $\hat\om$, defined as follows. Let
$(x_1,\ldots,x_m)$ be local coordinates on $N$. Extend them to
local coordinates $(x_1,\ldots,x_m,y_1,\ldots,y_m)$ on $T^*N$
such that $(x_1,\ldots,y_m)$ represents the 1-form $y_1\d x_1+
\cdots+y_m\d x_m$ in $T_{(x_1,\ldots,x_m)}^*N$. Then~$\hat\om=
\d x_1\w\d y_1+\cdots+\d x_m\w\d y_m$.

Identify $N$ with the zero section in $T^*N$. Then $N$ is a
{\it Lagrangian submanifold\/} of $T^*N$. The {\it Lagrangian
Neighbourhood Theorem} \cite[Th.~3.33]{McSa}, due to Weinstein
\cite{Wein}, shows that any compact Lagrangian submanifold $N$
in a symplectic manifold looks locally like the zero section
in~$T^*N$.

\begin{thm} Let\/ $(M,\om)$ be a symplectic manifold and\/
$N\subset M$ a compact Lagrangian submanifold. Then there
exists an open tubular neighbourhood\/ $U$ of the zero
section $N$ in $T^*N$, and an embedding $\Phi:U\ra M$ with\/
$\Phi\vert_N=\id:N\ra N$ and\/ $\Phi^*(\om)=\hat\om$, where
$\hat\om$ is the canonical symplectic structure on~$T^*N$.
\label{cr4thm1}
\end{thm}

We shall need the following variation of this, which may be
deduced from the proof of a result of Weinstein \cite[Th.~7.1]{Wein}
on {\it Lagrangian foliations}.

\begin{thm} Let\/ $(M,\om)$ be a $2m$-dimensional symplectic
manifold and\/ $N\subset M$ an embedded\/ $m$-dimensional submanifold.
Suppose $\{L_x:x\in N\}$ is a smooth family of embedded, noncompact
Lagrangian submanifolds in $M$ parametrized by $x\in N$, such that
for each\/ $x\in N$ we have $x\in L_x$, and\/~$T_xL_x\cap T_xN=\{0\}$.

Then there exists an open neighbourhood\/ $U$ of the zero section
$N$ in $T^*N$ such that the fibres of the natural projection
$\pi:U\ra N$ are connected, and a unique embedding $\Phi:U\ra M$
with\/ $\Phi\bigl(\pi^{-1}(x)\bigr)\subset L_x$ for each\/ $x\in N$,
$\Phi\vert_N=\id:N\ra N$ and\/ $\Phi^*(\om)=\hat\om+\pi^*(\om\vert_N)$,
where $\hat\om$ is the canonical symplectic structure on~$T^*N$.
\label{cr4thm2}
\end{thm}

In particular, if $N$ is compact and Lagrangian in Theorem
\ref{cr4thm2} then making $U$ smaller we can suppose it is
a tubular neighbourhood, and then $U,\Phi$ satisfy the
conditions of Theorem \ref{cr4thm1}. The important point is
that in Theorem \ref{cr4thm1}, the subsets $L_x=\Phi\bigl(
\pi^{-1}(x)\bigr)$ form a smooth family of noncompact Lagrangian
submanifolds of $M$, and $L_x$ intersects $N$ transversely at $x$.
Theorem \ref{cr4thm2} says that any such family $\{L_x:x\in N\}$
locally comes from a unique Lagrangian neighbourhood map~$\Phi$.

The goal of this section is to extend Theorem \ref{cr4thm1}
to SL cones $C$ in $\C^m$ and to SL $m$-folds $X$ with conical
singularities in an almost Calabi--Yau $m$-fold $M$. As this
involves {\it noncompact\/} Lagrangian $m$-folds $C',X'$, we
need to impose suitable asymptotic conditions on the Lagrangian
neighbourhood at the noncompact ends of $C',X'$. Throughout
we suppose~$m>2$.

\subsection{Dilation-equivariant neighbourhoods of cones}
\label{cr41}

We first extend Theorem \ref{cr4thm1} to SL cones in $\C^m$.
Most of the theorem is notation, not requiring proof. We
have to extend from a compact $N$ to the noncompact
$\Si\t(0,\iy)$, and include equivariance properties under
dilations on~$\C^m$.

\begin{thm} Let\/ $C$ be an SL cone in $\C^m$ with isolated
singularity at\/ $0$, and set\/ $\Si=C\cap{\cal S}^{2m-1}$.
Define $\iota:\Si\t(0,\iy)\ra\C^m$ by $\iota(\si,r)=r\si$,
with image $C\sm\{0\}$. For $\si\in\Si$, $\tau\in T_\si^*\Si$,
$r\in(0,\iy)$ and\/ $u\in\R$, let\/ $(\si,r,\tau,u)$ represent
the point\/ $\tau+u\,\d r$ in $T^*_{(\si,r)}\bigl(\Si\!\t\!
(0,\iy)\bigr)$. Identify $\Si\t(0,\iy)$ with the zero section
$\tau\!=\!u\!=\!0$ in $T^*\bigl(\Si\t(0,\iy)\bigr)$. Define an
action of\/ $(0,\iy)$ on $T^*\bigl(\Si\!\t\!(0,\iy)\bigr)$ by
\e
t:(\si,r,\tau,u)\longmapsto (\si,tr,t^2\tau,tu)
\quad\text{for $t\in(0,\iy)$,}
\label{cr4eq1}
\e
so that\/ $t^*(\hat\om)\!=\!t^2\hat\om$, for $\hat\om$ the
canonical symplectic structure on~$T^*\bigl(\Si\!\t\!(0,\iy)\bigr)$.

Then there exists an open neighbourhood\/ $U_{\sst C}$ of\/
$\Si\t(0,\iy)$ in $T^*\bigl(\Si\t(0,\iy)\bigr)$ invariant under
\eq{cr4eq1} given by
\e
U_{\sst C}=\bigl\{(\si,r,\tau,u)\in T^*\bigl(\Si\t(0,\iy)\bigr):
\bmd{(\tau,u)}<2\ze r\bigr\}\quad\text{for some $\ze>0$,}
\label{cr4eq2}
\e
where $\md{\,.\,}$ is calculated using the cone metric $\iota^*(g')$
on $\Si\t(0,\iy)$, and an embedding $\Phi_{\sst C}:U_{\sst C}\ra\C^m$
with\/ $\Phi_{\sst C}\vert_{\Si\t(0,\iy)}=\iota$, $\Phi_{\sst
C}^*(\om')=\hat\om$ and\/ $\Phi_{\sst C}\circ t=t\,\Phi_{\sst C}$
for all\/ $t>0$, where $t$ acts on $U_{\sst C}$ as in \eq{cr4eq1}
and on $\C^m$ by multiplication.
\label{cr4thm3}
\end{thm}

\begin{proof} For each $(\si,r)\in\Si\t(0,\iy)$, define
$L_{(\si,r)}$ to be the unique affine subspace $\R^m$ in $\C^m$
passing through $r\si$ and normal to $C$ there. Then $L_{(\si,r)}$
is a {\it Lagrangian plane} in $\C^m$, as $C$ is Lagrangian.
This defines a family $\bigl\{L_{(\si,r)}:(\si,r)\in\Si
\t(0,\iy)\bigr\}$ of Lagrangian submanifolds of $\C^m$
with $r\si\in L_{(\si,r)}$ and $T_{r\si}L_{(\si,r)}\cap
T_{r\si}C'=\{0\}$. We can therefore apply Theorem~\ref{cr4thm2}.

We have defined dilation actions of $\R_+$ on $T^*\bigl(\Si\t
(0,\iy)\bigr)$ and $\C^m$, and it is easy to see that we may
choose $U$ to be {\it dilation-invariant}, and then $\Phi$ is
{\it dilation-equivariant}, in the sense that $\Phi\circ t=t\,\Phi$.
It remains to show that we can take $U$ to be $U_{\sst C}$ in
\eq{cr4eq2} for some $\ze>0$. This is true if $U_{\sst C}\subset
U$. As $U,U_{\sst C}$ are both dilation-invariant, it is enough
for $U_{\sst C}\subset U$ to hold on the hypersurface $r=1$, that
is, over the compact subset $\Si\t\{1\}$. The existence of some
small $\ze>0$ with $U_{\sst C}\subset U$ then follows by compactness.
\end{proof}

Theorem \ref{cr4thm3} can also be proved by applying the {\it
Legendrian Neighbourhood Theorem} to $\Si$ in ${\cal S}^{2m-1}$.
This is the analogue of Theorem \ref{cr4thm1} for Legendrian
submanifolds in contact manifolds, and is described briefly
in~\cite[p.~107]{McSa}.

\subsection{Distinguished coordinates on $X'$ near $x_i$}
\label{cr42}

We shall use Theorem \ref{cr4thm3} to construct a particular
choice of $\phi_i$ in Definition~\ref{cr3def5}.

\begin{thm} Let\/ $(M,J,\om,\Om)$, $\psi,X,n,x_i,\up_i,C_i,\Si_i,
\mu_i,R,\Up_i$ and\/ $\iota_i$ be as in Definition \ref{cr3def5}.
Theorem \ref{cr4thm3} gives $\ze>0$, neighbourhoods $U_{\sst C_i}$
of\/ $\Si_i\t(0,\iy)$ in $T^*\bigl(\Si_i\t(0,\iy)\bigr)$ and
embeddings $\Phi_{\sst C_i}:U_{\sst C_i}\ra\C^m$ for~$i=1,\ldots,n$.

Then for sufficiently small\/ $R'\in(0,R]$ there exist unique
closed\/ $1$-forms $\eta_i$ on $\Si_i\t(0,R')$ for $i=1,\ldots,n$
written $\eta_i(\si,r)=\eta_i^1(\si,r)+\eta_i^2(\si,r)\d r$ for
$\eta_i^1(\si,r)\in T_\si^*\Si_i$ and\/ $\eta_i^2(\si,r)\in\R$,
and satisfying $\md{\eta_i(\si,r)}<\ze r$ and
\e
\bmd{\na^k\eta_i}=O(r^{\mu_i-1-k})
\quad\text{as $r\ra 0$ for $k=0,1,$}
\label{cr4eq3}
\e
computing $\na,\md{\,.\,}$ using the cone metric $\iota_i^*(g')$,
such that the following holds.

Define $\phi_i:\Si_i\t(0,R')\ra B_R$ by $\phi_i(\si,r)=\Phi_{\sst
C_i}\bigl(\si,r,\eta_i^1(\si,r),\eta_i^2(\si,r)\bigr)$. Then
$\Up_i\circ\phi_i:\Si_i\t(0,R')\ra M$ is a diffeomorphism
$\Si_i\t(0,R')\ra S_i$ for open sets $S_1,\ldots,S_n$ in $X'$
with\/ $\bar S_1,\ldots,\bar S_n$ disjoint, and\/ $K=X'\sm(S_1
\cup\cdots\cup S_n)$ is compact. Also $\phi_i$ satisfies
\eq{cr3eq5}, so that\/ $R',\phi_i,S_i,K$ satisfy
Definition~\ref{cr3def5}.
\label{cr4thm4}
\end{thm}

\begin{proof} As $X$ has a conical singularity at $x_i$ it
follows from \eq{cr3eq5} that near 0 in $B_R$ we can write
$\Up_i^*(X')$ as the image under $\Phi_{\sst C_i}$ of the graph of
a smooth 1-form $\eta_i$ on $\Si_i\t(0,R')$ for small $R'\in(0,R]$.
This just means that $\Up_i^*(X')$ intersects the Lagrangian ball
$\Phi_{\sst C_i}\bigl(T^*_{(\si,r)}(\Si_i\t(0,\iy))\cap
U_{\sst C_i}\bigr)$ transversely in exactly one point for
$(\si,r)\in\Si_i\t(0,R')$, and we define $\eta_i$ such that
this point is $\Phi_{\sst C_i}\bigl(\eta_i(\si,r)\bigr)$.
Since $\om\vert_{X'}\equiv 0$ and $\Up_i^*(\om)=\om'$,
$\Phi_{\sst C_i}^*(\om')=\hat\om$ we see that $\hat\om$
restricted to the graph of $\eta_i$ in $T^*\bigl(\Si_i\t
(0,R')\bigr)$ is zero. By a well-known fact in symplectic
geometry, this implies that $\eta_i$ is {\it closed}.

Now define $\phi_i:\Si_i\t(0,R')\ra B_R$ by $\phi_i(\si,r)=
\Phi_{\sst C_i}\bigl(\si,r,\eta_i^1(\si,r),\eta_i^2(\si,r)\bigr)$.
Then $\phi_i$ is an embedding, and by definition $\Up_i\circ\phi_i$
maps $\Si_i\t(0,R')\ra X'$. Define $S_i=\Up_i\circ\phi_i\bigl(\Si_i
\t(0,R')\bigr)$ and $K=X\sm(S_1\cup\cdots\cup S_n)$. Making $R'$
smaller if necessary we can arrange that $\bar S_1,\ldots,\bar S_n$
are disjoint. Then $\Up_i\circ\phi_i$ is a diffeomorphism
$\Si_i\t(0,R')\ra S_i$, and $S_i$ is an open set in $X'$, and
$K$ is the complement of open neighbourhoods of $x_1,\ldots,x_n$
in the compact space $X$, so $K$ is compact.

We have not yet shown that $\phi_1,\ldots,\phi_n$ satisfy
\eq{cr3eq5}. By Definition \ref{cr3def5} there must exist
some $\phi_1',\ldots,\phi_n'$ satisfying the conditions,
including \eq{cr3eq5}. Then $\phi_1,\ldots,\phi_n$ are
obtained from $\phi_1',\ldots,\phi_n'$ by a kind of
{\it projection}. What happens is that $\phi_i(\si,r)=
\phi_i'(\si',r')$ if $(\si,r)$, $(\si',r')$ are close
in $\Si_i\t(0,R')$ and $\phi_i'(\si',r')$ lies in the
affine normal subspace to $C_i$ at~$(\si,r)$.

For small $R''\in(0,R']$ define $\Xi_i:\Si_i\t(0,R'')\ra
\Si_i(0,R')$ by $\Xi_i(\si',r')=(\si,r)$. Then \eq{cr3eq5}
for $\phi_i'$ implies that
\e
\na^k(\Xi_i-\id)=O\bigl((r')^{\mu_i-1-k}\bigr)
\quad\text{as $r'\ra 0$ for $k=0,1$.}
\label{cr4eq4}
\e
But $\phi_i=\phi_i'\circ\Xi_i^{-1}$ for small $r$, so
combining \eq{cr4eq4} and \eq{cr3eq5} for $\phi_i'$
implies \eq{cr3eq5} for~$\phi_i$.

Equation \eq{cr3eq5} for $\phi_i$ and properties of
$\Phi_{\sst C_i}$ easily imply \eq{cr4eq3}. Finally,
as $\mu_i>2$ by Definition \ref{cr3def5}, equation
\eq{cr4eq3} implies that $\md{\eta_i}=o(r)$ for small
$r$. Therefore, making $R'$ smaller if necessary, we
can suppose that $\md{\eta_i(\si,r)}<\ze r$ for
$(\si,r)\in\Si_i\t(0,R')$.
\end{proof}

We can integrate the 1-forms $\eta_i$ in Theorem~\ref{cr4thm4}.

\begin{lem} In Theorem \ref{cr4thm4} we have $\eta_i=\d A_i$
for $i=1,\ldots,n$, where $A_i:\Si_i\t(0,R')\ra\R$ is given by
$A_i(\si,r)=\int_0^r\eta_i^2(\si,s)\d s$ and satisfies
\e
\bmd{\na^kA_i}=O(r^{\mu_i-k})
\quad\text{as $r\ra 0$ for $k=0,1,2,$}
\label{cr4eq5}
\e
computing $\na$ and\/ $\md{\,.\,}$ using the cone
metric~$\iota^*(g')$.
\label{cr4lem}
\end{lem}

\begin{proof} From \eq{cr4eq3} we deduce that $\md{\na^k\eta_i^2}
=O(r^{\mu_i-1-k})$ as $r\ra 0$ for $k=0,1$.  Integrating this and
using $\mu_i>2$ shows that $A_i(\si,r)=\int_0^r\eta_i^2(\si,s)\d s$
is well-defined and \eq{cr4eq5} holds for $k=0,1$. The $\d r$
component in $\d A_i$ is $\eta_i^2$, so that $\eta_i-\d A_i$ is a
closed 1-form on $\Si_i\t(0,R')$ with no $\d r$ component, and is
therefore independent of $r$. But \eq{cr4eq3} for $k=0$ and
\eq{cr4eq5} for $k=1$ imply that $\eta_i-\d A_i=O(r^{\mu_i-1})$ in the
cone metric on $\Si_i\t(0,R')$, so $\eta_i-\d A_i=O(r^{\mu_i-2})$
in the cylinder metric, and taking the limit $r\ra 0$ gives
$\eta_i-\d A_i=0$ as $\mu_i>2$. Hence $\eta_i=\d A_i$, and
\eq{cr4eq3} for $k=1$ then yields \eq{cr4eq5} for~$k=2$.
\end{proof}

\subsection{A Lagrangian Neighbourhood Theorem for $X'$}
\label{cr43}

Here is an analogue of Theorem \ref{cr4thm1} for SL $m$-folds
$X$ with conical singularities. We construct a Lagrangian
neighbourhood of $X'$ compatible with the distinguished
coordinates of Theorem \ref{cr4thm4}. The theorem will be
an important tool in \cite{Joyc3,Joyc4,Joyc5}, where we study
{\it deformations} and {\it desingularizations} of~$X$.

\begin{thm} Suppose $(M,J,\om,\Om)$ is an almost Calabi--Yau
$m$-fold and\/ $X$ a compact SL\/ $m$-fold in $M$ with conical
singularities at\/ $x_1,\ldots,x_n$. Let the notation $\psi,\up_i,
C_i,\Si_i,\mu_i,R,\Up_i$ and\/ $\iota_i$ be as in Definition
\ref{cr3def5}, and let\/ $\ze,U_{\sst C_i},\allowbreak
\Phi_{\sst C_i},\allowbreak R',\allowbreak \eta_i,\allowbreak
\eta_i^1,\eta_i^2,\phi_i,S_i$ and\/ $K$ be as in Theorem~\ref{cr4thm4}.

Then making $R'$ smaller if necessary, there exists an open tubular
neighbourhood\/ $U_\sXp\subset T^*X'$ of the zero section
$X'$ in $T^*X'$, such that under $\d(\Up_i\circ\phi_i):T^*\bigl(
\Si_i\t(0,R')\bigr)\ra T^*X'$ for $i=1,\ldots,n$ we have
\e
\bigl(\d(\Up_i\circ\phi_i)\bigr)^*(U_\sXp)=\bigl\{(\si,r,\tau,u)
\in T^*\bigl(\Si_i\t(0,R')\bigr):\bmd{(\tau,u)}<\ze r\bigr\},
\label{cr4eq6}
\e
and there exists an embedding $\Phi_\sXp:U_\sXp\ra M$ with\/
$\Phi_\sXp\vert_{X'}=\id:X'\ra X'$ and\/ $\Phi_\sXp^*(\om)
=\hat\om$, where $\hat\om$ is the canonical symplectic structure on
$T^*X'$, such that 
\e
\Phi_\sXp\circ\d(\Up_i\circ\phi_i)(\si,r,\tau,u)\equiv\Up_i\circ
\Phi_{\sst C_i}\bigl(\si,r,\tau+\eta_i^1(\si,r),u+\eta_i^2(\si,r)\bigr)
\label{cr4eq7}
\e
for all\/ $i=1,\ldots,n$ and\/ $(\si,r,\tau,u)\in T^*\bigl(\Si_i\t(0,R')
\bigr)$ with\/ $\bmd{(\tau,u)}<\ze r$. Here $\md{(\tau,u)}$ is computed
using the cone metric $\iota_i^*(g')$ on~$\Si_i\t(0,R')$.
\label{cr4thm5}
\end{thm}

\begin{proof} Let us regard \eq{cr4eq6} and \eq{cr4eq7} as
{\it definitions} of $U_\sXp$ and $\Phi_\sXp$ over
the subset $S_i$ of $X'$ for $i=1,\ldots,n$. Since
$\md{\eta_i(\si,r)}<\ze r$ by Theorem \ref{cr4thm4},
$\md{(\tau,u)}<\ze r$ in \eq{cr4eq6} and $\Phi_{\sst C_i}
(\si,r,\tau',u')$ is defined provided $\md{(\tau',u')}<2\ze r$
by \eq{cr4eq2}, we see that $\Phi_{\sst C_i}(\ldots)$ is
well-defined in~\eq{cr4eq7}.

Making $R'$ smaller if necessary, we can ensure that $\Phi_{\sst
C_i}(\ldots)$ lies in $B_R$, and so \eq{cr4eq7} makes sense and
$U_\sXp,\Phi_\sXp$ are well-defined over $S_i$. As
$\Up_i^*(\om)=\om'$, $\Phi_{\sst C_i}^*(\om')=\hat\om$ and $\eta_i$
is closed, it easily follows that $\Phi_\sXp^*(\om)=\hat\om$
on these regions of $U_\sXp$. Also $\Phi_\sXp$ is an
embedding on these regions, as $\Up_i$ and $\Phi_{\sst C_i}$ are,
and is the identity on each $S_i$, by definition of $\phi_i$ in
Theorem \ref{cr4thm4}. It remains to extend $U_\sXp$ and
$\Phi_\sXp$ over the compact subset $K$ in~$X'$.

For $x\in S_i$ define $L_x=\Phi_\sXp(T^*_xX'\cap U_\sXp)$,
where $U_\sXp$, $\Phi_\sXp$ are defined over $S_i$ as above.
As $\Phi_\sXp$ is an embedding with $\Phi_\sXp^*(\om)=\hat\om$
we see that $L_x$ is an open Lagrangian ball in $M$ which meets $X'$
transversely at $x$, and depends smoothly on $x$. Extend this
family $\{L_x:x\in S_i$, $i=1,\ldots,n\}$ to a family
$\{L_x:x\in X'\}$ such that $L_x$ is an open Lagrangian ball in
$M$ which meets $X'$ transversely at $x$, and depends smoothly
on $x$. This is possible by standard symplectic geometry
techniques, as the extension is over a compact set~$K$.

Now apply Theorem \ref{cr4thm2} to the family $\{L_x:x\in X'\}$.
This gives an open neighbourhood $U$ of $X'$ in $T^*X'$, and a map
$\Phi:U\ra M$ with $\Phi\vert_{X'}=\id:X'\ra X'$ and $\Phi^*(\om)
=\hat\om$. By the local uniqueness of $\Phi$ in Theorem
\ref{cr4thm2} we see that $\Phi$ and $\Phi_\sXp$ defined
above coincide where they are both defined. 

Therefore we can take $U$ to be $U_\sXp$ and $\Phi$ to be
$\Phi_\sXp$ as defined above over $S_i$, for $i=1,\ldots,n$.
Choose an open tubular neighbourhood $U_\sXp$ of $X'$ in $U$,
which coincides with the previous definition of $U_\sXp$ over
$S_i$. This is possible as $U$ is open and it only remains to choose
$U_\sXp$ over the compact set $K$. Let $\Phi_\sXp$ be the
restriction of $\Phi$ to $U_\sXp\subseteq U$. Then $U_\sXp,
\Phi_\sXp$ satisfy all the conditions of the theorem.
\end{proof}

\subsection{Extending to families of almost Calabi--Yau $m$-folds}
\label{cr44}

In \cite{Joyc3,Joyc5} we will study SL $m$-folds not just
in one almost Calabi--Yau $m$-fold $(M,J,\om,\Om)$, but in
a smooth family of them.

\begin{dfn} Let $(M,J,\om,\Om)$ be an almost Calabi--Yau
$m$-fold. A {\it smooth family of deformations of\/}
$(M,J,\om,\Om)$ is a connected open set $\F\subset\R^d$
for $d\ge 0$ with $0\in\F$ called the {\it base space},
and a smooth family $\bigl\{(M,J^s,\om^s,\Om^s):s\in{\cal
F}\bigr\}$ of almost Calabi--Yau structures on $M$
with~$(J^0,\om^0,\Om^0)=(J,\om,\Om)$.
\label{cr4def}
\end{dfn}

We now extend the Lagrangian neighbourhood of an SL
$m$-fold $X$ with conical singularities in $(M,J,\om,\Om)$
constructed in Theorem \ref{cr4thm5} to a smooth family of
similar neighbourhoods of $X$ in $(M,J^s,\om^s,\Om^s)$ for
small $s$. If $\om^s\vert_{X'}$ is not exact then we cannot
deform $X'$ to a Lagrangian $m$-fold in $(M,\om^s)$.
Therefore we replace the condition $\Phi_\sXp^*(\om)
=\hat\om$ in Theorem \ref{cr4thm5} by $(\Phi^s_\sXp)^*
(\om^s)=\hat\om+\pi^*(\nu^s)$, where $\nu^s$ is a
compactly-supported closed 2-form on~$X'$.

\begin{thm} Let\/ $(M,J,\om,\Om)$ be an almost Calabi--Yau
$m$-fold and\/ $X$ a compact SL\/ $m$-fold in $M$ with conical
singularities at $x_1,\ldots,x_n$. Let the notation
$R,\Up_i,\ze,\Phi_{\sst C_i},R',\eta_i,\eta_i^1,\eta_i^2,
\phi_i,S_i,K$ be as in Theorem \ref{cr4thm4}, and let\/
$U_\sXp,\Phi_\sXp$ be as in Theorem \ref{cr4thm5}.
Suppose $\bigl\{(M,J^s,\om^s,\Om^s):s\in \F\bigr\}$ is
a smooth family of deformations of $(M,J,\om,\Om)$ with base
space $\F\subset\R^d$. Define $\psi^s:M\ra(0,\iy)$ for
$s\in\F$ as in \eq{cr3eq3}, but using~$\om^s,\Om^s$.

Then making $R,R'$ and\/ $U_\sXp$ smaller if necessary,
for some connected open $\F'\subseteq\F$ with\/
$0\in\F'$ and all\/ $s\in \F'$ there exist
\begin{itemize}
\item[{\rm(a)}] isomorphisms $\up_i^s:\C^m\ra T_{x_i}M$ for
$i=1,\ldots,n$ with\/ $\up_i^0=\up_i$, $(\up_i^s)^*(\om^s)=\om'$
and\/~$(\up_i^s)^*(\Om)=\psi^s(x_i)^m\Om'$,
\item[{\rm(b)}] embeddings $\Up_i^s:B_R\ra M$ for $i=1,\ldots,n$
with\/ $\Up_i^0=\Up_i$, $\Up_i^s(0)=x_i$, $\d\Up_i^s\vert_0=\up_i^s$
and\/~$(\Up_i^s)^*(\om^s)=\om'$,
\item[{\rm(c)}] a closed\/ $2$-form $\nu^s\in C^\iy(\La^2T^*X')$
supported in $K\subset X'$ with\/ $\nu^0=0$, and
\item[{\rm(d)}] an embedding $\Phi^s_\sXp\!:\!U_\sXp\!\ra\!M$
with\/ $\Phi_\sXp^0\!=\!\Phi_\sXp$ and\/~$(\Phi^s_\sXp
)^*(\om^s)\!=\!\hat\om\!+\!\pi^*(\nu^s)$,
\end{itemize}
all depending smoothly on $s\in\F'$ with
\e
\Phi_\sXp^s\circ\d(\Up_i\circ\phi_i)(\si,r,\tau,u)\equiv\Up_i^s\circ
\Phi_{\sst C_i}\bigl(\si,r,\tau+\eta_i^1(\si,r),u+\eta_i^2(\si,r)\bigr)
\label{cr4eq8}
\e
for all\/ $s\in \F'$, $i=1,\ldots,n$ and\/ $(\si,r,\tau,u)\in
T^*\bigl(\Si_i\t(0,R')\bigr)$ with\/~$\bmd{(\tau,u)}<\ze r$.
\label{cr4thm6}
\end{thm}

\begin{proof} Clearly, for some open neighbourhood $\F_1$ of 0
in $\F$ we can extend $\up_i$ to a smooth family $\up_i^s$ for
$s\in\F_1$ satisfying part (a). By a version of Darboux's
Theorem \cite[Th.~3.15]{McSa} for smooth families of symplectic
manifolds, making $R$ smaller if necessary, for some open
neighbourhood $\F_2$ of 0 in $\F_1$ we can extend $\Up_i:B_R\ra M$
to a smooth family of embeddings $\Up_i^s:B_R\ra M$ for $s\in\F_2$
satisfying (b). We then make $R'$ smaller if necessary so that
Theorem \ref{cr4thm5} holds.

Next, for some open neighbourhood $\F_3$ of 0 in $\F_2$,
choose a smooth family of embeddings $\chi^s:X'\ra M$ for
$s\in\F_2$ with $\chi^0=\id:X'\ra X'\subset M$ such that
\e
\chi^s\circ\Up_i\circ\phi_i\equiv\Up_i^s\circ\phi_i
\quad\text{on $\Si_i\t(0,R')$ for $i=1,\ldots,n$ and $s\in\F_3$.}
\label{cr4eq9}
\e
That is, we define $\chi^s$ to be $\Up_i^s\circ\Up_i^{-1}$ on
$S_i$ for $i=1,\ldots,n$, and then extend $\chi^s$ smoothly to
an embedding on $K$, the rest of $X'$. This is possible for
$s$ near 0, as $K$ is compact.

Now define $\nu^s=(\chi^s)^*(\om^s)\in C^\iy(\La^2T^*X')$ for
$s\in\F_3$. As $\chi^s,\om^s$ depend smoothly on $s$
so does $\nu^s$, and as $\chi^0=\id$ and $\om^0\vert_{X'}=
\om\vert_{X'}=0$ we have $\nu^0=0$. Also, as $(\Up_i^s)^*(\om^s)
=\om'$ we see from \eq{cr4eq9} that $\nu^s=(\chi^s)^*(\om^s)$
is independent of $s$ on $S_i$, the image of $\Up_i\circ\phi_i$,
so that $\nu^s=\nu^0=0$ on $S_i$, and $\nu^s$ is supported on
$K=X'\sm(S_1\cup\cdots\cup S_n)$. This gives part~(c).

Define $\Phi_\sXp^0=\Phi_\sXp$. As in the proof of
Theorem \ref{cr4thm5}, regard \eq{cr4eq8} as a {\it definition}
of $\Phi_\sXp^s$ on $\pi^*(S_i)\subset U_\sXp$ for
$s\in\F_3$. This is well-defined, and depends smoothly on
$s$. Since $(\Up_i^s)^*(\om^s)=\om'$ is independent of $s$, we
see from \eq{cr4eq8} that $(\Phi_\sXp^s)^*(\om^s)$ is
independent of $s$ on $\pi^*(S_i)\subset U_\sXp$. Thus
on $\pi^*(S_i)\subset U_\sXp$ we have
\e
(\Phi_\sXp^s)^*(\om^s)=(\Phi_\sXp^0)^*(\om^0)
=\Phi_\sXp^*(\om)=\hat\om=\hat\om+\pi^*(\nu^s),
\label{cr4eq10}
\e
since $\nu^s\equiv 0$ on $S_i$. It only remains to extend
$\Phi_\sXp^s$ over $\pi^*(K)\subset U_\sXp$ for~$s\ne 0$.

Generalizing the proof of Theorem \ref{cr4thm5}, for $x\in X'$ and
$s=0$, or for $x\in S_i$ and $s\in\F_3$, define $L_x^s=
\Phi_\sXp^s(T^*_xX'\cap U_\sXp)$, where $\Phi_\sXp^s$
is defined in these regions as above. Since $(\Phi_\sXp^s)^*
(\om^s)=\hat\om+\pi^*(\nu^s)$ wherever $\Phi_\sXp$ is defined
and $\hat\om+\pi^*(\nu^s)$ vanishes on $T_x^*X$, we see that $L_x^s$
is an open Lagrangian ball in $(M,\om^s)$ which meets $\chi^s(X')$
transversely at $\chi^s(x)$, and depends smoothly on~$x,s$.

Extend this family $\{L_x^0:x\in X'\}\cup\{L_x^s:x\in S_i$,
$i=1,\ldots,n$, $s\in\F_3\}$ to a family $\{L_x^s:x\in X'$,
$s\in\F_4\}$ for some open neighbourhood $\F_4$ of 0
in $\F_3$, such that $L^s_x$ is an open Lagrangian ball in
$(M,\om^s)$ containing $\chi^s(x)$, which meets $\chi^s(X')$
transversely at $\chi^s(x)$, and depends smoothly on $x,s$.
This is possible by standard symplectic geometry techniques,
as the extension is over $x$ in a compact set $K$ and for
small~$s$.

Now apply Theorem \ref{cr4thm2} to the family $\{L_x^s:x\in X'\}$
in $(M,\om^s)$ for $s\in\F_4$, replacing $X$ by $\chi^s(X')$.
Arguing as Theorem \ref{cr4thm5}, we get a tubular neighbourhood
$U_\sXp^s$ of $X'$ in $T^*X'$ with $U_\sXp^s\cap\pi^*(S_i)
=U_\sXp\cap\pi^*(S_i)$ and an embedding $\Phi_\sXp^s:
U_\sXp^s\ra M$ with $\Phi_\sXp^s\vert_{X'}=\chi^s$,
satisfying \eq{cr4eq8}. The formula $\Phi^*(\om)=\hat\om+\pi^*
(\om\vert_N)$ in Theorem \ref{cr4thm2} yields $(\Phi_\sXp^s)^*
(\om^s)=\hat\om+\pi^*(\nu^s)$, as we have to prove, since~$\om^s
\vert_{\chi^s(X')}=(\chi^s)^*(\om^s)=\nu^s$.

As the $L_x^s,\chi^s$ and $\om^s$ depend smoothly on $s$, so does
$\Phi_\sXp^s$. It remains to show that we may take the domain
$U_\sXp^s\subset T^*X'$ of $\Phi_\sXp^s$ to be $U_\sXp$,
independent of $s$. We can achieve this by making $U_\sXp$
smaller if necessary, though keeping it defined by \eq{cr4eq6} over
$S_i$, and restricting $s$ to a small connected open neighbourhood
$\F'$ of 0 in $\F_4$. This completes the proof.
\end{proof}

In the notation of \S\ref{cr24}, the 2-forms $\nu^s$ in Theorem
\ref{cr4thm6} define classes $[\nu^s]$ in $H^2_{\rm cs}(X',\R)$.
We investigate these classes, and the freedom to choose~$\nu^s$.

\begin{thm} In the situation of Theorem \ref{cr4thm6}, under
the isomorphism \eq{cr2eq18}, the class $[\nu^s]\in H^2_{\rm cs}
(X',\R)$ is identified with the map $H_2(X,\R)\ra\R$ given by
$\ga\mapsto\iota_*(\ga)\cdot[\om^s]$, where $\iota:X\ra M$ is
the inclusion, $\iota_*:H_2(X,\R)\ra H_2(M,\R)$ the induced map,
and\/ $[\om^s]\in H^2(M,\R)$. Thus $[\nu^s]$ depends only on
$X,M$ and\/~$[\om^s]\in H^2(M,\R)$.

Let\/ $V\cong H^2_{\rm cs}(X',\R)$ be a vector space of smooth
closed\/ $2$-forms on $X'$ supported in $K$ representing
$H^2_{\rm cs}(X',\R)$. Then making $\F'$ smaller if
necessary, we can choose $\Up_i^s,\nu^s$ and\/ $\Phi_\sXp^s$ in
Theorem \ref{cr4thm6} so that\/ $\nu^s\in V$ for all\/ $s\in\F'$.
In particular, if\/ $[\nu^s]=0$ in $H^2_{\rm cs}(X',\R)$ then we can
choose~$\nu^s=0$.
\label{cr4thm7}
\end{thm}

\begin{proof} As $\nu^s=(\chi^s)^*(\om^s)$, for $\ga\in H_2(X,\R)$
we have
\begin{equation*}
\ga\cdot[\nu^s]=\chi^s_*(\ga)\cdot[\om^s]=\iota_*(\ga)\cdot[\om^s],
\end{equation*}
since $\chi^s,\iota:X\ra M$ are isotopic as $\F'$ is connected,
and so $\chi^s_*(\ga)=\iota_*(\ga)$. This proves the first part,
and thus $[\nu^s]$ depends only on $X,M$ and~$[\om^s]$.

Now let $\F',\Up_i^s,\chi^s,\nu^s,\Phi_\sXp^s$ be
as in Theorem \ref{cr4thm6}, and let $V$ be a vector space of
closed $2$-forms on $X'$ supported in $K$ representing
$H^2_{\rm cs}(X',\R)$. We shall show that making $\F'$
smaller if necessary, we can modify $\chi^s,\nu^s,\Phi_\sXp^s$
to alternative choices $\ti\chi^s,\ti\nu^s,\ti\Phi_\sXp^s$
with $\ti\nu^s\in V$, keeping the same choice of~$\Up_i^s$.

For each $s\in\F'$, let $\ti\nu^s$ be the unique element
of $V$ with $[\ti\nu^s]=[\nu^s]$ in $H^2_{\rm cs}(X',\R)$. Then
$\ti\nu^s$ depends smoothly on $s$, as $[\nu^s]$ does. As
$[\nu^s-\ti\nu^s]=0$ in $H^2_{\rm cs}(X',\R)$ there exist
compactly-supported 1-forms $\be^s$ on $X'$ with $\d\be^s=
\nu^s-\ti\nu^s$. Since $\ti\nu^s,\nu^s$ are supported in $K$ we
can choose $\be^s$ supported in $K$. We can also choose $\be^s$
to depend smoothly on $s\in\F'$, as $\ti\nu^s,\nu^s$ do,
and choose $\be^0=0$, as~$\ti\nu^0=\nu^0=0$.

As $\be^0=0$ and $\be^s$ depends smoothly on $s$ and is supported
in $K$, making $\F'$ smaller if necessary we can suppose
that the graph $\Ga(\be^s)$ of $\be^s$ lies in $U_\sXp\subset
T^*X'$ for all $s\in\F'$. Define $\ti\chi^s:X'\ra M$ by
$\ti\chi^s=\Phi_\sXp^s\circ\be^s$, regarding $\be^s$ as a
map $X'\ra\Ga(\be^s)\subset U_\sXp$ in the obvious way.
Then $\ti\chi^s$ depends smoothly on $s\in\F'$ as $\be^s,
\Phi_\sXp^s$ do, and $\ti\chi^0=\chi^0=\id$ as $\be^0=0$,
and $\ti\chi^s=\chi^s$ on $S_i$ as $\be^s=0$ on $S_i$, so
$\ti\chi^s$ satisfies~\eq{cr4eq9}.

As $(\Phi_\sXp^s)^*(\om^s)=\hat\om+\pi^*(\nu^s)$, we see that
$(\ti\chi^s)^*(\om^s)=(\be^s)^*\bigl(\hat\om+\pi^*(\nu^s)\bigr)$,
where $\be^s$ maps $X'\ra\Ga(\be^s)\subset U_\sXp$. But
$(\be^s)^*(\hat\om)=-\d\be^s$ by a well-known fact in symplectic
geometry, and $(\be^s)^*\bigl(\pi^*(\nu^s)\bigr)=\nu^s$ as
$\pi\circ\be^s=\id:X'\ra X'$, so
\begin{equation*}
(\ti\chi^s)^*(\om^s)=-\d\be^s+\nu^s=-(\nu^s-\ti\nu^s)+\nu^s=\ti\nu^s.
\end{equation*}
Thus in the proof of Theorem \ref{cr4thm6} we are free to choose
$\ti\chi^s$ instead of $\chi^s$ such that $(\ti\chi^s)^*(\om^s)=
\ti\nu^s$ lies in $V$. The rest of the proof of Theorem
\ref{cr4thm6} then shows that we can choose $\ti\Phi_\sXp^s$
consistently with $\ti\chi^s,\ti\nu^s$. Finally, if $[\nu^s]=0$ in
$H^2_{\rm cs}(X',\R)$ then $\ti\nu^s=0$ in $V$, so we can
choose~$\nu^s=0$.
\end{proof}

\section{The asymptotic behaviour of $X$ near $x_i$}
\label{cr5}

We shall now show that the asymptotic condition \eq{cr3eq5}
in Definition \ref{cr3def5} can be strengthened in two ways:
we can make \eq{cr3eq5} hold for all $k\ge 0$ rather than
just $k=0,1$, and we can improve the asymptotic decay rates
$\mu_i$. On the way we will prove that compact SL $m$-folds
$X$ with conical singularities are automatically Riemannian
manifolds with conical singularities in the sense of
\S\ref{cr2}, so we deduce some analytic and Hodge theoretic
results on $X'$. Throughout we suppose~$m>2$.

\subsection{Regularity of higher derivatives}
\label{cr51}

We shall use the special Lagrangian condition to show that
\eq{cr3eq5}, \eq{cr4eq3} and \eq{cr4eq5} hold for all $k\ge 0$
for the $\phi_i$ constructed in Theorem \ref{cr4thm4}. Note that
this is {\it not\/} true for arbitrary $\phi_1,\ldots,\phi_n$
satisfying Definition~\ref{cr3def5}.

\begin{thm} In the situation of Theorem \ref{cr4thm4} and
Lemma \ref{cr4lem} we have
\e
\begin{gathered}
\bmd{\na^k(\phi_i-\iota_i)}=O(r^{\mu_i-1-k}),\quad
\bmd{\na^k\eta_i}=O(r^{\mu_i-1-k})\quad\text{and}\\
\bmd{\na^kA_i}=O(r^{\mu_i-k})
\quad\text{as $r\ra 0$ for all\/ $k\ge 0$ and\/ $i=1,\ldots,n$.}
\end{gathered}
\label{cr5eq1}
\e
Here $\na$ and\/ $\md{\,.\,}$ are computed using the cone
metric~$\iota_i^*(g')$.
\label{cr5thm1}
\end{thm}

\begin{proof} Let $\al$ be a smooth 1-form on $\Si_i\t(0,R')$
with $\md{\al(\si,r)}<\ze r$, written $\al(\si,r)=\al^1(\si,r)
+\al^2(\si,r)\d r$ for $\al^1(\si,r)\in T_\si^*\Si_i$ and
$\al^2(\si,r)\in\R$. Define a map $\Th_\al:\Si_i\t(0,R')\ra B_R$ by
$\Th_\al(\si,r)=\Phi_{\sst C_i}\bigl(\si,r,\al^1(\si,r),\al^2(\si,r)
\bigr)$. Define a smooth real function $F_i(\al)$ on
$\Si_i\t(0,R')$ by
\e
F_i(\al)\,\d V=\psi(x_i)^{-m}(\Up_i\circ\Th_\al)^*(\Im\Om),
\label{cr5eq2}
\e
where $\d V$ is the volume form of $\iota_i^*(g')$ on $\Si_i\t(0,R')$.
This defines a function $F_i$ from smooth 1-forms $\al$ on $\Si_i\t
(0,R')$ with $\md{\al}<\ze r$ to smooth functions on~$\Si_i\t(0,R')$.

The value of $F_i(\al)$ at $(\si,r)\in\Si_i\t(0,R')$ depends on
$\Up_i\circ\Th_\al$ and $\d(\Up_i\circ\Th_\al)$ at $(\si,r)$, which
depend on both $\al$ and $\na\al$ at $(\si,r)$. Hence $F_i(\al)$
depends pointwise on both $\al$ and $\na\al$, rather than just
$\al$. Define a map
\begin{gather}
\begin{split}
Q_i:\bigl\{(\si,r,y,z):\,&(\si,r)\in\Si_i\t(0,R'),\quad
y\in T^*_{(\si,r)}\bigl(\Si_i\t(0,R')\bigr),\\
&\md{y}<\ze r,\quad z\in\ot^2T^*_{(\si,r)}\bigl(\Si_i\t(0,R')
\bigr)\bigr\}\ra\R
\end{split}
\label{cr5eq3}\\
\text{by}\quad
Q_i\bigl(\si,r,\al(\si,r),\na\al(\si,r)\bigr)=
\bigl(\d^*\al+F_i(\al)\bigr)\bigl[(\si,r)\bigr]
\label{cr5eq4}
\end{gather}
for all 1-forms $\al$ on $\Si_i\t(0,R')$ with $\md{\al(\si,r)}<\ze r$
when $(\si,r)\in\Si_i\t(0,R')$. This is well-defined, as $F_i$ has
the right pointwise dependence in \eq{cr5eq4}, and the 1-forms $\al$
sweep out the domain of $Q_i$ in~\eq{cr5eq3}.

Let $\Om'$ be as in \eq{cr3eq1}, and rewrite \eq{cr5eq2} as
\e
F_i(\al)\,\d V=\Th_\al^*\bigl(\psi(x_i)^{-m}\Up_i^*(\Im\Om)
-\Im\Om'\bigr)+\Th_\al^*(\Im\Om').
\label{cr5eq5}
\e
As $\Up_i^*(\Im\Om)=\psi(x_i)^m\Im\Om'$ at $0\in B_R$ by Definition
\ref{cr3def5} and $\Up_i^*(\Im\Om)$ is smooth we see that
\e
\psi(x_i)^{-m}\Up_i^*(\Im\Om)-\Im\Om'=O(r)\quad\text{on $B_R$.}
\label{cr5eq6}
\e
Now consider the map $\al\mapsto\Th_\al^*(\Im\Om')$. When $\al=0$
this is the pull-back of $\Im\Om'\vert_{C_i'}$, which is zero as
$C_i'$ is special Lagrangian. So~$\Th_0^*(\Im\Om')=0$.

Following \cite[p.~721-2]{McLe} or \cite[Prop.~2.10]{Joyc3}, we find
that the linearization of $\Th_\al^*(\Im\Om')$ in $\al$ at $\al=0$
is $-(\d^*\al)\,\d V$. Thus we see that
\e
\Th_\al^*(\Im\Om')=-(\d^*\al)\,\d V+O(r^{-2}\ms{\al}+\ms{\na\al})
\label{cr5eq7}
\e
for $r^{-1}\md{\al},\md{\na\al}$ small, using the
dilation-equivariance properties of $\Th_\al$ to determine the
powers of $r$ in $O(r^{-2}\ms{\al}+r^0\ms{\na\al})$. Combining
\eq{cr5eq4}--\eq{cr5eq7} gives
\e
Q_i(\si,r,y,z)=O(r+r^{-2}\ms{y}+\ms{z})
\quad\text{when $\md{y}=O(r)$ and $\md{z}=O(1)$.}
\label{cr5eq8}
\e

When $\al=\eta_i=\d A_i$ we have $\Th_\al=\phi_i$, so
$\Up_i\circ\Th_\al$ maps $\Si_i\t(0,R')\ra S_i\subset X'$. Thus
$(\Up_i\circ\Th_\al)^*(\Im\Om)=0$ as $X'$ is special Lagrangian,
and $F_i(\d A_i)=0$. Hence \eq{cr5eq4} gives
\e
\De A_i-Q_i\bigl(\si,r,\d A_i(\si,r),\na^2A_i(\si,r)\bigr)=0
\label{cr5eq9}
\e
for $(\si,r)\in\Si_i\t(0,R')$. This is a {\it second-order
nonlinear elliptic equation} on $A_i$. We shall use {\it elliptic
regularity results} for \eq{cr5eq9} to prove~\eq{cr5eq1}.

For $t\in(0,R']$ define $\de^t:\Si_i\t(\frac{1}{2},1)\ra\Si_i\t(0,R')$
by $\de^t(\si,r)=(\si,tr)$. Define $Q_i^t(\si,r,y,z)=t^{2-\mu_i}Q_i
\bigl(\si,tr,t^{\mu_i}\de^t_*(y),t^{\mu_i}\de^t_*(z)\bigr)$, where
\e
\begin{split}
Q_i^t:\bigl\{(\si,r,y,z):\,&(\si,r)\in\Si_i\t(\ts\frac{1}{2},1),\quad
y\in T^*_{(\si,r)}\bigl(\Si_i\t(\ts\frac{1}{2},1)\bigr),\\
&\md{y}<t^{2-\mu_i}\ze r,\quad z\in\ot^2T^*_{(\si,r)}
\bigl(\Si_i\t(\ts\frac{1}{2},1)\bigr)\bigr\}\ra\R.
\end{split}
\label{cr5eq10}
\e
Define functions $A_i^t:\Si_i\t(\frac{1}{2},1)\ra\R$ by
\e
A_i^t(\si,r)=t^{-\mu_i}A_i(\si,tr).
\label{cr5eq11}
\e
Then \eq{cr5eq9} implies that for $(\si,r)\in\Si_i\t(\ts\frac{1}{2},1)$
we have
\e
\De A_i^t-Q_i^t\bigl(\si,r,\d A_i^t(\si,r),\na^2A_i^t(\si,r)\bigr)=0.
\label{cr5eq12}
\e
Also \eq{cr4eq5} shows that for some $C>0$ independent of
$i,t$ we have
\e
\bmd{\na^kA_i^t}\le C
\quad\text{on $\Si_i\t(\ts\frac{1}{2},1)$ for $k=0,1,2$
and $t\in(0,R']$.}
\label{cr5eq13}
\e

From \eq{cr5eq8}, noting that $\bmd{\de^t_*(y)}=t^{-1}\md{y}$
and $\bmd{\de^t_*(z)}=t^{-2}\md{z}$, we find that
\begin{equation*}
Q_i^t(\si,r,y,z)=O(t^{3-\mu_i}+t^{\mu_i-2}\ms{y}+t^{\mu_i-2}\ms{z})
\end{equation*}
when $\md{y}=O(t^{2-\mu_i})$ and $\md{z}=O(t^{2-\mu_i})$.
Thus $Q_i^t\ra 0$ as $t\ra 0$ uniformly on compact subsets of
the domain in \eq{cr5eq10}, since $2<\mu_i<3$. Furthermore, one
can show that all derivatives of $Q_i^t$ converge to 0 uniformly
on compact subsets as $t\ra 0$. Therefore for small $t$, equation
\eq{cr5eq12} approximates the much simpler {\it linear} elliptic
equation~$\De A_i^t=0$.

Now Ivanov \cite{Ivan} studies nonlinear elliptic equations
$F(x,u,\d u,\na^2u)=0$ for $x$ in a bounded domain $S$ in $\R^n$
and $u\in C^4(S)$, where $F(x,u,v,w)$ is a smooth function of its
arguments. When $T\subset S^\circ$ is an interior domain,
$\md{\na^ku}\le C$ for $k=0,1,2$ and $F$ is {\it close to
quasilinear}, in the sense that the second derivatives of $F$
in the $w$ variables are small compared to other constants
depending on $S,T,C$ and the first and second derivatives of
$F$, he proves \cite[Th.~2.2]{Ivan} a priori interior estimates
for the H\"older $C^{k+2,\al}$ norm of $u$ on $T$, depending on
the same constants and the $C^{k,\al}$ norm of $F$ on a compact
subset of its domain.

This generalizes immediately to interior estimates on Riemannian
manifolds. Thus we can apply it to \eq{cr5eq12} with $S=\Si_i\t
(\frac{1}{2},1)$, $T=\Si_i\t(\ts\frac{2}{3},\frac{3}{4})$,
$u=A_i^t$ and $C$ as in \eq{cr5eq13}. For small $t$, say when
$t\le\ka$ for $\ka\in(0,R']$, equation \eq{cr5eq12} is `close
to quasilinear' in the appropriate sense, and Ivanov's result
applies {\it uniformly in} $t$. Hence there exist constants
$C_k>0$ for $k\ge 0$ such that
\e
\bmd{\na^kA_i^t}\le C_k
\quad\text{on $\Si_i\t(\ts\frac{2}{3},\frac{3}{4})$ for all
$k\ge 0$ and $t\in(0,\ka]$.}
\label{cr5eq14}
\e

Combining \eq{cr5eq11} and \eq{cr5eq14} proves that $\md{\na^kA_i}
=O(r^{\mu_i-k})$ as $r\ra 0$ for all $k\ge 0$. As $\eta_i=\d A_i$
by Lemma \ref{cr4lem}, it immediately follows that $\md{\na^k\eta_i}
=O(r^{\mu_i-1-k})$ as $r\ra 0$ for all $k\ge 0$. Finally $\md{\na^k(
\phi_i-\iota_i)}=O(r^{\mu_i-1-k})$ follows from relationship between
$\eta_i$ and $\phi_i$ in Theorem \ref{cr4thm4}, and the dilation
equivariance properties of $\Phi_{\sst C_i}$. This completes the proof.
\end{proof}

\subsection{Treating $X$ as a manifold with conical singularities}
\label{cr52}

From Theorem \ref{cr5thm1} it follows that $g$ on $X$ satisfies
\eq{cr2eq1} with $\nu_i=\mu_i-2$. Therefore SL $m$-folds with
conical singularities fit into the framework of~\S\ref{cr2}.

\begin{thm} Suppose $(M,J,\om,\Om)$ is an almost Calabi--Yau
$m$-fold and\/ $X$ a compact SL\/ $m$-fold in $M$ with conical
singularities at $x_1,\ldots,x_n$ with rates $\mu_1,\ldots,\mu_n>2$,
as in Definition \ref{cr3def5}. Then $X$ with the induced metric
$d$ is a Riemannian manifold with conical singularities in the
sense of Definition \ref{cr2def1}, with\/~$\nu_i=\mu_i-2>0$.
\label{cr5thm2}
\end{thm}

There are a few small notational differences between \S\ref{cr2}
and \S\ref{cr33}. For instance, $\phi_i$ in \S\ref{cr2} is replaced
by $\Up_i\circ\phi_i$ in \S\ref{cr33}, $\ep$ in \S\ref{cr2} is
replaced by $R'$ in \S\ref{cr33}, and $S_i$ is defined to be
$\{y\in X:0<d(x_i,y)<\ep\}$ in \S\ref{cr2} and the image of
$\Up_i\circ\phi_i$ in \S\ref{cr33}. These differences are all
entirely superficial, so we will ignore them.

We can now use the analysis of \S\ref{cr2} to prove elliptic
regularity results on $X'$. However, rather than studying the
Laplacian $\De$ on $X$ we consider the operator $P:f\mapsto
\d^*(\psi^m\d f)$, as this is what we will need in~\cite{Joyc3}.

\begin{thm} Let\/ $(M,J,\om,\Om)$ be an almost Calabi--Yau
$m$-fold, and define $\psi:M\ra(0,\iy)$ as in \eq{cr3eq3}.
Suppose $X$ is a compact SL\/ $m$-fold in $M$ with conical
singularities at $x_1,\ldots,x_n$ with cones $C_i$ and rates
$\mu_i$. Define the Banach spaces $L^p_{k,{\bs\be}}(X')$
as in \S\ref{cr22}. Let\/ $p>1$ and\/ $k\ge 2$, and for
${\bs\be}\in\R^n$ define $P_{\bs\be}:L^p_{k,{\bs\be}}(X')\ra
L^p_{k-2,{\bs\be}-2}(X')$ by $P_{\bs\be}(f)=\d^*(\psi^m\d f)$.
Then
\begin{itemize}
\item[{\rm(a)}] $P_{\bs\be}$ is Fredholm if and only if\/
${\bs\be}\in\bigl(\R\sm\D_{\sst\Si_1}\bigr)\t\cdots\t
\bigl(\R\sm\D_{\sst\Si_n}\bigr)$, and then
\e
\ind(P_{\bs\be})=-\sum_{i=1}^nN_\sSii(\be_i).
\label{cr5eq15}
\e
\item[{\rm(b)}] If\/ $\be_i>0$ for all\/ $i$ then $P_{\bs\be}$
is injective.
\end{itemize}
\label{cr5thm3}
\end{thm}

\begin{proof} Define $\check g=\psi^{2m/(m-2)}g$, a Riemannian
metric on $X'$ conformally equivalent to $g$. This is well-defined
as $m>2$. Since $\psi\vert_{X'}\in C_{\bs 0}^\iy(X')$ and $\psi(x)\ra
\psi(x_i)>0$ as $x\ra x_i$ for $i=1,\ldots,n$, one can show as in
Theorem \ref{cr5thm2} that $\check g$ induces a metric $\check d$
on $X$ and $(X,\check d\,)$ is a Riemannian manifold with conical
singularities at~$x_1,\ldots,x_n$.

Furthermore, $(X,\check d\,)$ has the {\it same} cones $C_i$ and
rates $\nu_i$ as does $(X,d)$ induced by $g$. The cones $C_i$
do not change because as Riemannian cones they are rescaled by
a homothety multiplying distances by $\psi(x_i)^{m/(m-2)}$,
but this gives the same Riemannian cone. As vector spaces of
functions $L^p_{k,{\bs\be}}(X')$ and $L^p_{k-2,{\bs\be}-2}(X')$
are the same for $g$ and $\check g$, with equivalent norms.

Write $\d_g^*,\d_{\check g}^*$ for $\d^*$ computed using
$g,\check g$ respectively. Let $\check\De=\d_{\check g}^*\d$
be the Laplacian of $\check g$ on functions. An elementary
calculation shows that
\begin{equation*}
\d_g^*(\psi^m\d f)=\psi^{m^2/(m-2)}\check\De f
\end{equation*}
for twice differentiable functions $f$ on $X'$.
Thus~$P_{\bs\be}=\psi^{\smash{m^2/(m-2)}}\check\De^p_{k,{\bs\be}}$.

Now multiplication by $\psi^{\smash{m^2/(m-2)}}$ gives an
automorphism of $L^p_{\smash{k-2,{\bs\be}-2}}(X')$. So $P_{\bs\be}$
is Fredholm, or injective, if and only if $\check\De^p_{k,{\bs\be}}$
is. Therefore (a) follows from Theorems \ref{cr2thm3} and
\ref{cr2thm6}, and (b) from part (a) of Lemma~\ref{cr2lem4}.
\end{proof}

By a similar proof we modify Theorem \ref{cr2thm8}, giving a
result needed in~\cite{Joyc3,Joyc5}.

\begin{thm} Let\/ $(M,J,\om,\Om)$ be an almost Calabi--Yau $m$-fold,
and define $\psi:M\ra(0,\iy)$ as in \eq{cr3eq3}. Suppose $X$ is a
compact SL\/ $m$-fold in $M$ with conical singularities at
$x_1,\ldots,x_n$, and let\/ $X',K,R',\Si_i,\Up_i,\phi_i,S_i$
and\/ $\mu_i$ be as in Definition \ref{cr2def5}, $\D_\sSii$
as in Definition \ref{cr2def3}, and\/ $\rho$ as in Definition
\ref{cr2def4}. Define
\e
\begin{split}
Y_\sXp=\bigl\{\al\in C^\iy(T^*X'):\,&\d\al=0,\quad
\d^*(\psi^m\al)=0,\\
&\text{$\md{\na^k\al}=O(\rho^{-1-k})$ for $k\ge 0$}\bigr\}.
\end{split}
\label{cr5eq16}
\e
Then $\pi:Y_\sXp\!\ra\!H^1(X',\R)$ given by $\pi:\al\!\mapsto\![\al]$
is an isomorphism. Furthermore:
\begin{itemize}
\item[{\rm(a)}] Fix $\al\in Y_\sXp$. By Hodge theory there exists
a unique $\ga_i\in C^\iy(T^*\Si_i)$ with\/ $\d\ga_i=\d^*\ga_i=0$
for $i=1,\ldots,n$, such that the image of\/ $\pi(\al)$ under the
map $H^1(X',\R)\ra\bigoplus_{i=1}^nH^1(\Si_i,\R)$ of\/ \eq{cr2eq16}
is $\bigl([\ga_1],\ldots,[\ga_n]\bigr)$. There exist unique $T_i\in
C^\iy\bigl(\Si_i\t(0,R')\bigr)$ for $i=1,\ldots,n$ such that
\end{itemize}
\ea
(\Up_i\circ\phi_i)^*(\al)&=\pi_i^*(\ga_i)+\d T_i
\quad \text{on $\Si_i\t(0,R')$ for $i=1,\ldots,n$, and}
\label{cr5eq17}\\
\na^kT_i(\si,r)&=O(r^{\nu_i-k})\qquad
\begin{aligned}
&\text{as\/ $r\ra 0$, for all\/ $k\ge 0$ and}\\
&\text{$\nu_i\in(0,\mu_i-2)$ with\/ $(0,\nu_i]\cap\D_\sSii=\emptyset$.}
\end{aligned}
\label{cr5eq18}
\ea
\begin{itemize}
\item[{\rm(b)}] Suppose $\ga_i\in C^\iy(T^*\Si_i)$ with\/
$\d\ga_i=\d^*\ga_i=0$ for $i=1,\ldots,n$, and the image
of\/ $\bigl([\ga_1],\ldots,[\ga_n]\bigr)$ under
$\bigoplus_{i=1}^nH^1(\Si_i,\R)\ra H^2_{\rm cs}(X',\R)$ in
\eq{cr2eq16} is $[\be]$ for some exact\/ $2$-form $\be$ on $X'$
supported on $K$. Then there exists $\al\in C^\iy(T^*X')$ with\/
$\d\al=\be$, $\d^*(\psi^m\al)=0$ and\/ $\md{\na^k\al}=O(\rho^{-1-k})$
for $k\ge 0$, such that\/ \eq{cr5eq17} and\/ \eq{cr5eq18} hold
for~$T_i\in C^\iy\bigl(\Si_i\t(0,R')\bigr)$.
\item[{\rm(c)}] Let\/ $f\in C^\iy(X')$ with\/ $\md{\na^kf}=
O(\rho^{{\bs\mu}-4-k})$ for $k\ge 0$ and\/ $\int_{X'}f\,\d V=0$.
Then there exists a unique exact\/ $1$-form $\al$ on $X'$ with\/
$\d^*(\psi^m\al)=f$ and\/ $\md{\na^k\al}=O(\rho^{-1-k})$ for $k\ge 0$,
such that\/ \eq{cr5eq17} and\/ \eq{cr5eq18} hold for $\ga_i=0$
and\/~$T_i\in C^\iy\bigl(\Si_i\t(0,R')\bigr)$.
\end{itemize}
\label{cr5thm4}
\end{thm}

\subsection{Improving the rates of convergence $\mu_i$}
\label{cr53}

We shall use the analysis results of \S\ref{cr2} to show that we
can improve the {\it rate} $\mu_i$ of the conical singularity $x_i$
in $X$ to all possibilities allowed by Definition~\ref{cr3def5}.

\begin{thm} In the situation of Theorem \ref{cr4thm4}
and Lemma \ref{cr4lem} suppose $\mu_i'\in(2,3)$ with\/
$(2,\mu_i']\cap\D_\sSii=\emptyset$ for
$i=1,\ldots,n$. Then
\e
\begin{gathered}
\bmd{\na^k(\phi_i-\iota_i)}=O(r^{\mu_i'-1-k}),\quad
\bmd{\na^k\eta_i}=O(r^{\mu_i'-1-k})\quad\text{and}\\
\bmd{\na^kA_i}=O(r^{\mu_i'-k})
\quad\text{as $r\ra 0$ for all\/ $k\ge 0$ and\/ $i=1,\ldots,n$.}
\end{gathered}
\label{cr5eq19}
\e

Hence $X$ has conical singularities at $x_i$ with cone $C_i$
and rate $\mu_i'$, for all possible rates $\mu_i'$ allowed by
Definition \ref{cr3def5}. Therefore, the definition of
conical singularities is essentially independent of the
choice of rate~$\mu_i$.
\label{cr5thm5}
\end{thm}

\begin{proof} Define a smooth function $A:X'\ra\R$ by
\e
A\bigl(\Up_i\circ\phi_i(\si,r)\bigr)=A_i(\si,r)
\quad\text{on $S_i$ for $i=1,\ldots,n$,}
\label{cr5eq20}
\e
and extend $A$ smoothly over $K=X'\sm(S_1\cup\cdots\cup S_n)$.
Then $A\in C^\iy_{\bs\mu}(X')$ by \eq{cr5eq1}. Let $Q_i'$ and
$g_i$ be the push-forwards of $Q_i$ in \eq{cr5eq3}--\eq{cr5eq4}
and $\iota_i^*(g')$ from $\Si_i\t(0,R')$ to $S_i$ under
$\Up_i\t\phi_i$. Then \eq{cr5eq9} implies that
\begin{equation*}
\d^*_{g_i}\d A[x]=Q_i'\bigl(x,\d A(x),\na^2A(x)\bigr)
\quad\text{for all $x\in S_i$.}
\end{equation*}

Here $\d^*_{g_i}$ is computed using the exactly conical metric
$g_i$ on $S_i$, rather than the asymptotically conical metric
$g$. Rearranging yields
\e
\De A[x]=Q_i'\bigl(x,\d A(x),\na^2A(x)\bigr)+(\d^*_g-\d^*_{g_i})\d A[x]
\label{cr5eq21}
\e
for all $x\in S_i$, where $\De=\d^*_g\d$ is the Laplacian of $g$.

We shall prove the theorem by using an inductive argument
to improve the decay rate of $A$ and its derivatives step by
step until we show that $A\in C^\iy_{{\bs\mu}'}(X')$ for all
${\bs\mu}'$ satisfying the conditions of the theorem. The
next two lemmas will be needed for the `inductive step'.

\begin{lem} Let\/ $\la_i\in(2,3)$ and define $\hat\la_i=\min(3,2\la_i-2)$
for $i=1,\ldots,n$. Then if\/ $A\in C^\iy_{\bs\la}(X')$,
then~$\De A\in \smash{C^\iy_{\hat{\bs\la}-2}}(X')$.
\label{cr5lem1}
\end{lem}

\begin{proof} Suppose $A\in C^\iy_{\bs\la}(X')$. Then from
\eq{cr5eq8} we find that
\begin{equation*}
Q_i'\bigl(x,\d A(x),\na^2A(x)\bigr)=O(\rho)+O(\rho^{2\la_i-4})
+O(\rho^{2\la_i-4})\quad\text{for $x\in S_i$.}
\end{equation*}
As the asymptotic behaviour of $g$ on $S_i$ as $\rho\ra 0$ depends
on $\eta_i=\d A_i=\d A$ we have $g-g_i=O(\rho^{\la_i-2})$ on $S_i$.
Thus $(\d^*_g-\d^*_{g_i})\d A=O(\rho^{2\la_i-4})$. Combining these
with \eq{cr5eq21} gives $\De A=O(\rho)+O(\rho^{2{\bs\la}-4})=
O(\rho^{\smash{\hat{\bs\la}-2}})$, by definition of $\hat{\bs\la}$.
The argument easily extends to derivatives of $\De A$, and
so~$\De A\in C^\iy_{\hat{\bs\la}-2}(X')$.
\end{proof}

\begin{lem} Suppose $\la_i,\hat\la_i\in(2,3)$ with\/ $(2,\la_i]\cap
\D_\sSii=(2,\hat\la_i]\cap\D_\sSii=\emptyset$
for $i=1,\ldots,n$. Let\/ $p>1$ and\/ $k\ge 2$. Then if\/ $A\in
L^p_{k,{\bs\la}}(X')$ and\/ $\De A\in L^p_{k-2,\smash{\hat{\bs\la}}-2}(X')$,
then~$A\in L^p_{k,\smash{\hat{\bs\la}}}(X')$.
\label{cr5lem2}
\end{lem}

\begin{proof} Define $q>1$ by $\frac{1}{p}+\frac{1}{q}=1$.
By Theorem \ref{cr2thm5}, $\De A$ is orthogonal to
$\Ker(\De^q_{2,-{\bs\la}+2-m})$ as $A\in L^p_{k,{\bs\la}}(X')$.
But the conditions on ${\bs\la},\hat{\bs\la}$ imply that
${\bs\la},\hat{\bs\la}$ lie in the same connected component of
\eq{cr2eq9}, and therefore $-{\bs\la}+2-m,-\hat{\bs\la}+2-m$ also
lie in the same connected component of \eq{cr2eq9}. Hence
$\Ker(\De^q_{2,-{\bs\la}+2-m})=\Ker(\De^q_{2,-\smash{\hat{\bs\la}}+2-m})$
by Theorem~\ref{cr2thm4}.

Therefore $\De A$ lies in $L^p_{k-2,\smash{\hat{\bs\la}}-2}(X')$ by
assumption and is orthogonal to $\Ker(\De^q_{2,-\smash{\hat{\bs\la}}+2-m})$.
So by Theorem \ref{cr2thm5}, $\De A$ lies in the image of
$\De^p_{k,\smash{\hat{\bs\la}}}$, and $\De A=\De A'$ for some
$A'\in L^p_{k,\smash{\hat{\bs\la}}}(X')$. Thus $\De(A-A')=0$ and
$A-A'$ is harmonic. But $A-A'=O(\rho^2)$ as $\la_i,\hat\la_i>2$,
so $(A-A')(x)\ra 0$ as $x\ra x_i$ in $X'$. Hence using the
maximum principle \cite[\S 3]{GiTr} we see that $A-A'\equiv 0$,
giving $A=A'$ and~$A\in L^p_{k,\smash{\hat{\bs\la}}}(X')$.
\end{proof}

Now we can prove the theorem. As $A\in C^\iy_{\bs\mu}(X')$ from
above, Lemma \ref{cr5lem1} shows that $\De A\in C^\iy_{\hat{\bs\mu}
-2}(X')$, where $\hat\mu_i=\min(3,2\mu_i-2)$ for $i=1,\ldots,n$.
Note that $\hat\mu_i>\mu_i$ as $\mu_i\in(2,3)$.
Therefore if $p>1$, $k\ge 2$ and $2<\la_i<\mu_i$, $2<\hat\la_i<\hat\mu_i$
for $i=1,\ldots,n$ we see that $A\in L^p_{k,{\bs\la}}(X')$ and
$\De A\in L^p_{k-2,\smash{\hat{\bs\la}}-2}(X')$, since $C^\iy_{\bs\mu}(X')
\subset L^p_{k,{\bs\la}}(X')$ and~$C^\iy_{\hat{\bs\mu}-2}(X')\subset
L^p_{k-2,\smash{\hat{\bs\la}}-2}(X')$.

Hence Lemma \ref{cr5lem2} shows that for all $\hat\la_i\in(2,\hat\mu_i)$
with $(2,\hat\la_i]\cap\D_\sSii=\emptyset$ for $i=1,\ldots,n$
we have $A\in L^p_{k,\smash{\hat{\bs\la}}}(X')$. As this holds for all $k
\ge 2$, Theorem \ref{cr2thm1} then proves that $A\in C^\iy_{\hat{\bs\la}}
(X')$. Thus starting with $A\in C^\iy_{\bs\mu}(X')$ we have shown that
$A\in C^\iy_{\hat{\bs\la}}(X')$ for all $\hat{\bs\la}=(\hat\la_1,\ldots,
\hat\la_n)$ with $2<\hat\la_i<\hat\mu_i$ and~$(2,\hat\la_i]\cap
\D_\sSii=\emptyset$.

Since $\hat\mu_i>\mu_i$, this is an improvement in the rate
of convergence of $A$. If $\hat\mu_i=3$ or $(2,\hat\mu_i]\cap
\D_\sSii\ne\emptyset$ then we have proved what we want
for convergence of $A$ on $S_i$. Otherwise $\hat\mu_i=2\mu_i-2$,
so that $\hat\mu_i-2=2(\mu_i-2)$. Applying the same argument $j$
times, we find that either we prove what we want for the
convergence of $A$ on $S_i$, or else $A\in C^\iy_{\hat{\bs\la}}(X')$
for all $\hat{\bs\la}$ with $2<\hat\la_i<\ti\mu_i<3$ with $\ti\mu_i-2=
2^j(\mu_i-2)$.

If $2^j(\mu_i-2)\ge 1$ this gives $\ti\mu_i\ge 3$, a contradiction, so
the process must terminate, and therefore for all ${\bs\mu}'$ satisfying
the conditions of the theorem we have $A\in C^\iy_{{\bs\mu}'}(X')$.
Equation \eq{cr5eq20} then gives $\md{\na^kA_i}=O(r^{\mu_i'-k})$
as $r\ra 0$ for all $k\ge 0$ and $i=1,\ldots,n$, the final equation
of \eq{cr5eq19}. The first two equations of \eq{cr5eq19} then
follow as for \eq{cr5eq1}. This completes the proof of
Theorem~\ref{cr5thm5}.
\end{proof}

\section{Geometric Measure Theory and tangent cones}
\label{cr6}

We now review some {\it Geometric Measure Theory}, and apply
it to special Lagrangian geometry. An introduction to the
subject is provided by Morgan \cite{Morg} and an in-depth
(but dated) treatment by Federer \cite{Fede}, and Harvey and
Lawson \cite[\S II]{HaLa} relate Geometric Measure Theory to
calibrated geometry.

Geometric Measure Theory studies measure-theoretic
generalizations of submanifolds called {\it integral currents},
which may be very singular, and is particularly powerful for
{\it minimal\/} submanifolds. We shall distinguish between
submanifolds or currents which are {\it volume-minimizing}
(local minima of the volume functional), and those which are
{\it minimal\/} (stationary points of the volume functional).
Stronger results are available for the volume-minimizing case.

We can consider {\it special Lagrangian integral currents},
a natural class of singular SL $m$-folds with strong compactness
properties, which are automatically volume-minimizing. Our main
result, Theorem \ref{cr6thm5}, says that if the {\it tangent
cones} of an SL integral current $T$ satisfy a certain condition
then $T$ is actually an SL $m$-fold with conical singularities,
in the sense of \S\ref{cr33}. Throughout we suppose~$m>2$.

\subsection{Introduction to Geometric Measure Theory}
\label{cr61}

Let $(M,g)$ be a complete Riemannian manifold. One defines a
class of $m$-dimensional {\it rectifiable currents} in $M$,
which are measure-theoretic generalizations of compact,
oriented $m$-submanifolds $N$ with boundary $\pd N$ in $M$,
with integer {\it multiplicities}. Here $N$ with multiplicity
$k$ is like $k$ copies of $N$ superimposed, and changing the
orientation of $N$ changes the sign of the multiplicity. This
enables us to add and subtract submanifolds.

If $T$ is an $m$-dimensional rectifiable current, one can define
the {\it volume} $\vol(T)$ of $T$, by Hausdorff $m$-measure. If
$\vp$ is a compactly-supported $m$-form on $M$ then one can
define $\int_T\vp$. Thus we can regard $T$ as a {\it current},
that is, an element $\vp\mapsto\int_T\vp$ of the dual space
$(\D^m)^*$ of the vector space $\D^m$ of smooth
compactly-supported $m$-forms on $M$. This induces a topology
on the space of rectifiable currents in~$M$.

Let $T$ be a $m$-dimensional rectifiable current, and define
an $(m\!-\!1)$-current $\pd T$ by ${\pd T}\cdot\al=\int_T\d\al$
for $\al\in\D^{m-1}$. We call $T$ an {\it integral
current} if $\pd T$ is a rectifiable current. By 
\cite[5.5]{Morg}, \cite[4.2.17]{Fede}, integral currents
have strong {\it compactness properties}.

Harvey and Lawson \cite[\S II]{HaLa} discuss calibrated geometry
and Geometric Measure Theory. They show that on a Riemannian
manifold $(M,g)$ with calibration $\vp$ one can define
{\it integral $\vp$-currents}, that is, integral currents
which are calibrated w.r.t.\ $\vp$, and that they are
{\it volume-minimizing} in their homology class.

In particular, as in \S\ref{cr3} SL $m$-folds in $\C^m$ and in
an almost Calabi--Yau manifold $M$ may be defined as calibrated
submanifolds, using the conformally rescaled metric $\ti g$ on $M$.
Therefore we can define {\it special Lagrangian integral currents}
in $\C^m$ and in almost Calabi--Yau manifolds $M$, and they are
{\it volume-minimizing currents} w.r.t.\ an appropriate metric.

\subsection{Tangent cones}
\label{cr62}

Next we discuss {\it tangent cones} of volume-minimizing integral
currents, a generalization of tangent spaces of submanifolds,
as in \cite[9.7]{Morg}. Define the {\it interior} $T^\circ$ of
$T$ to be $T\sm\pd T$ (that is,~$\supp T\sm\supp\pd T$).

\begin{dfn} An integral current $C$ in $\R^n$ is called a
{\it cone} if $C=tC$ for all $t>0$, where $t:\R^n\ra\R^n$
acts by dilations in the obvious way. Let $T$ be an integral
current in $\R^n$, and let ${\bf x}\in T^\circ$. We say
that $C$ is a {\it tangent cone} to $T$ at $\bf x$ if there
exists a decreasing sequence $r_1>r_2>\cdots$ tending to
zero such that $r_j^{-1}(T-{\bf x})$ converges to $C$ as
an integral current as~$j\ra\iy$.

More generally, if $(M,g)$ is a complete Riemannian
$n$-manifold, $T$ is an integral current in $M$, and
$x\in T^\circ$, then one can define a {\it tangent cone}
$C$ to $T$ at $x$, which is an integral current cone
in the Euclidean vector space $T_xM$. Identifying $M$
with $\R^n$ near $x$ using a coordinate system, the
two notions of tangent cone coincide.
\label{cr6def1}
\end{dfn}

The next result follows from Morgan \cite[p.~94-95]{Morg}, Federer
\cite[5.4.3]{Fede} and Harvey and Lawson~\cite[Th.~II.5.15]{HaLa}.

\begin{thm} Let\/ $(M,g)$ be a complete Riemannian manifold,
and\/ $T$ a volume-minimizing integral current in $M$. Then for all\/
$x\in T^\circ$, there exists a tangent cone $C$ to $T$ at\/
$x$. Moreover $C$ is itself a volume-minimizing integral current in
$T_xM$ with\/ $\pd C=\emptyset$, and if\/ $T$ is calibrated
with respect to a calibration $\vp$ on $(M,g)$, then $C$ is
calibrated with respect to the constant calibration
$\vp\vert_x$ on~$T_xM$.
\label{cr6thm1}
\end{thm}

Note that the theorem does {\it not\/} claim that the tangent
cone $C$ is unique, and in fact it is an important open question
whether a volume-minimizing integral current has a unique tangent cone at
each point of $T^\circ$. However, Leon Simon \cite{Simo1,Simo2},
improving an earlier result of Allard and Almgren \cite{AlAl},
shows that if some tangent cone $C$ is nonsingular and multiplicity
1 away from 0, then $C$ is the unique tangent cone, and $T$
converges to $C$ in a $C^1$ sense. For later use we model
the result on the notation of Definition~\ref{cr3def5}.

\begin{thm} Let\/ $C$ be an $m$-dimensional oriented minimal
cone in $\R^n$ with\/ $C'=C\sm\{0\}$ nonsingular, so that\/
$\Si=C\cap{\cal S}^{n-1}$ is a compact, oriented, nonsingular,
embedded, minimal $(m\!-\!1)$-submanifold of\/ ${\cal S}^{n-1}$.
Define $\iota:\Si\t(0,\iy)\ra C'\subset\R^n$ by $\iota(\si,r)=r\si$.
Let\/ $(M,g)$ be a complete Riemannian $n$-manifold and\/ $x\in M$.
Fix an isometry $\up:\R^n\ra T_xM$, and choose an embedding
$\Up:B_R\ra M$ with\/ $\Up(0)=x$ and\/ $\d\Up\vert_0=\up$,
where $B_R$ is the ball of radius $R>0$ about\/~$0\in\R^n$.

Suppose that\/ $T$ is a minimal integral current in $M$ with\/
$x\in T^\circ$, and that\/ $\up_*(C)$ is a tangent cone to $T$
at $x$ with multiplicity $1$. Then $\up_*(C)$ is the unique
tangent cone to $T$ at\/ $x$. Furthermore there exists
$R'\in(0,R]$ and an embedding $\phi:\Si\t(0,R')\ra B_{R'}
\subseteq B_R$ with
\e
\bmd{\phi(\si,r)}\equiv r,\quad
\bmd{\phi-\iota}=o(r)\quad\text{and}\quad
\bmd{\na(\phi-\iota)}=o(1)
\quad\text{as $r\ra 0$,}
\label{cr6eq1}
\e
such that\/ $T\cap\bigl(\Up(B_{R'})\sm\{x\}\bigr)$ is the
embedded submanifold\/ $\Up\circ\phi\bigl(\Si\t(0,R')\bigr)$,
with multiplicity~$1$.
\label{cr6thm2}
\end{thm}

\begin{proof} This follows from \cite[Cor., p.~564]{Simo1} and
\cite[Th.~5.7]{Simo2}, which are equivalent results, the latter
more explicit. Simon claims only that $\phi$ is $C^2$ rather
than smooth, but smoothness follows from standard regularity
results for minimal submanifolds.
\end{proof}

We define {\it Jacobi fields} on $\Si$, following
Lawson~\cite[p.~46-52]{Laws}.

\begin{dfn} Let $\Si$ be a compact, minimal submanifold in
the unit sphere ${\cal S}^{n-1}$ in $\R^n$. Let $\nu$ be the
normal bundle of $\Si$ in ${\cal S}^{n-1}$, so that
$T{\cal S}^{n-1}\vert_\Si=\nu\op T\Si$ is an orthogonal
splitting. Let $g$ and $g_\sSi$ be the Riemannian
metrics on ${\cal S}^{n-1}$ and $\Si$ induced by the
Euclidean metric on~$\R^n$.

Let $\na^\nu$ be the connection on $\nu$ defined by projecting
the Levi-Civita connection of $g$ on $T{\cal S}^{n-1}\vert_\Si$
to $\nu$. Let $\De^\nu:C^\iy(\nu)\ra C^\iy(\nu)$ be the Laplacian
$(\na^\nu)^*\na^\nu$ defined using $\na^\nu,g$ and $g_\sSi$.
Define maps ${\cal R},{\cal B}:C^\iy(\nu)\ra C^\iy(\nu)$ by
\begin{equation*}
{\cal R}(w)^i=\pi^\nu\bigl(R^i_{\phantom{i}jkl}g_\sSi^{jk}w^l\bigr)
\quad\text{and}\quad
{\cal B}(w)^a=B^a_{bc}B^i_{jk}g_\sSi^{bj}g_\sSi^{ck}g_{il}w^l,
\end{equation*}
using the index notation for tensors, where $R^i_{\phantom{i}jkl}$ is
the Riemann curvature of $g$, $B^i_{jk}\in C^\iy(\nu\ot S^2T^*\Si)$ is
the second fundamental form of $\Si$ in ${\cal S}^{n-1}$, and $\pi^\nu$
is the orthogonal projection from $T{\cal S}^{n-1}$ to~$\nu$.

We call a normal vector field $w\in C^\iy(\nu)$ to $\Si$ in
${\cal S}^{n-1}$ a {\it Jacobi field\/} if
\e
\De^\nu w-{\cal R}(w)+{\cal B}(w)=0.
\label{cr6eq2}
\e
Jacobi fields are zeroes of the linearization at $\Si$ of the
Euler--Lagrange equation for the volume of submanifolds $\Si'$
in ${\cal S}^{n-1}$. Therefore a Jacobi field is an {\it
infinitesimal deformation of\/ $\Si$ as a minimal submanifold},
a null direction of the second variation of volume for submanifolds.

In particular, the Lie algebra $\so(n)$ of isometries of
${\cal S}^n$ clearly induce infinitesimal deformations of $\Si$
as a minimal submanifold, and so Jacobi fields. Regarding
$v\in\so(n)$ as a vector field on ${\cal S}^{n-1}$, the
corresponding Jacobi field on $\Si$ is $w=\pi^\nu(v\vert_\Si)$.
However, for some $\Si$ not all Jacobi fields come from $\so(n)$
in this way. Note that as $\Si$ is compact and \eq{cr6eq2} is
an elliptic equation, the Jacobi fields form a finite-dimensional
vector space.
\label{cr6def2}
\end{dfn}

Now by Allard and Almgren \cite[p.~215]{AlAl}, or equivalently
by Adams and Simon \cite[Th.~1]{AdSi}, if the Jacobi fields
on $\Si$ satisfy a condition then we can strengthen the rate
of convergence in~\eq{cr6eq1}.

\begin{thm} Let\/ $C$ be an $m$-dimensional oriented minimal
cone in $\R^n$ with\/ $C'=C\sm\{0\}$ nonsingular, and set\/
$\Si=C\cap{\cal S}^{n-1}$. Suppose that\/ $\Si$ satisfies
\begin{itemize}
\item[$(*)$] Each Jacobi field\/ $w$ of\/ $\Si$ in ${\cal S}^{n-1}$
exponentiates to a smooth\/ $1$-parameter family $\bigl\{\Si_t:t\in
(-\ep,\ep)\bigr\}$ of minimal submanifolds in ${\cal S}^{n-1}$ for
$\ep>0$, with\/ $\Si_0=\Si$ and velocity $w$ at\/~$t=0$.
\end{itemize}
Then for some $\mu>2$, the map $\phi$ of Theorem \ref{cr6thm2}
satisfies
\e
\bmd{\phi-\iota}=O(r^{\mu-1})\quad\text{and}\quad
\bmd{\na(\phi-\iota)}=O(r^{\mu-2})
\quad\text{as $r\ra 0$.}
\label{cr6eq3}
\e
\label{cr6thm3}
\end{thm}

Adams and Simon \cite[Th.~1(ii)]{AdSi} also study the case
when condition $(*)$ does not hold, and prove:

\begin{thm} Let\/ $C$ be an $m$-dimensional oriented minimal
cone in $\R^n$ with\/ $C'=C\sm\{0\}$ nonsingular, and set\/
$\Si=C\cap{\cal S}^{n-1}$. Suppose that condition $(*)$ of
Theorem \ref{cr6thm3} does not hold, and also that a certain sign
condition {\rm\cite[p.~232]{AdSi}} holds for some Jacobi field.

Then there exist large families of minimal integral currents
$T$ in $\R^n$ with\/ $0\in T^\circ$ such that\/ $C$ is a
tangent cone to $T$ at\/ $0$ with multiplicity $1$, and for
some $\al\in(0,1]$ the map $\phi$ in Theorem \ref{cr6thm2}
with\/ $\Up=\id:B_R\ra B_R\subset\R^n$ decays exactly at rate
\e
\bmd{\phi-\iota}=O\bigl(r\md{\log r}^{-\al}\bigr)
\quad\text{and}\quad
\bmd{\na(\phi-\iota)}=O\bigl(\md{\log r}^{-\al}\bigr)
\quad\text{as $r\ra 0$.}
\label{cr6eq4}
\e
\label{cr6thm4}
\end{thm}

Here is what we mean by the `sign condition' above. If
condition $(*)$ fails then there exists an integer $p>2$
and a nonzero homogeneous degree $p$ real polynomial $P$
on the Jacobi fields. If $P(w)>0$ then we can construct
minimal integral currents $T$ near 0 in $\R^n$ for which
\e
\phi(\si,r)=\iota(\si,r)+r\md{\log r}^{-1/(p-2)}w(\si)+
\text{lower order terms}
\label{cr6eq5}
\e
as $r\ra 0$. Thus \eq{cr6eq4} holds exactly for $\al=1/(p-2)$.
If $P(w)<0$ then we can instead construct {\it Asymptotically
Conical\/} minimal integral currents $T$ near $\iy$ in $\R^n$
for which \eq{cr6eq5} holds as $r\ra\iy$. We need $P(w)>0$
for some Jacobi field $w$, which is automatic when $p$ is
odd, and hence in the `most generic' case~$p=3$.

\subsection{Tangent cones of special Lagrangian $m$-folds}
\label{cr63}

We shall now specialize the results of \S\ref{cr62} to
the case when $T$ is a {\it special Lagrangian integral
current} in an almost Calabi--Yau $m$-fold $(M,J,\om,\Om)$.
Our aim is to prove that if the tangent cones of $T$ satisfy
certain conditions then $T$ satisfies Definition \ref{cr3def5},
and so is an {\it SL\/ $m$-fold with conical singularities}.

By restricting to special Lagrangian currents we can strengthen
Theorem \ref{cr6thm3}, as condition $(*)$ need not hold for all
Jacobi fields $w$, but only for those which represent infinitesimal
deformations of $C$ {\it as a special Lagrangian cone}, rather than
as a minimal cone.

\begin{dfn} Let $C$ be an SL cone in $\C^m$ with $C'=C\sm\{0\}$
nonsingular, and set $\Si=C\cap{\cal S}^{2m-1}$. Then $\Si$ is a
compact, nonsingular, minimal Legendrian submanifold of
${\cal S}^{2m-1}$. Define $\iota:\Si\t(0,\iy)\ra\C^m$ by
$\iota(\si,r)=r\si$, with image $C'$. Let $g_\sSi=g'\vert_\Si$
be the metric on $\Si$ and $\De_\sSi$ the Laplacian on~$\Si$.

Suppose $v\in C^\iy(\Si)$ is an eigenfunction of $\De_\sSi$
with eigenvalue $2m$. Then $u:r\si\mapsto r^2v(\si)$ is a
homogeneous harmonic function on $C'$ of order 2, by Lemma
\ref{cr2lem1}. Thus $\d u$ is a homogeneous closed and
coclosed 1-form on $C'$ of order 1. Let $\nu\ra C'$ be the
normal bundle of $C'$ in $C^m$. Then $\nu\cong T^*C'$ by
the usual isomorphism. So $\d u$ corresponds to a homogeneous
section of $\nu$ of order 1, which is an {\it infinitesimal
deformation of\/ $C$ as an SL cone}.

Define $w_v$ to be the restriction of this section to
$\Si\subset C'$. Then $w_v$ is a smooth section of the
normal bundle of $\Si$ in ${\cal S}^{2m-1}$, and is a
{\it Jacobi field} on $\Si$ in the sense of Defintion
\ref{cr6def2}. Define a {\it special Lagrangian Jacobi
field} to be a Jacobi field $w_v$ on $\Si$ constructed
from a $\De_\sSi$ $2m$-eigenfunction $v\in C^\iy(\Si)$
in this way.

We call $C$ {\it Jacobi integrable} if it satisfies the condition
\begin{itemize}
\item[$(**)$] Each special Lagrangian Jacobi field $w_v$ of
$\Si$ in ${\cal S}^{2m-1}$ exponentiates to a smooth 1-parameter
family $\bigl\{\Si_t:t\in(-\ep,\ep)\bigr\}$ for $\ep>0$ with
$\Si_0=\Si$ and velocity $w_v$ at $t=0$, where $\Si_t=C_t\cap
{\cal S}^{2m-1}$ for $C_t$ a special Lagrangian cone in~$\C^m$.
\end{itemize}
That is, each special Lagrangian Jacobi field should be integrable.

Each element $x$ of the Lie algebra $\su(m)$, regarded as vector
field on ${\cal S}^{2m-1}$, induces an infinitesimal deformation
of $C$ as a special Lagrangian cone, so that $\pi^\nu(x\vert_\Si)$
is a special Lagrangian Jacobi field $w_v$ on $\Si$. The
corresponding eigenfunction $v\in C^\iy(\Si)$ is the restriction
to $\Si$ of the unique moment map $\mu:\C^m\ra\R$ of $x$ with
$\mu(0)=0$. Now Jacobi fields $w_v$ constructed from $x\in\su(m)$
in this way automatically satisfy $(**)$, as $\exp(tx)\in\SU(m)$
for $t\in\R$, so $C_t=\exp(tx)C$ is a special Lagrangian cone,
and $\Si_t=C_t\cap{\cal S}^{2m-1}$ satisfies the conditions.

Define $C$ to be {\it rigid\/} if all special Lagrangian Jacobi
fields $w_v$ on $\Si$ come from $\su(m)$ as above. Then $C$
{\it rigid implies $C$ Jacobi integrable}, from above. There
is a simple test for rigidity: let $G$ be the Lie subgroup of
$\SU(m)$ preserving $C$, and $\g$ the Lie algebra of $G$.
Then the special Lagrangian Jacobi fields on $\Si$ from
$\su(m)$ are a vector space isomorphic to $\su(m)/\g$,
with dimension $m^2-1-\dim G$. Therefore $C$ {\it is rigid
if and only if the multiplicity of the eigenvalue $2m$ of\/
$\De_\sSi$ is} $m^2-1-\dim G$, that is, if $m_\sSi(2)
=m^2-1-\dim G$ in the notation of Definition \ref{cr2def3}.
This may be taken as an alternative definition of rigidity.
\label{cr6def3}
\end{dfn}

Now we can prove the main result of this section.

\begin{thm} Let\/ $(M,J,\om,\Om)$ be an almost Calabi--Yau
$m$-fold and define $\psi:M\ra(0,\iy)$ as in \eq{cr3eq3}.
Let\/ $x\in M$ and fix an isomorphism $\up:\C^m\ra T_xM$
with\/ $\up^*(\om)=\om'$ and\/ $\up^*(\Om)=\psi(x)^m\Om'$,
where $\om',\Om'$ are as in~\eq{cr3eq1}.

Suppose that\/ $T$ is a special Lagrangian integral current in
$M$ with\/ $x\in T^\circ$, and that\/ $\up_*(C)$ is a multiplicity
$1$ tangent cone to $T$ at\/ $x$, where $C$ is a Jacobi integrable
special Lagrangian cone in $\C^m$, in the sense of Definition
\ref{cr6def3}. Then $T$ has a conical singularity at\/ $x$, in
the sense of Definition~\ref{cr3def5}.

Suppose that\/ $T$ is a special Lagrangian integral current in
$M$ with\/ $\pd T=\emptyset$, and that every singular point of\/
$T$ has a Jacobi integrable multiplicity $1$ special Lagrangian
tangent cone. Then $T$ is a compact SL\/ $m$-fold in $M$ with
conical singularities, in the sense of Definition~\ref{cr3def5}.
\label{cr6thm5}
\end{thm}

\begin{proof} Let $(M,J,\om,\Om),x,\up$  and $T$ be as in the
first part of the theorem, and choose an embedding $\Up:B_R\ra M$
with $\Up(0)=x$, $\d\Up\vert_0=\up$ and $\Up^*(\om)=\om'$, as in
Definition \ref{cr3def5}. Then Theorem \ref{cr6thm2} applies, and
gives $R'\in(0,R]$ and an embedding $\phi:\Si\t(0,R')\ra B_R$ 
satisfying \eq{cr6eq1} such that $\Up\circ\phi$ parametrizes
$T$ near~$x$.

We would like to apply Theorem \ref{cr6thm3} to deduce that
$\phi$ satisfies \eq{cr6eq3}. Now following the proof of Theorem
\ref{cr6thm3} in \cite{AdSi}, we find that {\it either} \eq{cr6eq3}
holds, {\it or} we can construct a Jacobi field $w$ from $T$ by a
limiting process, which does not satisy $(*)$. Since $T$ is special
Lagrangian it turns out that $w$ must be a {\it special Lagrangian
Jacobi field}, and so does not satisfy~$(**)$.

But as $C$ is Jacobi integrable, condition $(**)$ holds for all
such $w$. Therefore $\phi$ satisfies \eq{cr6eq3} for some
$\mu>2$. Making $\mu$ smaller if necessary we can suppose
$\mu\in(2,3)$ and $\mu$ satisfies \eq{cr3eq4}. Then
\eq{cr6eq3} is equivalent to \eq{cr3eq5}, so $T$ satisfies
Definition \ref{cr3def5} near $x$, and has a conical
singularity at $x$ with identification $\up$, cone $C$ and
rate $\mu$. This completes the first part of the theorem.

For the second part, note that by the first part every
singular point of $T$ is a conical singularity, and so
is isolated. Thus by compactness of $M$ there are only
finitely many singular points $x_1,\ldots,x_n$ of $T$,
and it quickly follows that $T$ is a compact SL $m$-fold
with conical singularities.
\end{proof}

This is a {\it weakening} of Definition \ref{cr3def5}, in
that if $T$ satisfies the apparently much weaker condition
of having a certain kind of tangent cone at $x$, then $T$
actually has a conical singularity at~$x$.

Finally we discuss singularities $x$ of SL $m$-folds $X$
modelled on multiplicity one SL cones $C$ with $C\sm\{0\}$
nonsingular, but where $C$ is not Jacobi integrable. Then Theorem
\ref{cr6thm2} shows that $X$ can be parametrized near $x$
using a map $\phi:\Si\t(0,R')\ra B_R$ satisfying~\eq{cr6eq1}.

However, Theorem \ref{cr6thm3} suggests that the asymptotic
behaviour we should expect of $\phi$, at least for $X$ suitably
generic, is exactly that of \eq{cr6eq4} for some $\al\in(0,1]$.
This does not satisfy \eq{cr3eq5}, and so such singular points
will {\it not\/} be conical singularities in our sense.

This indicates that for SL cones $C$ which are not Jacobi
integrable, Definition \ref{cr3def5} is actually {\it too strong},
in that there should exist examples of singular SL $m$-folds
with tangent cone $C$ which are not covered by Definition
\ref{cr3def5}, since the decay conditions in \eq{cr3eq5} are
too strict. Nevertheless, we will continue to use Definition
\ref{cr3def5} in the sequels \cite{Joyc3,Joyc4,Joyc5,Joyc6},
because without it we will be unable to use the powerful
analytic framework of~\S\ref{cr2}.

\section{Asymptotically Conical SL $m$-folds}
\label{cr7}

Let $C$ be an SL cone in $\C^m$ with an isolated singularity
at 0. Sections \ref{cr3}--\ref{cr6} considered SL $m$-folds
with conical singularities, which are asymptotic to $C$ at 0.
We now discuss {\it Asymptotically Conical\/} SL $m$-folds
$L$ in $\C^m$, which are asymptotic to $C$ at infinity. Here
is the definition.

\begin{dfn} Let $C$ be a closed SL cone in $\C^m$ with isolated
singularity at 0 for $m>2$, and let $\Si=C\cap{\cal S}^{2m-1}$,
so that $\Si$ is a compact, nonsingular $(m-1)$-manifold, not
necessarily connected. Let $g_\sSi$ be the metric on $\Si$
induced by the metric $g'$ on $\C^m$ in \eq{cr3eq1}, and $r$ the
radius function on $\C^m$. Define $\iota:\Si\t(0,\iy)\ra\C^m$ by
$\iota(\si,r)=r\si$. Then the image of $\iota$ is $C\sm\{0\}$,
and $\iota^*(g')=r^2g_\sSi+\d r^2$ is the cone metric
on~$C\sm\{0\}$.

Let $L$ be a closed, nonsingular SL $m$-fold in $\C^m$. We
call $L$ {\it Asymptotically Conical (AC)} with {\it rate}
$\la<2$ and {\it cone} $C$ if there exists a compact subset
$K\subset L$ and a diffeomorphism $\vp:\Si\t(T,\iy)\ra L\sm K$
for some $T>0$, such that
\e
\bmd{\na^k(\vp-\iota)}=O(r^{\la-1-k})
\quad\text{as $r\ra\iy$ for $k=0,1$.}
\label{cr7eq1}
\e
Here $\na,\md{\,.\,}$ are computed using the cone metric~$\iota^*(g')$.
\label{cr7def1}
\end{dfn}

This is very similar to Definition \ref{cr3def5}, and in
fact there are strong similarities between the theories of SL
$m$-folds with conical singularities and of Asymptotically
Conical SL $m$-folds. Note that we do {\it not\/} impose any
condition on $\la$ analogous to \eq{cr3eq4}, although we
could. We continue to assume~$m>2$.

We begin in \S\ref{cr71} by defining {\it cohomological
invariants} $Y(L),Z(L)$ of $L$ in $H^*(\Si,\R)$, which
have no parallel in the conical singularities case. Then
\S\ref{cr72} and \S\ref{cr73} develop the analogues of
parts of \S\ref{cr4} and \S\ref{cr5} for AC SL $m$-folds.
The {\it deformation theory} of AC SL $m$-folds is studied
by Marshall \cite{Mars}, and {\it examples} of AC SL
$m$-folds will be discussed in~\cite[\S 6.4]{Joyc6}.

\subsection{Cohomological invariants of AC SL $m$-folds}
\label{cr71}

Let $L$ be an AC SL $m$-fold in $\C^m$ with cone $C$, and set
$\Si=C\cap{\cal S}^{2m-1}$. Using the notation of \S\ref{cr24},
as in \eq{cr2eq15} there is a long exact sequence
\e
\cdots\ra
H^k_{\rm cs}(L,\R)\ra H^k(L,\R)\ra H^k(\Si,\R)\ra
H^{k+1}_{\rm cs}(L,\R)\ra\cdots.
\label{cr7eq2}
\e
We shall define {\it cohomological invariants} $Y(L),Z(L)$ of~$L$.

\begin{dfn} Let $L$ be an AC SL $m$-fold in $\C^m$ with cone $C$,
and let $\Si=C\cap{\cal S}^{2m-1}$. As $\om',\Im\Om'$ in \eq{cr3eq1}
are closed forms with $\om'\vert_L\equiv\Im\Om'\vert_L\equiv 0$, they
define classes in the relative de Rham cohomology groups $H^k(\C^m;
L,\R)$ for $k=2,m$. But for $k>1$ we have the exact sequence
\begin{equation*}
0=H^{k-1}(\C^m,\R)\ra H^{k-1}(L,\R){\buildrel\cong\over\longra}
H^k(\C^m;L,\R)\ra H^k(\C^m,\R)=0.
\end{equation*}
Define $Y(L)\in H^1(\Si,\R)$ to be the image of $[\om']$
in $H^2(\C^m;L,\R)\cong H^1(L,\R)$ under the map
$H^1(L,\R)\ra H^1(\Si,R)$ of \eq{cr7eq2}, and $Z(L)\in
H^{m-1}(\Si,\R)$ to be the image of $[\Im\Om']$ in
$H^m(\C^m;L,\R)\cong H^{m-1}(L,\R)$ under
$H^{m-1}(L,\R)\ra H^{m-1}(\Si,R)$ in~\eq{cr7eq2}.
\label{cr7def2}
\end{dfn}

Here are some conditions for $Y(L)$ or $Z(L)$ to be zero.

\begin{prop} Let\/ $L$ be an AC SL\/ $m$-fold in $\C^m$ with
cone $C$ and rate $\la$, and let\/ $\Si=C\cap{\cal S}^{2m-1}$.
If\/ $\la<0$ or $b^1(L)=0$ then $Y(L)=0$. If\/ $\la<2-m$ or
$b^0(\Si)=1$ then~$Z(L)=0$.
\label{cr7prop1}
\end{prop}

\begin{proof} Let $C,\Si,K,T,\vp,\iota$ be as in Definition
\ref{cr7def1}. Suppose $Y(L)\ne 0$. Then there exists
$\ga\in H_1(\Si,\Z)$ with $Y(L)\cdot\ga\ne 0$. Choose
a closed 1-chain $\de$ in $\Si$ with $[\de]=\ga$. Let $r>T$.
Then $\de\t\{r\}$ is a closed 1-chain in $\Si\t(T,\iy)$, and
so $\vp\bigl(\de\t\{r\}\bigr)$, $\iota\bigl(\de\t\{r\}\bigr)$
are closed 1-chains in~$\C^m$.

Suppose $S$ is a 2-chain in $\C^m$ with $\pd S=\vp\bigl(\de\t
\{r\}\bigr)$. Then one can show using Definition \ref{cr7def2}
that $\int_S\om'=Y(L)\cdot\ga$. Also, if $T$ is a 2-chain in
$\C^m$ with $\pd T=\iota\bigl(\de\t\{r\}\bigr)$ then
$\int_T\om'=0$. This is because the answer is independent of
$T$, and we can choose $T$ inside the cone $C$, so that
$\om'\vert_T\equiv 0$ as $C$ is Lagrangian.

Therefore if $U$ is a 2-chain in $\C^m$ with $\pd U=
\vp\bigl(\de\t\{r\}\bigr)-\iota\bigl(\de\t\{r\}\bigr)$ then
$\int_U\om'=Y(L)\cdot\ga\ne 0$. Now $\bmd{\int_U\om'}\le\vol(U)$
as $\om'$ is a calibration. But by \eq{cr7eq1} we can choose
$U$ with $\vol(U)=O(r^\la)$ for large $r$. Thus
$\bmd{\int_U\om'}=O(r^\la)$ for large $r$, and also
$\int_U\om'\ne 0$ is independent of $r$. Together these
force $\la\ge 0$. Hence if $\la< 0$ then $Y(L)=0$. The
$\la<2-m$ case is similar, using $\Im\Om'$ instead of
$\om'$, and is left as an exercise.

If $b^1(L)=0$ then $H^1(L,\R)=0$, so $Y(L)=0$, as it lies in the
image of $H^1(L,\R)\ra H^1(\Si,\R)$ from Definition \ref{cr7def2}.
If $b^0(\Si)=1$ then $H_{m-1}(\Si,\R)=\an{[\Si]}$. Now $\Si$ is a
boundary in $L$, so the map $H_{m-1}(\Si,\R)\ra H_{m-1}(L,\R)$ is
zero, and the dual map $H^{m-1}(L,\R)\ra H^{m-1}(\Si,\R)$ also
zero. But $Z(L)$ lies in the image of this, so~$Z(L)=0$.
\end{proof}

\subsection{Lagrangian Neighbourhood Theorems}
\label{cr72}

Next we give versions of parts of \S\ref{cr4} for AC SL
$m$-folds rather than SL $m$-folds with conical singularities.
Here is an analogue of Theorem~\ref{cr4thm4}.

\begin{thm} Let\/ $C$ be an SL cone in $\C^m$ with isolated
singularity at\/ $0$, and set\/ $\Si=C\cap{\cal S}^{2m-1}$.
Define $\iota:\Si\t(0,\iy)\ra\C^m$ by $\iota(\si,r)=r\si$. Let\/
$\ze$, $U_{\sst C}\subset T^*\bigl(\Si\t(0,\iy)\bigr)$ and\/
$\Phi_{\sst C}:U_{\sst C}\ra\C^m$ be as in Theorem~\ref{cr4thm3}.

Suppose $L$ is an AC SL\/ $m$-fold in $\C^m$ with cone $C$
and rate $\la<2$. Then there exists a compact\/ $K\subset L$
and a diffeomorphism $\vp:\Si\t(T,\iy)\ra L\sm K$ for some
$T>0$ satisfying \eq{cr7eq1}, and a closed\/ $1$-form $\chi$ on
$\Si\t(T,\iy)$ written $\chi(\si,r)=\chi^1(\si,r)+\chi^2(\si,r)\d r$
for $\chi^1(\si,r)\in T_\si^*\Si$ and\/ $\chi^2(\si,r)\in\R$,
satisfying
\e
\begin{gathered}
\bmd{\chi(\si,r)}<\ze r,\quad \vp(\si,r)\equiv
\Phi_{\sst C}\bigl(\si,r,\chi^1(\si,r),\chi^2(\si,r)\bigr)\\
\text{and}\quad\bmd{\na^k\chi}=O(r^{\la-1-k})
\quad\text{as $r\ra\iy$ for $k=0,1,$}
\end{gathered}
\label{cr7eq3}
\e
computing $\na,\md{\,.\,}$ using the cone metric~$\iota^*(g')$.
\label{cr7thm1}
\end{thm}

\begin{proof} As $L$ is Asymptotically Conical with cone $C$ it
follows from \eq{cr7eq1} that near infinity in $\C^m$ we can write
$L$ as the image under $\Phi_{\sst C}$ of the graph of a 1-form
$\chi$ on $\Si\t(T,\iy)$ for large $T>0$. This just means that $L$
intersects the Lagrangian ball $\Phi_{\sst C}\bigl(T^*_{(\si,r)}
(\Si\t(0,\iy))\cap U_{\sst C}\bigr)$ in exactly one point for
$(\si,r)\in\Si\t(T,\iy)$, and we define $\chi$ such that this
point is~$\Phi_{\sst C}\bigl(\chi(\si,r)\bigr)$.

Now define $\vp:\Si\t(T,\iy)\ra L$ by $\vp(\si,r)=
\Phi_{\sst C}\bigl(\si,r,\chi^1(\si,r),\chi^2(\si,r)\bigr)$ and
$K=L\sm\Image\vp$. Then $K$ is compact and $\vp:\Si\t(T,\iy)
\ra L\sm K$ is a diffeomorphism. These $T,\vp,K$ are a valid
choice of $T,\vp,K$ for $L$ in Definition \ref{cr7def1}. In
particular, $\vp$ satisfies \eq{cr7eq1}. One can show this by
starting with $T',\vp',K'$ satisfying Definition \ref{cr7def1},
regarding $\vp$ as obtained from $\vp'$ by a kind of projection,
and showing that $\bmd{\na^k(\vp'-\vp)}=O(r^{\la-1-k})$ as
$r\ra\iy$ for~$k=0,1$.

As $\om'\vert_L\equiv 0$ and $\Phi_{\sst C}^*(\om')=\hat\om$ we
see that $\hat\om$ restricted to the graph of $\chi$ is zero.
By a well-known fact in symplectic geometry, this implies that
$\chi$ is {\it closed}. Equation \eq{cr7eq1} and the properties
of $\Phi_{\sst C}$ imply that $\md{\na^k\chi}=O(r^{\la-1-k})$ as
$r\ra\iy$ for $k=0,1$. As $\la<2$ this gives $\md{\chi}=o(r)$,
and so by making $T,K$ larger if necessary we can suppose that
$\bmd{\chi(\si,r)}<\ze r$ for $(\si,r)\in\Si\t(T,\iy)$. This
completes the proof.
\end{proof}

Here is the analogue of Theorem \ref{cr4thm5}. Its proof is a
straightforward modification of that of Theorem \ref{cr4thm5},
and we leave it as an exercise.

\begin{thm} Suppose $L$ is an AC SL\/ $m$-fold in $\C^m$ with cone
$C$. Let\/ $\Si,
\allowbreak
\iota,
\allowbreak
\ze,
\allowbreak
U_{\sst C},
\allowbreak
\Phi_{\sst C},K,T,\vp,\chi,\chi^1,\chi^2$ be as in Theorem \ref{cr7thm1}.
Then making $T,K$ larger if necessary, there exists an open tubular
neighbourhood\/ $U_{\sst L}\subset T^*L$ of the zero section $L$ in $T^*L$,
such that under $\d\vp:T^*\bigl(\Si\t(T,\iy)\bigr)\ra T^*L$ we have
\e
(\d\vp)^*(U_{\sst L})=\bigl\{(\si,r,\tau,u)\in
T^*\bigl(\Si\t(T,\iy)\bigr):\bmd{(\tau,u)}<\ze r\bigr\},
\label{cr7eq4}
\e
and there exists an embedding $\Phi_{\sst L}:U_{\sst L}\ra\C^m$ with\/
$\Phi_{\sst L}\vert_L=\id:L\ra L$ and\/ $\Phi_{\sst L}^*(\om')=\hat\om$,
where $\hat\om$ is the canonical symplectic structure on $T^*L$,
such that
\e
\Phi_{\sst L}\circ\d\vp(\si,r,\tau,u)\equiv
\Phi_{\sst C}\bigl(\si,r,\tau+\chi^1(\si,r),u+\chi^2(\si,r)\bigr)
\label{cr7eq5}
\e
for all\/ $(\si,r,\tau,u)\!\in\!T^*\bigl(\Si\!\t\!(T,\iy)\bigr)$ with\/
$\md{(\tau,u)}<\ze r$, computing $\md{\,.\,}$ using~$\iota^*(g')$.
\label{cr7thm2}
\end{thm}

We can decompose $\chi$ in Theorem \ref{cr7thm1}, in a similar
way to Lemma~\ref{cr4lem}.

\begin{prop} In the situation of Theorem \ref{cr7thm1} we have
$[\chi]=Y(L)$ in $H^1\bigl(\Si\t(T,\iy),\R\bigr)\cong H^1(\Si,\R)$,
where $Y(L)$ is as in Definition \ref{cr7def2}. Let\/ $\ga$ be the
unique $1$-form on $\Si$ with\/ $\d\ga=\d^*\ga=0$ and\/ $[\ga]=Y(L)
\in H^1(\Si,\R)$, which exists by Hodge theory. Then we may write
$\chi=\pi^*(\ga)+\d E$, where $\pi:\Si\t(T,\iy)\ra\Si$ is the
projection and $E\in C^\iy\bigl(\Si\t(T,\iy)\bigr)$, such that
\begin{itemize}
\setlength{\itemsep}{0pt}
\setlength{\parsep}{0pt}
\item[{\rm(a)}] If\/ $\la<0$ then $Y(L)=\ga=0$ and\/ $E$
is given by $E(\si,r)=-\int_r^\iy\chi^2(\si,s)\d s$ and
satisfies $\md{\na^kE}=O(r^{\la-k})$ for $k=0,1,2$
as~$r\ra\iy$.
\item[{\rm(b)}] If\/ $\la=0$ then $\md{E}=O\bigl(\md{\log r}\bigr)$
and\/ $\md{\na^kE}=O(r^{-k})$ for~$k=1,2$.
\item[{\rm(c)}] If\/ $\la>0$ then $\md{\na^kE}=O(r^{\la-k})$
for $k=0,1,2$ as~$r\ra\iy$.
\end{itemize}
Here we compute $\na,\md{\,.\,}$ using the cone metric
$\iota^*(g')$ on~$\Si\t(T,\iy)$.
\label{cr7prop2}
\end{prop}

\begin{proof} The proof that $[\chi]=Y(L)$ is similar to
Proposition \ref{cr7prop1}, and we leave it as an exercise.
Let $\ga$ be as in the proposition. Then $\pi^*(\ga)$ is a
closed 1-form on $\Si\t(T,\iy)$ with $\bigl[\pi^*(\ga)\bigr]
=Y(L)=[\chi]\in H^1\bigl(\Si\t(T,\iy),\R\bigr)$. Thus
$\chi-\pi^*(\ga)$ is an {\it exact\/} 1-form, and we may
write $\chi-\pi^*(\ga)=\d E$ for some $E\in C^\iy\bigl(
\Si\t(T,\iy)\bigr)$, unique up to addition of constants.

For part (a), if $\la<0$ then $Y(L)=0$ by Proposition
\ref{cr7prop1}, so $\ga=0$. By \eq{cr7eq3} we see that
$E'(\si,r)=-\int_r^\iy\chi^2(\si,s)\d s$ is well-defined.
The $\d r$ component in $\d E'$ is $\chi^2$, so that
$\chi-\d E'$ is a closed 1-form on $\Si\t(T,\iy)$ with no
$\d r$ component, and is therefore independent of $r$. But
\eq{cr7eq3} implies that $\chi-\d E'=O(r^{\la-1})$ in the
cone metric on $\Si\t(T,\iy)$, so $\chi-\d E'=O(r^\la)$ in
the cylinder metric, and taking the limit $r\ra\iy$ gives
$\chi-\d E'=0$ as $\la<0$. Thus we may take $E=E'$, and
\eq{cr7eq3} then yields $\md{\na^kE}=O(r^{\la-k})$ as
$r\ra\iy$ for~$k=0,1,2$.

For parts (b) and (c), using $\bmd{\na^k\pi^*(\ga)}=O(r^{-1-k})$,
equation \eq{cr7eq3} and $\chi=\pi^*(\ga)+\d E$, we find that if
$\la\ge 0$ then $\md{\na^kE}=O(r^{\la-k})$ for $k=1,2$. But
\e
E(\si,r)=E(\si,T+1)+\int_{T+1}^r\frac{\d E}{\d r}(\si,s)\d s
\label{cr7eq6}
\e
for $r\ge T+1$, and $\md{\frac{\d E}{\d r}(\si,s)}\le \md{\na
E(\si,s)}=O(s^{\la-1})$. Substituting this into \eq{cr7eq6} gives
$\md{E}=O\bigl(\md{\log r}\bigr)$ for $\la=0$ and $\md{E}=O(r^\la)$
for $\la>0$, which completes the proof.
\end{proof}

\subsection{The asymptotic behaviour of $L$ at infinity}
\label{cr73}

Finally we give analogues of the material of \S\ref{cr5} for
AC SL $m$-folds. Here is the analogue of Theorem \ref{cr5thm1}.
Stephen Marshall also has an independent proof.

\begin{thm} In the situation of Theorem \ref{cr7thm1} and
Proposition \ref{cr7prop2} we have
\e
\begin{gathered}
\bmd{\na^k(\vp-\iota)}=O(r^{\la-1-k}),\quad
\bmd{\na^k\chi}=O(r^{\la-1-k})\quad\text{for all $k\ge 0$}\\
\text{and}\quad\bmd{\na^kE}=O(r^{\la-k})
\quad\text{for all\/ $k\ge 1$ as $r\ra\iy$.}
\end{gathered}
\label{cr7eq7}
\e
Here $\na,\md{\,.\,}$ are computed using the cone metric
$\iota^*(g')$ on~$\Si\t(T,\iy)$.
\label{cr7thm3}
\end{thm}

\begin{proof} We modify the proof of Theorem \ref{cr5thm1}.
Let $\al$ be a smooth 1-form on $\Si\t(T,\iy)$ with
$\md{\al(\si,r)}<\ze r$, written $\al(\si,r)=\al^1(\si,r)+
\al^2(\si,r)\d r$ for $\al^1(\si,r)\in T_\si^*\Si$ and
$\al^2(\si,r)\in\R$. Define a map $\Th_\al:\Si\t(T,\iy)\ra
\C^m$ by $\Th_\al(\si,r)=\Phi_{\sst C}\bigl(\si,r,\al^1(\si,r),
\al^2(\si,r)\bigr)$. Define a smooth real function $F(\al)$ on
$\Si\t(T,\iy)$ by $F(\al)\,\d V=\Th_\al^*(\Im\Om')$, where
$\d V$ is the volume form of $\iota^*(g')$ on $\Si\t(T,\iy)$.
As in \eq{cr5eq3}--\eq{cr5eq4}, define
\begin{gather}
\begin{split}
Q:\bigl\{(\si,r,y,z):\,&(\si,r)\in\Si\t(T,\iy),\quad
y\in T^*_{(\si,r)}\bigl(\Si\t(T,\iy)\bigr),\\
&\md{y}<\ze r,\quad z\in\ot^2T^*_{(\si,r)}\bigl(\Si\t(T,\iy)
\bigr)\bigr\}\ra\R
\end{split}
\label{cr7eq8}\\
\text{by}\quad
Q\bigl(\si,r,\al(\si,r),\na\al(\si,r)\bigr)=
\bigl(\d^*\al+F(\al)\bigr)\bigl[(\si,r)\bigr]
\label{cr7eq9}
\end{gather}
for all 1-forms $\al$ on $\Si\t(T,\iy)$ with $\md{\al(\si,r)}<\ze r$
when $(\si,r)\in\Si\t(T,\iy)$. Then $F,Q$ are well-defined, and
analogous to $F_i,Q_i$ in the proof of Theorem~\ref{cr5thm1}.

As $\d^*\pi^*(\ga)=0$ on $\Si\t(T,\iy)$, the proofs of
\eq{cr5eq8} and \eq{cr5eq9} give
\begin{gather}
Q(\si,r,y,z)=O(r^{-2}\ms{y}+\ms{z})
\quad\text{when $\md{y}=O(r)$ and $\md{z}=O(1)$, and}
\label{cr7eq10}\\
\De E(\si,r)-Q\bigl(\si,r,\pi^*(\ga)(\si,r)+\d E(\si,r),
\na\pi^*(\ga)+\na^2E(\si,r)\bigr)=0.
\label{cr7eq11}
\end{gather}
Adapting the rest of the proof of Theorem \ref{cr5thm1} now proves
the result. Note that the $\pi^*(\ga)$ terms in \eq{cr7eq11} are
zero when $\la<0$, and when $\la\ge 0$ they can be absorbed into
the other estimates, and do not cause a problem.
\end{proof}

From \eq{cr7eq7} we deduce that on $\Si\t(T,\iy)$ we have
\e
\bmd{\na^k\bigl(\vp^*(g')-\iota^*(g')\bigr)}=O(r^{\la-2-k})
\quad\text{as $r\ra\iy$ for all $k\ge 0$,}
\label{cr7eq12}
\e
computing $\na,\md{\,.\,}$ using $\iota^*(g')$. This is an
analogue of equation \eq{cr2eq1}. As in Theorem \ref{cr5thm2},
equation \eq{cr7eq12} implies that $(L,g')$ is an {\it
Asymptotically Conical Riemannian manifold\/} with cone $C$ and
rate $\la-2<0$, in a sense analogous to Definition \ref{cr2def1}.
Therefore we can develop a theory of {\it analysis on AC SL\/
$m$-folds}, similar to \S\ref{cr2}. Here is the analogue of
Definitions \ref{cr2def4} and~\ref{cr2def5}.

\begin{dfn} Let $L$ be an AC SL $m$-fold in $\C^m$, as in
Definition \ref{cr7def1}. Define the {\it radius function}
$\rho:L\ra[1,\iy)$ by $\rho(x)=(1+\ms{x})^{1/2}$. For
$\be\in\R$ and $k\ge 0$ define $C^k_{\smash\be}(L)$ to be the
space of continuous functions $f$ on $L$ with $k$ continuous
derivatives, such that $\md{\rho^{-\be+j}\na^jf}$ is bounded
on $L$ for $j=0,\dots,k$. Define the norm $\nm{\,.\,}_{\smash{
C^k_\be}}$ on $C^k_{\smash\be}(L)$ by $\nm{f}_{\smash{C^k_\be}}
=\sum_{j=0}^k\sup_L\md{\rho^{-\be+j}\na^jf}$. Then
$C^k_{\smash\be}(L)$ is a Banach space. Define~$C^\iy_{
\smash\be}(L)=\bigcap_{k\ge 0}C^k_{\smash\be}(L)$.

For $p\ge 1$, $\be\in\R$ and $k\ge 0$ define the {\it weighted
Sobolev space} $L^p_{k,\be}(L)$ to be the set of functions $f$
on $L$ that are locally integrable and $k$ times weakly
differentiable, and for which the norm $\snm{f}p{k,\be}=
\bigl(\sum_{j=0}^k\int_L\md{\rho^{-\be+j}\na^jf}^p\rho^{-m}\,
\d V_g\bigr)^{\smash{1/p}}$ is finite. Then $L^p_{k,\be}(L)$
is a Banach space, and $L^2_{k,\be}(L)$ a Hilbert space. 
\label{cr7def3}
\end{dfn}

We can now develop the theory of \S\ref{cr22}--\S\ref{cr23} for
these spaces. This is done in detail by Marshall \cite[\S 4]{Mars}.
The basic references \cite{Lock,LoMc} apply to Asymptotically
Conical Riemannian manifolds just as for Riemannian manifolds
with conical singularities. Theorem \ref{cr2thm1} holds for
$C^k_{\smash\be}(L),L^p_{k,\be}(L)$ except that the directions
of the inequalities ${\bs\be}\ge{\bs\ga}$, ${\bs\be}>{\bs\ga}$
must be reversed. As in \eq{cr2eq8}, let $\De=\d^*\d$ be the
Laplacian on $L$, and for $p>1$, $k\ge 2$ and $\be\in\R$ write
$\De^p_{k,\be}$ for the map
\e
\De^p_{k,\be}=\De:L^p_{k,\be}(L)\ra L^p_{k-2,\be-2}(L).
\label{cr7eq13}
\e
Then we can prove the following condensation of the analogue
of~\S\ref{cr23}:

\begin{thm} Let\/ $L$ be an AC SL\/ $m$-fold in $\C^m$, with
cone $C$, and set\/ $\Si=C\cap{\cal S}^{2m-1}$. Let\/ $\D_\sSi$
and\/ $N_\sSi$ be as in Definition \ref{cr2def3}. Let\/ $p>1$,
$k\ge 2$ and\/ $\be\in\R$, and define $q>1$ by $\frac{1}{p}+
\frac{1}{q}=1$. Let\/ $\De^p_{k,\be}$ be as in \eq{cr7eq13}. Then
\begin{itemize}
\setlength{\itemsep}{0pt}
\setlength{\parsep}{0pt}
\item[{\rm(a)}] $\De^p_{k,\be}$ is Fredholm if and only
if\/~$\be\notin\D_\sSi$.
\item[{\rm(b)}] If\/ $\De^p_{k,\be}$ is Fredholm 
then~$\ind\bigl(\De^p_{k,\be}\bigr)=N_\sSi(\be)$.
\item[{\rm(c)}] $\Ker\bigl(\De^p_{k,\be}\bigr)$ is a
finite-dimensional subspace of\/ $C^\iy_{\smash\be}(L)$, and
independent of\/ $k$. If\/ $\be\notin\D_\sSi$ then $\Ker\bigl(
\De^p_{k,\be}\bigr)$ is independent of\/ $p$, and depends only
on the connected component of\/ $\R\sm\D_\sSi$ containing $\be$.
If\/ $\be<0$ then~$\Ker\bigl(\De^p_{k,\be}\bigr)\!=\!0$.
\item[{\rm(d)}] Suppose $\De^p_{k,\be}$ is Fredholm. Then
$u\in L^p_{k-2,\be-2}(L)$ lies in the image of $\De^p_{k,\be}$
if and only if\/ $\int_Luv\,\d V=0$ for all\/~$v\in\Ker\bigl(
\De^q_{2,-\be+2-m}\bigr)$.
\end{itemize}
\label{cr7thm4}
\end{thm}

We study the vector space $V$ of {\it bounded harmonic functions\/}
on $L$. A similar result is proved by Marshall~\cite[\S 5.1.3]{Mars}.

\begin{thm} Suppose $L$ is an AC SL\/ $m$-fold in $\C^m$, with cone
$C$. Let\/ $\Si,T$ and\/ $\vp$ be as in Theorem \ref{cr7thm1}. Let\/
$l=b^0(\Si)$, and\/ $\Si^1,\ldots,\Si^l$ be the connected components
of\/ $\Si$. Let\/ $V$ be the vector space of bounded harmonic functions
on $L$. Then $\dim V=l$, and for each\/ ${\bf c}=(c^1,\ldots,c^l)\in\R^l$
there exists a unique $v^{\bf c}\in V$ such that for all\/ $j=1,\ldots,l$,
$k\ge 0$ and\/ $\be\in(2-m,0)$ we have
\e
\na^k\bigl(\vp^*(v^{\bf c})-c^j\,\bigr)=O\bigl(\md{{\bf c}}
r^{\be-k}\bigr)\quad\text{on $\Si^j\t(T,\iy)$ as $r\ra\iy$.}
\label{cr7eq14}
\e
Note also that\/ $V=\{v^{\bf c}:{\bf c}\in\R^l\}$
and\/~$v^{(1,\ldots,1)}\equiv 1$.
\label{cr7thm5}
\end{thm}

\begin{proof} Let $\D_\sSi$ and $N_\sSi$ be as in Definition
\ref{cr2def3}, and choose $p>1$ and $0<\ga<\min\bigl(\D_\sSi
\cap(0,\iy)\bigr)$. Let $v\in V$. Then $v\in L^p_{0,\ga}(L)$
as $v$ is bounded and $\ga>0$. Also $v$ is smooth, as it is
harmonic. By the analogue of Theorem \ref{cr2thm2} for AC
SL $m$-folds we find that $v\in L^p_{k,\ga}(L)$ for all
$k\ge 0$. Fix $k\ge 2$. As $\De v=0$ we have $v\in\Ker
(\De^p_{k,\ga})$, so that~$V\subseteq\Ker(\De^p_{k,\ga})$.

Part (a) of Theorem \ref{cr2thm4} shows that $\De^p_{k,\ga}$
is Fredholm as $\ga\notin\D_\sSi$, so
\e
\ind\bigl(\De^p_{k,\ga}\bigr)=N_\sSi(\ga)=
N_\sSi(0)=b^0(\Si)=l,
\label{cr7eq15}
\e
by part (b) of Theorem \ref{cr2thm4}. Here $N_\sSi(\ga)=
N_\sSi(0)$ as $\D_\sSi\cap(0,\ga]=\emptyset$ and $N_\sSi$ is
upper semicontinuous and discontinuous exactly on $\D_\sSi$,
and $N_\sSi(0)=m_\sSi(0)=b^0(\Si)$ is the multiplicity of the
eigenvalue 0 of $\De_\sSi$. Now $\De^p_{k,\ga}$ is surjective
by parts (c) and (d) of Theorem \ref{cr7thm4}. Thus $\dim\Ker
(\De^p_{k,\ga})=l$ by \eq{cr7eq15} and $V\subseteq\Ker(
\De^p_{k,\ga})$, which proves that~$\dim V\le l$.

Let ${\bf c}=(c^1,\ldots,c^l)\in\R^l$ and $\be\in(2-m,0)$,
and choose a smooth function $\hat v^{\bf c}$ on $L$ with
$\hat v^{\bf c}\equiv c^j$ on $\vp\bigl(\Si^j\t(T+1,\iy)\bigr)$ for
$j=1,\ldots,l$. Clearly this is possible. Then $\De\hat v^{\bf c}$
is smooth and compactly-supported on $L$, so $\De\hat v^{\bf c}\in
L^p_{k-2,\be-2}(L)$. Now as $\be\in(2-m,0)$, $k\ge 2$ and $p>1$
parts (a), (c) and (d) of Theorem \ref{cr7thm4} show that
$\De^p_{k,\be}:L^p_{k,\be}(L)\ra L^p_{k-2,\be-2}(L)$ is an
isomorphism.

Hence there exists a unique $\ti v^{\bf c}\in L^p_{k,\be}(L)$ with
$\De\ti v^{\bf c}=\De\hat v^{\bf c}$. This $\ti v^{\bf c}$ is
independent of $k\ge 2$ and $\be$, so $\ti v^{\bf c}\in C^\iy_\be(L)$
by the analogue of Theorem \ref{cr2thm1}, for all $\be\in(2-m,0)$.
Define $v^{\bf c}=\hat v^{\bf c}-\ti v^{\bf c}$. Then $\De v^{\bf c}=
\De\hat v^{\bf c}-\De\hat v^{\bf c}=0$, and \eq{cr7eq14} holds as
$\hat v^{\bf c}\equiv c^j$ on $\vp\bigl(\Si^j\t(T+1,\iy)\bigr)$ and
$\ti v^{\bf c}\in C^\iy_\be(L)$. Thus $v^{\bf c}$ is harmonic and
bounded, so~$v^{\bf c}\in V$.

All possible such functions $v^{\bf c}$ for ${\bf c}\in\R^l$
generate a vector subspace of $V$, of dimension at least $l$.
But $\dim V\le l$ from above. Hence each $v^{\bf c}$ is unique,
and $V=\{v^{\bf c}:{\bf c}\in\R^l\}$, and $\dim V=l$, as we have
to prove. Finally, $1\in V$, and clearly~$v^{(1,\ldots,1)}=1$.
\end{proof}

We finish by proving a version of Theorem \ref{cr5thm5} for AC
SL $m$-folds, improving the rate of convergence~$\la$.

\begin{thm} Let\/ $L$ be an AC SL\/ $m$-fold in $\C^m$ with
cone $C$ and rate $\la$. Set\/ $\Si=C\cap{\cal S}^{2m-1}$,
and let\/ $\D_\sSi,N_\sSi$ be as in Definition \ref{cr2def3}.
Let\/ $\iota,T,\vp,\chi$ be as in Theorem \ref{cr7thm1},
and\/ $Y(L),\ga,E$ as in Proposition \ref{cr7prop2}. Then
\begin{itemize}
\setlength{\itemsep}{0pt}
\setlength{\parsep}{0pt}
\item[{\rm(a)}] Suppose $\la,\la'$ lie in the same
connected component of\/ $\R\sm\D_\sSi$. Then
\e
\begin{gathered}
\bmd{\na^k(\vp-\iota)}=O(r^{\la'-1-k}),\quad
\bmd{\na^k\chi}=O(r^{\la'-1-k})\quad\text{and}\\
\bmd{\na^kE}=O(r^{\la'-k})
\quad\text{as $r\ra\iy$ for all\/ $k\ge 0$.}
\end{gathered}
\label{cr7eq16}
\e
Hence $L$ is an AC SL\/ $m$-fold with rate $\la'$. In particular,
if\/ $\la\in(2-m,0)$ then $L$ is an AC SL\/ $m$-fold with rate
$\la'$ for all\/~$\la'\in(2-m,0)$.
\item[{\rm(b)}] Suppose $0\le\la<\min\bigl(\D_\sSi\cap(0,\iy)\bigr)$.
Then adding a constant to $E$ if necessary, for all\/ $\la'\in\bigl(
\max(-2,2-m),0\bigr)$ we have
\e
\bmd{\na^kE}=O(r^{\la'-k})\quad\text{as $r\ra\iy$ for all\/ $k\ge 0$.}
\label{cr7eq17}
\e
Thus if\/ $Y(L)=0=\ga$ then $L$ is an AC SL\/ $m$-fold with rate
$\la'$, and if\/ $Y(L)\ne 0\ne\ga$ then $L$ is an AC SL\/ $m$-fold
with rate~$0$.
\end{itemize}
\label{cr7thm6}
\end{thm}

\begin{proof} Use the notation of the rest of this section. Define
a smooth function $\ti E:L\ra\R$ by $\ti E\bigl(\vp(\si,r)\bigr)=
E(\si,r)$ on $L\sm K$, and extend $\ti E$ smoothly over $K$. Then
$\ti E\in C^\iy_\la(L)$ by \eq{cr7eq7}. Write $g_{\sst L}$ for the
metric $g'\vert_L$ on $L$, and $g_{\sst C}$ for the cone metric
$\vp_*\bigl(\iota^*(g')\bigr)$ on $L\sm K\cong\Si\t(T,\iy)$. Let
$\d^*_{\sst L},\d^*_{\sst C}$ be the $\d^*$ operators w.r.t.\
$g_{\sst L},g_{\sst C}$ on $L,L\sm K$. Let $\De_{\sst L}=
\d^*_{\sst L}\d$ be the Laplacian on $L$. Let $Q$ be as in
\eq{cr7eq8}--\eq{cr7eq9}, and $Q'$ the push-forward of $Q$ to
$L\sm K$ under $\vp$. Let $\ti\ga$ be the push-forward of
$\pi^*(\ga)$ to $L\sm K$ under~$\vp$.

Following the proof of \eq{cr5eq21}, equation \eq{cr7eq11}
implies that for $x\in L\sm K$
\e
\De_{\sst L}\ti E=Q'\bigl(x,\ti\ga+\d\ti E(x),\na\ti\ga+
\na^2\ti E(x)\bigr)+(\d^*_{\sst L}-\d^*_{\sst C})\d\ti E[x].
\label{cr7eq18}
\e
From \eq{cr7eq7}, \eq{cr7eq10} and similar
estimates on the derivatives of $Q$ we deduce that
\e
\bmd{\na^k(\De_{\sst L}\ti E)}=
\begin{cases} O(\rho^{-4-k})+O(\rho^{2\la-4-k}), & \ga\ne 0, \\
O(\rho^{2\la-4-k}), & \ga=0, \end{cases}
\label{cr7eq19}
\e
as $\md{\na^k\ti\ga}=O(\rho^{-1-k})$. Since $\la\ge 0$ when
$\ga\ne 0$, this gives $\De_{\sst L}\ti E\in C^\iy_{2\la-4}(L)$,
which is the analogue of Lemma \ref{cr5lem1}. Here is the analogue
of Lemma~\ref{cr5lem2}.

\begin{lem} Let\/ $p>1$, $k\ge 2$ and\/ $\la,\la'$ lie in the
same connected component of\/ $\R\sm\D_\sSi$. If\/ $\ti E\in
L^p_{k,\la}(L)$ and\/ $\De_{\sst L}\ti E\in L^p_{k-2,\la'-2}(L)$,
then~$\ti E\in L^p_{k,\la'}(L)$.
\label{cr7lem}
\end{lem}

\begin{proof} This is trivial for $\la'\ge\la$, so suppose
$\la'<\la$. As $\la\in\R\sm\D_\sSi$ part (a) of Theorem
\ref{cr7thm4} shows that $\De^p_{k,\la}$ is Fredholm, and
thus part (d) that $\De_{\sst L}\ti E$ is $L^2$-orthogonal
to $\Ker(\De^q_{2,-\la+2-m})$. But $\la,\la'$ lie in the
same connected component of $\R\sm\D_\sSi$, and $\D_\sSi$
is preserved by the involution $\be\mapsto-\be+2-m$ by
\eq{cr2eq4}, so $-\la+2-m,-\la'+2-m$ lie in the same
connected component of~$\R\sm\D_\sSi$.

Hence $\Ker(\De^q_{2,-\la+2-m})=\Ker(\De^q_{2,-\la'+2-m})$ by
part (c), so $\De_{\sst L}\ti E$ lies in $L^p_{k-2,\la'-2}(L)$ by
assumption and is $L^2$-orthogonal to $\Ker(\De^q_{2,-\la'+2-m})$.
Thus $\De_{\sst L}\ti E$ lies in the image of $\De^p_{k,\la'}$
by part (d) of Theorem~\ref{cr7thm4}.

Therefore $\De_{\sst L}\ti E=\De_{\sst L}E'$ for some $E'\in
L^p_{k,\la'}(L)$. Then $E'\in L^p_{k,\la}(L)$ as $\la'<\la$,
so $\ti E-E'\in\Ker(\De^p_{k,\la})$. But $\Ker(\De^p_{k,\la})
=\Ker(\De^p_{k,\la'})$ by part (c), so both $E'$ and $\ti E-E'$
lie in $L^p_{k,\la'}(L)$, and~$\ti E\in L^p_{k,\la'}(L)$.
\end{proof}

We can now use the method of Theorem \ref{cr5thm5} to decrease
the rate $\la$ by an inductive process. Applying Lemma \ref{cr7lem}
repeatedly $j$ times as in the proof of Theorem \ref{cr5thm5} shows
that \eq{cr7eq16} holds for all $\la'$ in the same connected component
of $\R\sm\D_\sSi$ as $\la$ with $\la'-2>2^j(\la-2)$. But
$2^j(\la-2)\ra-\iy$ as $j\ra\iy$ since $\la<2$, so this proves
part (a) of Theorem~\ref{cr7thm6}.

Now suppose that $L$ has rate $\la$ with $0\le\la<\min\bigl(
\D_\sSi\cap(0,\iy)\bigr)$, as in part (b). Then (a) implies
that $L$ has rate $\mu$ for all $\mu$ with $0<\mu<\min\bigl(
\D_\sSi\cap(0,\iy)\bigr)$. Thus \eq{cr7eq19} gives
$\De_{\sst L}\ti E\in C^\iy_{2\mu-4}(L)$. Therefore
$\De_{\sst L}\ti E\in C^\iy_{\la'-2}(L)$ for all $\la'>-2$,
and $\De_{\sst L}\ti E\in L^p_{k-2,\la'-2}(L)$ for all
$p>1,k\ge 2$ and~$\la'>-2$.

Let $p>1,k\ge 2$ and $\la'\in\bigl(\max(-2,2-m),0\bigr)$. Then
$\De_{\sst L}\ti E\in L^p_{k-2,\la'-2}(L)$, and $\la'>2-m$
implies $-\la'+2-m<0$, so that $\Ker(\De^q_{2,-\la'+2-m})=0$
by part (c) of Theorem \ref{cr7thm4}. Thus $\De_{\sst L}\ti E$
is trivially $L^2$-orthogonal to $\Ker(\De^q_{2,-\la'+2-m})$.
Also $\De^p_{k,\la'}$ is Fredholm by part (a) of Theorem
\ref{cr7thm4}, as $\la'\in(2-m,0)$ and $\D_\sSi\cap(2-m,0)=
\emptyset$. Therefore part (d) of Theorem \ref{cr7thm4} shows
that $\De_{\sst L}\ti E$ lies in the image of $\De^p_{k,\la'}$,
so $\De_{\sst L}\ti E=\De_{\sst L}E'$ for some~$E'\in
L^p_{k,\la'}(L)$.

As $\la'<\la$ we have $E'\in L^p_{k,\la}(L)$, and so
$\ti E-E'\in\Ker(\De^p_{k,\la})$. Increasing $\la$ if $\la=0$
we may take $0<\la<\min\bigl(\D_\sSi\cap(0,\iy)\bigr)$, so
that $\la\notin\D_\sSi$ and $N_\sSi(\la)=m_\sSi(0)=b^0(\Si)$.
Parts (a), (c) and (d) then show that $\De^p_{k,\la}$ is
Fredholm and surjective with~$\ind(\De^p_{k,\la})=N_\sSi(\la)$.

Hence $\dim\Ker(\De^p_{k,\la})=b^0(\Si)$. If $b^0(\Si)=1$ then
$\Ker(\De^p_{k,\la})=\an{1}$, the constant functions on $L$. More
generally, if $b^0(\Si)=l$ then $L$ has $l$ ends at infinity, and
elements of $\Ker(\De^p_{k,\la})$ are harmonic functions on $L$
which are asymptotic to $O(\rho^{2-m})$ at infinity to a constant
$c_i$ for $i=1,\ldots,l$ on each of the $l$ ends. The values of
$c_1,\ldots,c_l$ parametrize~$\Ker(\De^p_{k,\la})\cong\R^l$.

Now the function $E\in C^\iy\bigl(\Si\t(T,\iy)\bigr)$ was defined
in Proposition \ref{cr7prop2} to satisfy $\d E=\chi$. Thus, if
$b^0(\Si)=l$ then $E$ is unique up to the addition of a constant
on each of the $l$ components of $\Si\t(T,\iy)$. By choosing
these constants appropriately we can set to zero the constants
$c_1,\ldots,c_l$ that $\ti E-E'$ is asymptotic to on the $l$
ends of $L$. Then $\ti E-E'=0$, as $c_1,\ldots,c_l$
parametrize~$\Ker(\De^p_{k,\la})$.

Thus adding a constant to $E$ if necessary we have $\ti E=E'$,
so $\ti E\in\smash{L^p_{k,\la'}}(L)$. As this holds for all $p>1$ and
$k\ge 2$ we have $\ti E\in C^\iy_{\la'}(L)$ by the Asymptotically
Conical version of Theorem \ref{cr2thm1}. But $\vp^*(\ti E)=E$ on
$\Si\t(T,\iy)$, so this implies \eq{cr7eq17}. The last part is
immediate. This concludes the proof of Theorem~\ref{cr7thm6}.
\end{proof}


\begin{thebibliography}{99}

\bibitem{AdSi} D. Adams and L. Simon, {\it Rates of Asymptotic
Convergence Near Isolated Singularities of Geometric Extrema},
Indiana Univ. Math. J. 37 (1988), 225--254.

\bibitem{AlAl} W.K. Allard and F.J. Almgren, {\it On the radial
behaviour of minimal surfaces and the uniqueness of their tangent
cones}, Annals of Math. 113 (1981), 215--265.

\bibitem{Bart} R. Bartnik, {\it The Mass of an Asymptotically
Flat Manifold}, Comm. Pure Appl. Math. 39 (1986), 661--693.

\bibitem{Bred} G.E. Bredon, {\it Topology and Geometry}, Graduate
Texts in Mathematics 139, Springer-Verlag, Berlin, 1993.

\bibitem{Fede} H. Federer, {\it Geometric Measure Theory},
Grundlehren der math. Wiss. 153, Springer-Verlag, Berlin, 1969.

\bibitem{GiTr} D. Gilbarg and N.S. Trudinger, {\it Elliptic Partial
Differential Equations of Second Order}, Classics in Mathematics,
Springer-Verlag, Berlin, 2001.

\bibitem{HaLa} R. Harvey and H.B. Lawson, {\it Calibrated geometries},
Acta Mathematica 148 (1982), 47--157.

\bibitem{Ivan} A.V. Ivanov, {\it A priori estimates for solutions
of nonlinear second-order elliptic equations}, J. Soviet Math. 10 (1978),
217--240.

\bibitem{Joyc1} D.D. Joyce, {\it On counting special Lagrangian homology
$3$-spheres}, pages 125--151 in {\it Topology and Geometry: Commemorating
SISTAG}, editors A.J. Berrick, M.C. Leung and X.W. Xu, Contemporary
Mathematics 314, A.M.S., Providence, RI, 2002. hep-th/9907013.

\bibitem{Joyc2} D.D. Joyce, {\it Lectures on Calabi--Yau and special
Lagrangian geometry}, math.DG/0108088, 2001. Published, with extra
material, as Part I of M. Gross, D. Huybrechts and D. Joyce,
{\it Calabi--Yau Manifolds and Related Geometries}, Universitext
series, Springer, Berlin, 2003.

\bibitem{Joyc3} D.D. Joyce, {\it Special Lagrangian submanifolds
with conical singularities. II. Moduli spaces}, math.DG/0211295,
version 3, 2003.

\bibitem{Joyc4} D.D. Joyce, {\it Special Lagrangian submanifolds with
isolated conical singularities. III. Desingularization, the unobstructed
case}, math.DG/0302355, version 2, 2003.

\bibitem{Joyc5} D.D. Joyce, {\it Special Lagrangian submanifolds with
isolated conical singularities. IV. Desingularization, obstructions
and families}, math.DG/0302356, version 2, 2003.

\bibitem{Joyc6} D.D. Joyce, {\it Special Lagrangian submanifolds
with isolated conical singularities. V. Survey and applications},
math.DG/0303272, 2003.

\bibitem{Laws} H.B. Lawson, {\it Lectures on Minimal Submanifolds,
Volume 1}, Publish or Perish, Berkeley, CA, 1980.

\bibitem{Lock} R. Lockhart, {\it Fredholm, Hodge and Liouville Theorems
on noncompact manifolds}, Trans. A.M.S. 301 (1987), 1--35.

\bibitem{LoMc} R.B. Lockhart and R.C. McOwen, {\it Elliptic
Differential Operators on Noncompact Manifolds}, Annali della Scuola
normal superiore di Pisa, Classe di scienze 12 (1987), 409--447.

\bibitem{Mars} S.P. Marshall, {\it Deformations of special Lagrangian
submanifolds}, Oxford D.Phil. thesis, 2002.

\bibitem{McSa} D. McDuff and D. Salamon, {\it Introduction to
symplectic topology}, second edition, OUP, Oxford, 1998.

\bibitem{McLe} R.C. McLean, {\it Deformations of calibrated submanifolds},
Communications in Analysis and Geometry 6 (1998), 705--747.

\bibitem{Morg} F. Morgan, {\it Geometric Measure Theory, A Beginner's
Guide}, Academic Press, San Diego, 1995.

\bibitem{Morr} C.B. Morrey, {\it Multiple Integrals in the Calculus
of Variations}, Grundlehren der math. Wiss. 130, Springer-Verlag,
Berlin, 1966.

\bibitem{Simo1} L. Simon, {\it Asymptotics for a class of non-linear
evolution equations, with applications to geometric problems},
Annals of Math. 118 (1983), 525--571.

\bibitem{Simo2} L. Simon, {\it Isolated singularities of extrema
of geometric variational problems}, pages 206--277 in E. Giusti,
editor, {\it Harmonic Mappings and Minimal Immersions}, Springer
Lecture Notes in Math. 1161, Springer, Berlin, 1985.

\bibitem{SYZ} A. Strominger, S.-T. Yau, and E. Zaslow, {\it Mirror 
symmetry is T-duality}, Nuclear Physics B479 (1996), 243--259.
hep-th/9606040.

\bibitem{Wein} A. Weinstein, {\it Symplectic manifolds and their
Lagrangian submanifolds}, Advances in Math. 6 (1971), 329--346.

\end{thebibliography}
\end{document}